\documentclass[12pt,a4paper,twoside,titlepage]{report}
\usepackage{latexsym,amsmath,amssymb,amsfonts,amsthm}
\usepackage[mathscr]{eucal}
\usepackage{makeidx}

\newtheorem{definition}{Definition}[chapter]

\newtheorem{theorem}{Theorem}[chapter]
\def\nabla{\bigtriangledown}
\makeindex

\begin{document}

\title{ Nonlinear Connections and  Clifford Structures}
\author{Sergiu I.\ Vacaru \thanks{%
e--mail: vacaru@fisica.ist.utl.pt, sergiu$_{-}$vacaru@yahoo.com,\ $%
^\diamondsuit$ nadejda$_{-}$vicol@yahoo.com }, and Nadejda A. Vicol $%
^\diamondsuit$ \\
\\
{\small * Centro Multidisciplinar de Astrofisica - CENTRA,}\\
{\small Departamento de Fisica, Instituto Superior Tecnico,}\\
{\small Av. Rovisco Pais 1, Lisboa, 1049-001, Portugal }\\
{} \\
{\small $\diamondsuit$ Faculty of Mathematics and Informatics, gr. 33 MI,}\\
{\small State University of Moldova, Mateevici str. 60,}\\
{\small Chi\c sin\v au MD2009, Republic of Moldova} }
\date{May 17, 2002}
\maketitle
\newpage
--
\newpage

\begin{abstract}
We present an introduction to the geometry of higher order vector
and co--vector bundles (including higher order generalizations of
the Finsler geometry and Kaluza--Klein gravity)
\cite{mhss,m1,m2,vstr2,vbook} and review the basic results on
Clifford and spinor structures on spaces with generic local
anisotropy modeled by higher order nonlinear connections
\cite{vhsp,vp}. Geometric applications in locally anisotropic
gravity and matter field interactions are considered \cite{vd}.
This article contains the results outlined by authors and P.
Stavrinos in theirs lectures \cite{vsv,vvcook}.



\end{abstract}

\newpage
--
\newpage

\tableofcontents

\vskip0.5cm

\section{Introduction}

\bigskip The spinors studied by mathematicians and physicists are connected
with the general theory of Clifford spaces introduced in 1876
\cite{clifford}. The theory of spinors and Clifford algebras play
a major role in
contemporary physics and mathematics. The  spinors were discovered by \`{E}%
lie Cartan in 1913 in  mathematical form in his researches on representation
group theory \cite{car38};  he showed that spinors furnish a linear
representation of the groups of rotations of a space of arbitrary
dimensions. Physicists Pauli \cite{pauli} and Dirac \cite{dirac} (in 1927,
respectively, for the three--dimensional and four--dimensional space--time)
introduced spinors for the representation of the wave functions.

In general relativity theory spinors and the Dirac equations on (pseudo)
Riemannian spaces were defined in 1929 by H. Weyl \cite{weyl}, V. Fock \cite%
{foc} and E. Schr\"{o}dinger \cite{schr}. The books \cite{pen,penr1,penr2}
by R. Penrose and W. Rindler monograph  summarize the spinor and twistor
methods in space--time geometry (see additional references \cite%
{hladik,benn,morand,lue,tur,car} on Clifford structures and spinor theory).

Spinor variables were introduced in Finsler geometries by Y. Takano in 1983 %
\cite{t1} where he dismissed anisotropic dependencies not only on vectors on
the tangent bundle but on some spinor variables in a spinor bundle on a
space--time manifold. Then generalized Finsler geometries, with spinor
variables, were developed by T. Ono and Y. Takano in a series of
publications during 1990--1993 \cite{ot1,ot2,ot3,ono}. The next steps were
investigations of anisotropic and deformed geometries with spinor and vector
variables and applications in gauge and gravity theories elaborated by P.
Stavrinosand his students, S. Koutroubis, P. Manouselis, and  V. Balan
beginning 1994 \cite{sk,sm2,sm3,sbmp,sbmp1}. In those works the authors
assumed that some spinor variables may be introduced in a Finsler-like way
but they did not relate the Finlser metric to a Clifford structure and
restricted the spinor--gauge Finsler constructions only for antisymmetric
spinor metrics on two--spinor fibers with possible generalizations to four
dimensional Dirac spinors.

Isotopic spinors, related with $SU(2)$ internal structural groups, were
considered in generalized Finsler gravity and gauge theories also by G.
Asanov and S. Ponomarenko \cite{asa88} in 1988. In that book, and in other
papers on Finsler geometry with spinor variables, the authors did not
investigate the possibility of introducing a rigorous mathematical
definition of spinors on spaces with generic local anisotropy.

An alternative approach to spinor differential geometry and generalized
Finsler spaces was elaborated, beginning 1994, in a series of papers and
communications by S. Vacaru with participation of S. Ostaf \cite%
{vdeb,vo1,vrom,vbm}. This direction originates from Clifford algebras and
Clifford bundles \cite{kar,tur} and Penrose's spinor and twistor space--time
geometry \cite{pen,penr1,penr2}, which were re--considered for the case of
nearly autoparallel maps (generalized conformal transforms) in Refs. \cite%
{v87,vkaz,vit}. In the works \cite{viasm1,vjmp,vsp2}, a rigorous definition
of spinors for Finsler spaces, and their generalizations, was given. It was
proven that a Finsler, or Lagrange, metric (in a tangent, or, more
generally, in a vector bundle) induces naturally a distinguished Clifford
(spinor) structure which is locally adapted to the nonlinear connection
structure. Such spinor spaces could be defined for arbitrary dimensions of
base and fiber subspaces, their spinor metrics are symmetric, antisymmetric
or nonsymmetric, depending on the corresponding base and fiber dimensions.
This work resulted in the formation of spinor differential geometry of
generalized Finsler spaces and developed a number of geometric applications
to the theory of gravitational and matter field interactions with generic
local anisotropy.

The geometry of anisotropic spinors and (distinguished by nonlinear
connections) Clifford structures was elaborated for higher order anisotropic
spaces \cite{vsp1,vhsp,vbook} and, more recently, for Hamilton and Lagrange
spaces \cite{vsv}.

We emphasize that the theory of anisotropic spinors may be related not only
to generalized Finsler, Lagrange, Cartan and Hamilton spaces or their higher
order generalizations, but also to anholonomic frames with associated
nonlinear connections which appear naturally even in (pseudo) Riemannian
geometry if off--diagonal metrics are considered \cite%
{vb,vbh,vkinet,vtheor,vtor}. In order to construct exact solutions of the
Einstein equations in general relativity and extra dimension gravity (for
lower dimensions see \cite{vtherm,vsgab,vsgon}), it is more convenient to
diagonalize space--time metrics by using some anholonomic transforms. As a
result one induces locally anisotropic structures on space--time which are
related to anholonomic (anisotropic) spinor structures.

The main purpose of the present review is to present a detailed summary and
new results on spinor differential geometry for generalized Finsler spaces
and (pseudo) Riemannian space--times provided with anholonomic frame and
associated nonlinear connection structure, to discuss and compare the
existing approaches and to consider applications to modern gravity and gauge
theories.

This article is organized in four Chapters:

In Chapter 1, we give the basic definitions from the theory of generalized
Finsler, Lagrange, Cartan and Hamilton spaces on vector and co--vector
(tangent and co--tangent spaces) and their generalizations for higher order
vector--covector bundles following the monographs \cite{ma94,mhss,vbook}.

Chapter 2 is a generalization of the results on Clifford
structures for higher order vector bundles.

Chapter 3 is devoted to the differential geometry of Spinors in
Higher Order Anisotropic Spaces.

Chapter 4 contains geometric applications to the theory of locally
anisotropic interactions  \cite{vsp2,vhsp}.

\subsection*{Acknowledgements}

~~ The S. V. work is supported by a NATO/Portugal fellowship at
CENTRA, Instituto Superior Tecnico, Lisbon. The author is
grateful to R. Ablamowicz, John Ryan and B. Fauser for
collaboration and support of his participation at ''The 6th
International Conference on Clifford Algebras'', Cookeville,
Tennessee, USA (May, 20-25, 2002). He would like to thank J. P.
S. Lemos, R. Miron and M. Anastasiei for hospitality and support.

\chapter{(Co-) Vector Bundles and Nonlinear Connections}

In this Chapter the space--time geometry is modeled not only on a (pseudo)
Riemannian manifold $V^{[n+m]}$ of dimension $n+m$ but it is considered on a
vector bundle (or its dual, covector bundle) being, for simplicity, locally
trivial with a base space $M$ of dimension $n$ and a typical fiber $F $
(cofiber $F^{*}$) of dimension $m,$ or as a higher order extended
vector/covector bundle (we follow the geometric constructions and
definitions of monographs \cite{ma94,ma87,mhss,m1,m2} which were generalized
for vector superbundles in Refs. \cite{vstr2,vbook}). Such fibered
space--times (in general, with extra dimensions and duality relations) are
supposed to be provided with compatible structures of nonlinear and linear
connections and (pseudo) Riemannian metric. For the particular cases when:\
a) the total space of the vector bundle is substituted by a
pseudo--Riemannian manifold of necessary signature we can model the usual
pseudo--Riemannian space--time from the Einstein gravity theory with field
equations and geometric objects defined with respect to some classes of
moving anholonomic frames with associated nonlinear connection structure; b)
if the dimensions of the base and fiber spaces are identical, $n=m,$ for the
first order anisotropy, we obtain the tangent bundle $TM.$

Such both (pseudo) Riemanian spaces and vector/covector (in particular
cases, tangent/cotangent) bundles of metric signature (-,+,...,+) enabled
with compatible fibered and/or anholonomic structures, the metric in the
total space being a solution of the Einstein equations, will be called
\textbf{anisotropic space--times}. If the anholonomic structure with
associated nonlinear connection is modeled on higher order vector/covector
bundles we shall use the term of \textbf{higher order anisotropic
space--time.}

The geometric constructions are outlined as to present the main concepts and
formulas in a unique way for both type of vector and covector structures. In
this part of the book we usually shall omit proofs which can be found in the
mentioned monographs \cite{ma87,ma94,m1,m2,mhss,vbook}.

\section{Vector and Covector Bundles}

In this Section we introduce the basic definitions and denotations for
vector and tangent (and theirs dual spaces) bundles and higher order
vector/covector bundle geometry.

\subsection{Vector and tangent bundles}

A locally trivial \textbf{vector bundle}, in brief,
\textbf{v--bundle}, \index{vector bundle} \index{v--bundle}
$\mathcal{E}=\left( E,\pi ,M,Gr,F\right) $ is introduced as a set
of spaces and surjective map with the properties that a real
vector space $F=\mathcal{R}^m$ of dimension $m$ ($\dim F=m,$
$\mathcal{R\ }$ denotes the real number field) defines the
typical fiber, the structural group
\index{v--bundle} is chosen to be the group of automorphisms of $\mathcal{R}%
^m,$ i. e. $Gr=GL\left( m,\mathcal{R}\right) ,\ $ and $\pi :E\rightarrow M$
is a differentiable surjection of a differentiable
\index{manifold} manifold $E$ (total space, $\dim E=n+m)$ to a
differentiable manifold $M$ $\left( \mbox{base
space, }\dim M=n\right) .$ Local coordinates on $\mathcal{E}$ are denoted $%
u^{{\alpha }}=\left( x^{{i}},y^{{a}}\right) ,$ or in brief
${u=\left( x,y\right) }$ (the Latin indices
${i,j,k,...}=1,2,...,n$ define coordinates of geometrical objects
with respect to a local frame on base space $M;$\ the Latin
indices ${a,b,c,...=}1,2,...,m$ define fiber coordinates of
geometrical objects and the Greek indices ${\alpha ,\beta ,\gamma
,...}$ are considered as cumulative ones for coordinates of
objects defined on the total space of a v-bundle).

Coordinate trans\-forms $u^{{\ \alpha ^{\prime }\ }} = u^{{\ \alpha ^{\prime
}\ }}\left( u^{{\ \alpha }}\right) $ on a v--bundle $\mathcal{E}$ are
defined as
\begin{equation*}
\left(x^{{\ i}}, y^{{\ a}}\right) \rightarrow \left( x^{{\ i^{\prime }\ }},
y^{{\ a^{\prime }}}\right) ,
\end{equation*}
where
\begin{equation}
x^{{i^{\prime }\ }} = x^{{\ i^{\prime }\ }}(x^{{\ i}}),\qquad y^{{a^{\prime
}\ }} = K_{{\ a\ }}^{{\ a^{\prime }}}(x^{{i\ }})y^{{a}}  \label{coordtr}
\end{equation}
and matrix $K_{{\ a\ }}^{{\ a^{\prime }}}(x^{{\ i\ }}) \in GL\left( m,%
\mathcal{R}\right) $ are functions of necessary smoothness class.

A local coordinate parametrization of v--bundle $\mathcal{E}$ naturally
defines a coordinate basis
\begin{equation}
\partial _{\alpha }=%
\frac{\partial }{\partial u^{\alpha }}=\left( \partial _{i}=\frac{\partial }{%
\partial x^{i}},\ \partial _{a}=\frac{\partial }{\partial y^{a}}\right) ,
\label{pder}
\end{equation}%
and the reciprocal to (\ref{pder}) coordinate basis
\begin{equation}
d^{\alpha }=du^{\alpha }=(d^{i}=dx^{i},\ d^{a}=dy^{a})  \label{pdif}
\end{equation}%
which is uniquely defined from the equations
\begin{equation*}
d^{\alpha }\circ \partial _{\beta }=\delta _{\beta }^{\alpha },
\end{equation*}%
where $\delta _{\beta }^{\alpha }$ is the Kronecher symbol and by ''$\circ $%
'' we denote the inner (scalar) product in the tangent bundle $\mathcal{TE}$

A \textbf{tangent bundle}
\index{tangent bundle} (in brief, \textbf{t--bundle})
\index{t--bundle} $(TM,\pi ,M)$ to a manifold $M$ can be defined as a
particular case of a v--bundle when the dimension of the base and fiber
spaces (the last one considered as the tangent subspace) are identic, $n=m.$
In this case both type of indices $i,k,...$ and $a,b,...$ take the same
values $1,2,...n$. For t--bundles the matrices of fiber coordinates
transforms from (\ref{coordtr}) can be written $K_{{\ i\ }}^{{\ i^{\prime }}%
} = {\partial x^{i^{\prime }}}/ {\partial x^i}.$

We shall distinguish the base and fiber indices and values which is
necessary for our further geometric and physical applications.

\subsection{Covector and cotangent bundles}

We shall also use the concept of \textbf{covector bundle},
\index{covector bundle} (in brief, \textbf{cv--bundles)}
\index{ cv--bundle}\newline
$%
\breve {\mathcal{E}}=\left({\breve E},\pi ^{*} ,M,Gr,F^{*}\right) $, which
is introduced as a dual vector bundle for which the typical fiber $F^{*}$
(cofiber)
\index{cofiber} is considered to be the dual vector space (covector space)
\index{covector space} to the vector space $F.$ The fiber coordinates $p_a$
of $\breve E$ are dual to $y^a$ in $E.$ The local coordinates on total space
$\breve E$ are denoted $\breve{u}=(x,p)=(x^i,p_a).$ The coordinate transform
on $\breve E,$
\begin{equation*}
{\breve{u}}=(x^i,p_a)\to {\breve{u}}^{\prime }=(x^{i^{\prime }},p_{a^{\prime
}}),
\end{equation*}
are written
\begin{equation}
x^{{i^{\prime }\ }} = x^{{\ i^{\prime }\ }}(x^{{\ i}}),\qquad p_{{a^{\prime
}\ }} = K^{{\ a\ }}_{{\ a^{\prime }}}(x^{{i\ }})p_{{a}}.  \label{coordtrd}
\end{equation}
The coordinate bases on $E^{*}$ are denoted
\begin{equation}
{%
\breve{\partial}}_\alpha =\frac{\breve{\partial}}{\partial u^\alpha }=\left(
\partial _i=\frac \partial {\partial x^i},{{\breve {\partial}}^a}=\frac{%
\breve{\partial}}{\partial p_a}\right)  \label{pderct}
\end{equation}
and
\begin{equation}
{\breve{d}}^\alpha ={\breve{d}}u^\alpha =\left( d^i=dx^i,{\breve{d}}_a=
dp_a\right) .  \label{pdifct}
\end{equation}
We shall use ''breve'' symbols in order to distinguish the geometrical
objects on a cv--bundle $\mathcal{E}^{*}$ from those on a v--bundle $%
\mathcal{E}$.

As a particular case with the same dimension of base space and cofiber one
obtains the \textbf{cotangent bundle}
\index{cotangent bundle} $(T^{*}M,\pi ^{*},M)$ , in brief, \textbf{%
ct--bundle,}
\index{ct--bundle} being dual to $TM.$ The fibre coordinates $p_i$ of $%
T^{*}M $ are dual to $y^i$ in $TM.$ The coordinate transforms (\ref{coordtrd}%
) on $T^{*}M$ are stated by some matrices $K^{k}_{{\ k^{\prime }}}(x^i)= {%
\partial x^k}/{\partial x^{k^{\prime }}}.$

In our further considerations we shall distinguish the base and cofiber
indices.

\subsection{Higher order vector/covector bundles}

The geometry of higher order tangent and cotangent bundles provided with
nonlinear connection structure was elaborated in Refs. \cite%
{m1,m2,mirata,mhss} following the aim of geometrization of higher order
Lagrange and Hamilton mechanics. In this case we have base spaces and fibers
of the same dimension. In order to develop the approach to modern high
energy physics (in superstring and Kaluza--Klein theories) one had to
introduce (in Refs \cite{vsp1,vhsp,vbook,vstr2}) the concept of higher order
vector bundle with the fibers consisting from finite 'shells" of vector, or
covector, spaces of different dimensions not obligatory coinciding with the
base space dimension.

\begin{definition}
A distinguished vector/covector space, in brief dvc--space,
\index{dvc--space} of type
\begin{equation}
{%
\tilde{F}}=F[v(1),v(2),cv(3),...,cv(z-1),v(z)]  \label{orient}
\end{equation}
is a vector space decomposed into an invariant oriented direct sum
\begin{equation*}
{\tilde{F}}=F_{(1)}\oplus F_{(2)}\oplus F_{(3)}^{*}\oplus ...\oplus
F_{(z-1)}^{*}\oplus F_{(z)}
\end{equation*}
of vector spaces $F_{(1)},F_{(2)},...,F_{(z)}$ of respective dimensions
\begin{equation*}
dimF_{(1)}=m_1,dimF_{(2)}=m_2,...,dimF_{(z)}=m_z
\end{equation*}
and of covector spaces $F_{(3)}^{*},...,F_{(z-1)}^{*}$ of respective
dimensions
\begin{equation*}
dimF_{(3)}^{*}=m_3^{*},...,dimF_{(z-1)}^{*}=m_{(z-1)}^{*}.
\end{equation*}
\end{definition}

As a particular case we obtain a distinguished vector space, in brief
dv--space
\index{dv--space} (a distinguished covector space, in brief dcv--space),
\index{dcv--space} if all components of the sum are vector (covector)
spaces. We note that we have fixed for simplicity an orientation of
vector/covector subspaces like in (\ref{orient}); in general there are
possible various type of orientations, number of subspaces and dimensions of
subspaces.

Coordinates on ${\tilde F}$ are denoted
\begin{equation*}
\tilde{y} = (y_{(1)},y_{(2)},p_{(3)},...,p_{(z-1)},y_{(z)})= \{y^{<\alpha
_z>}\} = (y^{a_1},y^{a_2},p_{a_3},...,p_{a_{z-1}},y^{a_z}),
\end{equation*}
where indices run corresponding values:\
\begin{equation*}
a_1 = 1,2,...,m_1;\ a_2 = 1,2,...,m_2,\ ..., a_z = 1,2,...,m_z.
\end{equation*}

\begin{definition}
A higher order vector/covector bundle (in brief, hvc-\--bund\-le)
\index{hvc--bundle} of type ${%
\tilde{\mathcal{E}}}={\tilde{\mathcal{E}}}[v(1),v(2),cv(3),...,cv(z-1),v(z)]$
is a vector bundle ${\tilde{\mathcal{E}}}=({\tilde{E}},p^{<d>},{\tilde{F}}%
,M) $ with corresponding total, ${\tilde{E}}$, and base, $M,$ spaces,
surjective projection $p^{<d>}:\ {\tilde{E}}\to M$ and typical fiber ${%
\tilde{F}}.$
\end{definition}

We define higher order vector (covector) bundles, in brief, hv--bundles (in
brief, hcv--bundles), if the typical fibre is a dv--space (dcv--space) as
particular cases of hvc--bundles.

A hvc--bundle is constructed as an oriented set of enveloping 'shell by
shell' v--bundles and/or cv--bundles,
\begin{equation*}
p^{<s>}:\ {\tilde{E}}^{<s>}\rightarrow {\tilde{E}}^{<s-1>},
\end{equation*}%
where we use the index $<s>=0,1,2,...,z$ in order to enumerate the shells,
when ${\tilde{E}}^{<0>}=M.$ Local coordinates on ${\tilde{E}}^{<s>}$ are
denoted
\begin{equation}
{\tilde{u}}_{(s)}=(x,{\tilde{y}}%
_{<s>})=(x,y_{(1)},y_{(2)},p_{(3)},...,y_{(s)})=(x^{i},y^{a_{1}},y^{a_{2}},p_{a_{3}},...,y^{a_{s}}).
\notag
\end{equation}%
If $<s>=<z>$ we obtain a complete coordinate system on ${\tilde{\mathcal{E}}}
$ denoted in brief
\begin{equation*}
\tilde{u}=(x,{\tilde{y}})=\tilde{u}^{\alpha
}=(x^{i}=y^{a_{0}},y^{a_{1}},y^{a_{2}},p_{a_{3}},...,p_{a_{z-1}},y^{a_{z}}).
\end{equation*}%
We shall use the general commutative indices $\alpha ,\beta ,...$ for
objects on hvc---bundles which are marked by tilde, like $\tilde{u},\tilde{u}%
^{\alpha },...,$ ${\tilde{E}}^{<s>},....$

The coordinate transforms for a hvc--bundle ${\tilde{\mathcal{E}},}$
\begin{equation*}
\tilde{u}=(x,{\tilde{y}})\rightarrow \tilde{u}^{\prime }=(x^{\prime },{%
\tilde{y}}^{\prime })
\end{equation*}
are given by recurrent formulas
\begin{eqnarray*}
x^{i^{\prime }} &=&x^{i^{\prime }}\left( x^i\right) ,\ rank\left( \frac{%
\partial x^{i^{\prime }}}{\partial x^i}\right) =n; \\
y^{a_1^{\prime }} &=&K_{a_1}^{a_1^{^{\prime
}}}(x)y^{a_1},K_{a_1}^{a_1^{^{\prime }}}\in GL(m_1,\mathcal{R}); \\
y^{a_2^{\prime }} &=&K_{a_2}^{a_2^{^{\prime
}}}(x,y_{(1)})y^{a_2},K_{a_2}^{a_2^{^{\prime }}}\in GL(m_2,\mathcal{R}); \\
p_{a_3^{\prime }} &=&K_{a_3^{\prime
}}^{a_3}(x,y_{(1)},y_{(2)})p_{a_3},K_{a_3^{\prime }}^{a_3}\in GL(m_3,%
\mathcal{R}); \\
y^{a_4^{\prime }} &=&K_{a_4}^{a_4^{^{\prime
}}}(x,y_{(1)},y_{(2)},p_{(3)})y^{a_4},K_{a_4}^{a_4^{^{\prime }}}\in GL(m_4,%
\mathcal{R}); \\
&&................ \\
p_{a_{z-1}^{\prime }} &=&K_{a_{z-1}^{\prime
}}^{a_{z-1}}(x,y_{(1)},y_{(2)},p_{(3)},...,y_{(z-2)})p_{a_{z-1}},K_{a_{z-1}^{\prime }}^{a_{z-1}}\in GL(m_{z-1},%
\mathcal{R}); \\
y^{a_z^{\prime }}
&=&K_{a_z}^{a_z^{^{%
\prime}}}(x,y_{(1)},y_{(2)},p_{(3)},...,y_{(z-2)},p_{a_{z-1}})y^{a_z},K_{a_z}^{a_z^{^{\prime }}}\in GL(m_z,%
\mathcal{R}),
\end{eqnarray*}
where, for instance. by $GL(m_2,\mathcal{R})$ we denoted the group of linear
transforms of a real vector space of dimension $m_2.$

The coordinate bases on ${\tilde{\mathcal{E}}}$ are denoted
\begin{eqnarray}
{\tilde{\partial}}_{\alpha } &=&\frac{{\tilde{\partial}}}{\partial u^{\alpha
}}=\left( \partial _{i}=\frac{\partial }{\partial x^{i}},\partial _{a_{1}}=%
\frac{\partial }{\partial y^{a_{1}}},\partial _{a_{2}}=\frac{\partial }{%
\partial y^{a_{2}}},{{\breve{\partial}}}^{a_{3}}=\frac{\breve{\partial}}{%
\partial p_{a_{3}}},...,\partial _{a_{z}}=\frac{\partial }{\partial y^{a_{z}}%
}\right)  \label{pderho} \\
&&  \notag
\end{eqnarray}%
and
\begin{eqnarray}
{\tilde{d}}^{\alpha } &=&{\tilde{d}}u^{\alpha }=\left(
d^{i}=dx^{i},d^{a_{1}}=dy^{a_{1}},d^{a_{2}}=dy^{a_{2}},{\breve{d}}%
_{a_{3}}=dp_{a_{3}},...,d^{a_{z}}=dy^{a_{z}}\right) .  \label{pdifho} \\
&&  \notag
\end{eqnarray}

We end this subsection with two examples of higher order tangent /
co\-tan\-gent bundles (when the dimensions of fibers/cofibers coincide with
the dimension of bundle space, see Refs. \cite{m1,m2,mirata,mhss}).

\subsubsection{Osculator bundle}

\index{osculator} \label{oscsect}

The $k$--osculator bundle
\index{$k$--osculator} is identified with the $k$--tangent bundle\newline
$\left( T^kM,p^{(k)},M\right) $ of a $n$--dimensional manifold $M.$ We
denote the local coordinates
\begin{equation*}
{\tilde{u}}^\alpha =\left( x^i,y_{(1)}^i,...,y_{(k)}^i\right) ,
\end{equation*}
where we have identified $y_{(1)}^i\simeq y^{a_1},...,y_{(k)}^i\simeq
y^{a_k},k=z,$ in order to to have similarity with denotations from \cite%
{mhss}. The coordinate transforms
\begin{equation*}
{\tilde{u}}^{\alpha ^{\prime }}\rightarrow {\tilde{u}}^{\alpha ^{\prime
}}\left( {\tilde{u}}^\alpha \right)
\end{equation*}
preserving the structure of such higher order vector bundles are
parametrized
\begin{eqnarray*}
x^{i^{\prime }} &=&x^{i^{\prime }}\left( x^i\right) ,\det \left(
\frac{\partial x^{i^{\prime }}}{\partial x^i}\right) \neq 0, \\
y_{(1)}^{i^{\prime }} &=&\frac{\partial x^{i^{\prime }}}{\partial x^i}%
y_{(1)}^i, \\
2y_{(2)}^{i^{\prime }} &=&\frac{\partial y_{(1)}^{i^{\prime }}}{\partial x^i}%
y_{(1)}^i+2\frac{\partial y_{(1)}^{i^{\prime }}}{\partial y^i}y_{(2)}^i, \\
&&................... \\
ky_{(k)}^{i^{\prime }} &=&\frac{\partial y_{(1)}^{i^{\prime }}}{\partial x^i}%
y_{(1)}^i+...+k\frac{\partial y_{(k-1)}^{i^{\prime }}}{\partial y_{(k-1)}^i}%
y_{(k)}^i,
\end{eqnarray*}
where the equalities
\begin{equation*}
\frac{\partial y_{(s)}^{i^{\prime }}}{\partial x^i}=\frac{\partial
y_{(s+1)}^{i^{\prime }}}{\partial y_{(1)}^i}=...=\frac{\partial
y_{(k)}^{i^{\prime }}}{\partial y_{(k-s)}^i}
\end{equation*}
hold for $s=0,...,k-1$ and $y_{(0)}^i=x^i.$

The natural coordinate frame on $\left( T^kM,p^{(k)},M\right) $ is defined
\begin{equation*}
{\tilde{\partial}}_\alpha =\left( \frac \partial {\partial x^i},\frac
\partial {\partial y_{(1)}^i},...,\frac \partial {\partial y_{(k)}^i}\right)
\end{equation*}
and the coframe is
\begin{equation*}
{\tilde d}_\alpha =\left( dx^i,dy_{(1)}^i,...,dy_{(k)}^i\right) .
\end{equation*}
These formulas are respectively some particular cases of $\left( \ref{pderho}%
\right) $ and $\left( \ref{pdifho}\right) .$

\subsubsection{The dual bundle of k--osculator bundle}

\label{doscsect}

This higher order vector/covector bundle, denoted as $\left(
T^{*k}M,p^{*k},M\right) ,$ is defined as the dual bundle to the k--tangent
bundle $\left( T^kM,p^k,M\right) .$ The local coordinates (parametrized as
in the previous paragraph) are
\begin{equation*}
\tilde u=\left( x,y_{(1)},...,y_{(k-1)},p\right) =\left(
x^i,y_{(1)}^i,...,y_{(k-1)}^i,p_i\right) \in T^{*k}M.
\end{equation*}
The coordinate transforms on $\left( T^{*k}M,p^{*k},M\right) $ are
\begin{eqnarray*}
x^{i^{\prime }} &=&x^{i^{\prime }}\left( x^i\right) ,\det \left( \frac{%
\partial x^{i^{\prime }}}{\partial x^i}\right) \neq 0, \\
y_{(1)}^{i^{\prime }} &=&\frac{\partial x^{i^{\prime }}}{\partial x^i}%
y_{(1)}^i, \\
2y_{(2)}^{i^{\prime }} &=&\frac{\partial y_{(1)}^{i^{\prime }}}{\partial x^i}%
y_{(1)}^i+2\frac{\partial y_{(1)}^{i^{\prime }}}{\partial y^i}y_{(2)}^i, \\
&&................... \\
(k-1)y_{(k-1)}^{i^{\prime }} &=&\frac{\partial y_{(k-2)}^{i^{\prime }}}{%
\partial x^i}y_{(1)}^i+...+k\frac{\partial y_{(k-1)}^{i^{\prime }}}{\partial
y_{(k-2)}^i}y_{(k-1)}^i, \\
p_{i^{\prime }} &=&\frac{\partial x^i}{\partial x^{i^{\prime }}}p_i,
\end{eqnarray*}
where the equalities
\begin{equation*}
\frac{\partial y_{(s)}^{i^{\prime }}}{\partial x^i}=\frac{\partial
y_{(s+1)}^{i^{\prime }}}{\partial y_{(1)}^i}=...=\frac{\partial
y_{(k-1)}^{i^{\prime }}}{\partial y_{(k-1-s)}^i}
\end{equation*}
hold for $s=0,...,k-2$ and $y_{(0)}^i=x^i.$

The natural coordinate frame on $\left( T^{\ast k}M,p^{\ast (k)},M\right) $
is defined
\begin{equation*}
{\tilde{\partial}}_{\alpha }=\left( \frac{\partial }{\partial x^{i}},\frac{%
\partial }{\partial y_{(1)}^{i}},...,\frac{\partial }{\partial y_{(k-1)}^{i}}%
,\frac{\partial }{\partial p_{i}}\right)
\end{equation*}
and the coframe is
\begin{equation*}
{\tilde{d}}_{\alpha }=\left(
dx^{i},dy_{(1)}^{i},...,dy_{(k-1)}^{i},dp_{i}\right) .
\end{equation*}
These formulas are respectively another particular cases of $\left( \ref%
{pderho}\right) $ and $\left( \ref{pdifho}\right) .$

\section{Nonlinear Connections}

\label{subsnc}

The concept of \textbf{nonlinear connection,}
\index{nonlinear connection} in brief, N-connection,
\index{N--connection} is fundamental in the geometry of vector bundles and
anisotropic spaces (see a detailed study and basic references in \cite%
{ma87,ma94}). A rigorous mathematical definition is possible by using the
formalism of exact sequences of vector bundles.

\subsection{N--connections in vector bundles}

Let $\mathcal{E=}=(E,p,M)$ be a v--bundle with typical fiber
$\mathcal{R}^m$ and $\pi ^T:\ TE\to TM$ being the differential of
the map $P$ which is a fibre--preserving morphism of the tangent
bundle $TE,\tau _E,E)\to E$ and of tangent bundle $(TM,\tau
,M)\to M.$ The kernel of the vector bundle morphism, denoted as
$(VE,\tau _V,E),$ is called the \textbf{vertical
subbundle} over $E,$ which is a vector subbundle of the vector bundle $%
(TE,\tau _E,E).$

A vector $X_u$ tangent to a point $u\in E$ is locally written as
\begin{equation*}
(x,y,X,Y)=(x^i,y^a,X^i,Y^a),
\end{equation*}
where the coordinates $(X^i,Y^a)$ are defined by the equality
\begin{equation*}
X_u=X^i\partial _i+Y^a\partial _a.
\end{equation*}
We have\ $\pi ^T(x,y,X,Y)=(x,X).$ Thus the submanifold $VE$ contains the
elements which are locally represented as $(x,y,0,Y).$

\begin{definition}
\label{ncon} A nonlinear connection $\mathbf{N}$ in a vector bundle $%
\mathcal{E}=(E,\pi ,M)$ is the splitting on the left of the exact sequence
\begin{equation*}
0\mapsto VE\mapsto TE\mapsto TE/VE\mapsto 0
\end{equation*}
where $TE/VE$ is the factor bundle.
\end{definition}

By definition (\ref{ncon}) it is defined a morphism of vector bundles $C:\
TE\to VE$ such the superposition of maps $C\circ i$ is the identity on $VE,$
where $i:\ VE\mapsto VE.$ The kernel of the morphism $C$ is a vector
subbundle of $(TE,\tau _E,E)$ which is the horizontal subbundle, denoted by $%
(HE,\tau _H,E).$ Consequently, we can prove that in a v-bundle $\mathcal{E}$
a N--connection can be introduced as a distribution
\begin{equation*}
\{N:\ E_u\rightarrow H_uE,T_uE=H_uE\oplus V_uE\}
\end{equation*}
for every point $u\in E$ defining a global decomposition, as a
Whitney sum, into horizontal, $H\mathcal{E},$ and vertical,
$V\mathcal{E},$ subbundles of the tangent bundle $T\mathcal{E}$
\begin{equation}
T\mathcal{E}=H\mathcal{E}\oplus V\mathcal{E}.  \label{whitney}
\end{equation}

Locally a N-connection in a v--bundle $\mathcal{E}$ is given by its
coefficients\newline
$N_{{\ i}}^{{a}}({\ u})=N_{{i}}^{{a}}({x,y})$ with respect to bases
(\ref{pder}) and (\ref{pdif})
\begin{equation*}
\mathbf{N}=N_i^{~a}(u)d^i\otimes \partial _a.
\end{equation*}

We note that a linear connection in a v--bundle $\mathcal{E}$ can be
considered as a particular case of a N--connection when $%
N_i^{~a}(x,y)=K_{bi}^a\left( x\right) y^b,$ where functions $K_{ai}^b\left(
x\right) $ on the base $M$ are called the Christoffel coefficients.

\subsection{N--connections in covector bundles:}

A nonlinear connection in a cv--bundle ${%
\breve{\mathcal{E}}}$ (in brief a \v N--connection) can be introduces in a
similar fashion as for v--bundles by reconsidering the corresponding
definitions for cv--bundles. For instance, it is stated by a Whitney sum,
into horizontal,$H{\breve{\mathcal{E}}},\ $ and vertical, $V{\breve{\mathcal{%
E}}},$ subbundles of the tangent bundle $T{\breve{\mathcal{E}}}:$
\begin{equation}
T{\breve{\mathcal{E}}}=H{\breve{\mathcal{E}}}\oplus V{\breve{\mathcal{E}}}.
\label{whitneyc}
\end{equation}

Hereafter, for the sake of brevity we shall omit details on definition of
geometrical objects on cv--bundles if they are very similar to those for
v--bundles:\ we shall present only the basic formulas by emphasizing the
most important particularities and differences.

\begin{definition}
\label{ctvn} A \v{N}--connection on ${\breve{\mathcal{E}}}$ is a
differentiable distribution
\begin{equation*}
\breve{N}:\ {\breve{\mathcal{E}}}\rightarrow {\breve{N}}_{u}\in T_{u}^{\ast }%
{\breve{\mathcal{E}}}
\end{equation*}%
which is suplimentary to the vertical distribution $V,$ i. e.
\begin{equation*}
T_{u}{\breve{\mathcal{E}}}={\breve{N}}_{u}\oplus {\breve{V}}_{u},\forall {%
\breve{\mathcal{E}}}.
\end{equation*}
\end{definition}

The same definition is true for \v N--connections in ct--bundles, we have to
change in the definition (\ref{ctvn}) the symbol ${\breve{\mathcal{E}}}$
into $T^*M.$

A \v N--connection in a cv--bundle ${\breve{\mathcal{E}}}$ is given locally
by its coefficients\newline
${\breve N}_{{\ ia}}({\ u})= {\breve N}_{{ia}}({x,p})$ with respect to bases
(\ref{pder}) and (\ref{pdif})
\begin{equation*}
\mathbf{\breve N}={\breve N}_{ia}(u)d^i\otimes {{\breve{\partial}}^a}.
\end{equation*}

We emphasize that if a N--connection is introduced in a v--bundle
(cv--bundle) we have to adapt the geometric constructions to the
N--connection structure.

\subsection{N--connections in higher order bundles}

The concept of N--connection can be defined for higher order vector /
covec\-tor bundle in a standard manner like in the usual vector bundles:

\begin{definition}
A nonlinear connection ${\tilde{\mathbf{N}}}$ in hvc--bundle
\index{N--connection}
\begin{equation*}
{%
\tilde{\mathcal{E}}}={\tilde{\mathcal{E}}}[v(1),v(2),cv(3),...,cv(z-1),v(z)]
\end{equation*}
is a splitting of the left of the exact sequence
\begin{equation}
0\to V{\tilde{\mathcal{E}}}\to T{\tilde{\mathcal{E}}}\to T{\tilde{\mathcal{E}%
}}/V{\tilde{\mathcal{E}}}\to 0  \label{exacts}
\end{equation}
\end{definition}

We can associate sequences of type (\ref{exacts}) to every mappings of
intermediary subbundles. For simplicity, we present here the Whitney
decomposition
\begin{equation*}
T{\tilde{\mathcal{E}}}=H{\tilde{\mathcal{E}}}\oplus V_{v(1)}{\tilde{\mathcal{%
E}}}\oplus V_{v(2)}{\tilde{\mathcal{E}}}\oplus V_{cv(3)}^{*}{\tilde{\mathcal{%
E}}}\oplus ....\oplus V_{cv(z-1)}^{*}{\tilde{\mathcal{E}}}\oplus V_{v(z)}{%
\tilde{\mathcal{E}}}.
\end{equation*}
Locally a N--connection ${\tilde{\mathbf{N}}}$ in ${\tilde{\mathcal{E}}}$ is
given by its coefficients
\begin{equation}
\begin{array}{llllll}
N_i^{~a_1}, & N_i^{~a_2}, & N_{ia_3}, & ..., & N_{ia_{z-1}}, & N_i^{~a_z},
\\
0, & N_{a_1}^{~a_2}, & N_{a_1a_3}, & ..., & N_{a_1a_{z-1}}, & N_{a_1}^{~a_z},
\\
0, & 0, & N_{a_2a_3}, & ..., & N_{a_2a_{z-1}}, & N_{a_2}^{~a_z}, \\
..., & ..., & ..., & ..., & ..., & ..., \\
0, & 0, & 0, & ..., & N_{a_{z-2}~a_{z-1}}, & N_{a_{z-2}}^{~a_z}, \\
0, & 0, & 0, & ..., & 0, & N^{a_{z-1}a_z},%
\end{array}
\label{nconho}
\end{equation}
which are given with respect to the components of bases $\left( \ref{pderho}%
\right) $ and $\left( \ref{pdifho}\right) .$

We end this subsection with two examples of N--connections in
higher order vector/covector bundles:

\subsubsection{N--connection in osculator bundle}

\label{nconoscs}

Let us consider the second order of osculator bundle (see subsection (\ref%
{oscsect})) $T^2M=Osc^2M.$ A N--connection ${\tilde {\mathbf{N}}}$ in $%
Osc^2M $ is associated to a Whitney sum
\begin{equation*}
TT^2M=NT^2M\oplus VT^2M
\end{equation*}
which defines in every point $\tilde u\in T^2M$ a distribution
\begin{equation*}
T_uT^2M=N_0\left( \tilde u\right) \oplus N_1\left( \tilde u\right) \oplus
VT^2M.
\end{equation*}
We can parametrize ${\tilde {\mathbf{N}}}$ with respect to natural
coordinate bases as
\begin{equation}  \label{nconoscsf}
\begin{array}{ll}
N_i^{a_1}, & N_i^{a_2}, \\
0, & N_{a_1}^{a_2}.%
\end{array}%
\end{equation}
As a particular case we can consider $N_{a_1}^{a_2}=0.$

\subsubsection{N--connection in dual osculator bundle}

\label{nconoscds}

In a similar fashion we can take the bundle $\left( T^{*2}M,p^{*2},M\right) $
being dual bundle to the $Osc^2M$ (see subsection (\ref{doscsect})). We have
\begin{equation*}
T^{*2}M=TM\otimes T^{*}M.
\end{equation*}
The local coefficients of a N--connection in $\left( T^{*2}M,p^{*2},M\right)
$ are parametrizied
\begin{equation}  \label{nconoscsdf}
\begin{array}{ll}
N_i^{~a_1}, & N_{ia_2}, \\
0, & N_{a_1a_2}.%
\end{array}%
\end{equation}
We can choose a particular case when $N_{a_1a_2}=0.$

\subsection{Anholonomic frames and N--connections}

Having defined a N--connection structure in a (vector, covector, or higher
order vector / covenctor) bundle we can adapt to this structure, (by
'N--elonga\-ti\-on', the operators of partial derivatives and differentials
and to consider decompositions of geometrical objects with respect to
adapted bases and cobases.

\subsubsection{Anholonomic frames in v--bundles}

In a v--bunde $\mathcal{E}$ provided with a N-connection we can adapt to
this structure the geometric constructions by introducing locally adapted
basis (N--frame, or N--basis):
\begin{equation}
\delta _\alpha =\frac \delta {\delta u^\alpha }=\left( \delta _i=\frac
\delta {\delta x^i}=\partial _i-N_i^{~a}\left( u\right) \partial _a,\partial
_a=\frac \partial {\partial y^a}\right) ,  \label{dder}
\end{equation}
and its dual N--basis, (N--coframe, or N--cobasis),
\begin{equation}
\delta \ ^\alpha =\delta u^\alpha =\left( d^i=\delta x^i=dx^i,\delta
^a=\delta y^a+N_i^{~a}\left( u\right) dx^i\right) .  \label{ddif}
\end{equation}

The\textbf{\ anholonomic coefficients, } $\mathbf{w}=\{w_{\beta \gamma
}^\alpha \left( u\right) \},$ of N--frames are defined to satisfy the
relations
\index{anholonomic }
\begin{equation}  \label{anhol}
\left[ \delta _\alpha ,\delta _\beta \right] =\delta _\alpha \delta _\beta
-\delta _\beta \delta _\alpha =w_{\beta \gamma }^\alpha \left( u\right)
\delta _\alpha .
\end{equation}


A frame bases is holonomic is all anholonomy coefficients vanish (like for
usual coordinate bases (\ref{pdif})), or anholonomic if there are nonzero
values of $w_{\beta \gamma }^\alpha.$

So, we conclude that a N--connection structure splitting conventionally a
v--bundle $\mathcal{E}$ into some horizontal $H\mathcal{E}$ and vertical $V%
\mathcal{E} $ subbundles can be modeled by an anholonomic frame
structure with mixed holonomic $\{x^i\}$ and anholonomic
$\{y^a\}$ variables. This case differs from usual, for instance,
tetradic approach in general relativity when tetradic (frame)
fields are stated to have only for holonomic or only for
anholonomic variables. By using the N--connection formalism we
can investigate geometrical and physical systems when some degees
of freedoms (variables) are subjected to anholonomic constraints,
the rest of variables being holonomic.

The operators (\ref{dder}) and (\ref{ddif}) on a v--bundle $\mathcal{E} $
en\-abled with a N--connecti\-on can be considered as respective equivalents
of the operators of partial derivations and differentials:\ the existence of
a N--connection structure results in 'elongation' of partial derivations on $%
x$--variables and in 'elongation' of differentials on $y$--variables.


The \textbf{algebra of tensorial distinguished fields} $DT\left( \mathcal{E}%
\right) $ (d--fields, d--ten\-sors, d--objects) on $\mathcal{E}$ is
introduced as the tensor algebra $\mathcal{T} =\{ \mathcal{T}_{qs}^{pr}\}$
of the v--bundle
\begin{equation*}
\mathcal{E}_{\left( d\right) }= \left(H\mathcal{E}\oplus V\mathcal{E}, p_d,
\mathcal{E} \right),
\end{equation*}
where $p_d:\ H\mathcal{E}\oplus V\mathcal{E}\rightarrow \mathcal{E}.$

An element $\mathbf{t}\in \mathcal{T}_{qs}^{pr},$ d--tensor field of type $%
\left(
\begin{array}{cc}
p & r \\
q & s%
\end{array}%
\right) ,$ can be written in local form as
\begin{equation}
\mathbf{t}=t_{j_{1}...j_{q}b_{1}...b_{r}}^{i_{1}...i_{p}a_{1}...a_{r}}\left(
u\right) \delta _{i_{1}}\otimes ...\otimes \delta _{i_{p}}\otimes \partial
_{a_{1}}\otimes ...\otimes \partial _{a_{r}}\otimes d^{j_{1}}\otimes
...\otimes d^{j_{q}}\otimes \delta ^{b_{1}}...\otimes \delta ^{b_{r}}.
\notag
\end{equation}

We shall respectively use the denotations $\mathcal{X\left( E\right) }$ (or $%
\mathcal{X\ } {\left( M\right) ),\ } \Lambda ^p\left( \mathcal{E}\right) $
or \newline
$\left( \Lambda ^p\left( M\right) \right) $ and $\mathcal{F\left( E\right) }$
(or $\mathcal{F}$ $\left( M\right) $) for the module of d--vector fields on $%
\mathcal{E}$ (or\newline
$M$), the exterior algebra of p--forms on $\mathcal{E}$ (or $M$) and the set
of real functi\-ons on $\mathcal{E}$ (or $M).$

\subsubsection{Anholonomic frames in cv--bundles}

The anholnomic frames
\index{anholonomic} adapted to the \v N--connection structure are introduced
similarly to (\ref{dder}) and (\ref{ddif}):

the locally adapted basis (\v N--basis, or \v N--frame):
\begin{equation}
{\breve \delta}_\alpha =
\frac {\breve \delta}{\delta u^\alpha}= \left( \delta _i=\frac \delta
{\delta x^i}= \partial _i + {\breve N}_{ia}\left({\breve u}\right) {\breve
\partial}^a, {\breve \partial}^a = \frac \partial {\partial p_a}\right) ,
\label{ddercv}
\end{equation}

and its dual (\v N--cobasis, or \v N--coframe)
\index{N--cobasis}
\index{N--coframe} :
\begin{equation}
{\breve \delta}^\alpha ={\breve \delta} u^\alpha =\left( d^i=\delta
x^i=dx^i,\ {\breve \delta}_a= {\breve \delta} p_a= d p_a - {\breve N}%
_{ia}\left({\breve u}\right) dx^i\right) .  \label{ddifcv}
\end{equation}

We note that for the singes of \v N--elongations are inverse to
those for N--elongations.

The\textbf{\ anholonomic coefficients, } $\mathbf{\breve w}= \{{\breve w}%
_{\beta \gamma }^\alpha \left({\breve u}\right) \},$ of \v N--frames are
defined by the relations
\begin{equation}  \label{anhola}
\left[{\breve \delta} _\alpha ,{\breve \delta} _\beta \right] ={\breve \delta%
}_\alpha {\breve \delta}_\beta - {\breve \delta}_\beta {\breve \delta}%
_\alpha = {\breve w}_{\beta \gamma }^\alpha \left({\breve u}\right) {\breve
\delta}_\alpha .
\end{equation}

The \textbf{algebra of tensorial distinguished fields} $DT\left({%
\breve {\mathcal{E}}}\right) $ (d--fields, d--tensors, d--objects) on ${%
\breve {\mathcal{E}}}$ is introduced as the tensor algebra ${\breve {%
\mathcal{T}}} =\{ {\breve {\mathcal{T}}}_{qs}^{pr}\}$ of the cv--bundle
\begin{equation*}
{\breve {\mathcal{E}}}_{\left( d\right) }= \left(H{\breve {\mathcal{E}}}%
\oplus V{\breve {\mathcal{E}}},{\breve p}_d, {\breve {\mathcal{E}}} \right),
\end{equation*}
where ${\breve p}_d:\ H{\breve {\mathcal{E}}}\oplus V{\breve {\mathcal{E}}}%
\rightarrow {\breve {\mathcal{E}}}.$

An element ${\breve{\mathbf{t}}}\in {\breve{\mathcal{T}}}_{qs}^{pr},$
d--tensor field of type $\left(
\begin{array}{cc}
p & r \\
q & s%
\end{array}%
\right) ,$ can be written in local form as
\begin{equation}
{\breve{\mathbf{t}}}={\breve{t}}%
_{j_{1}...j_{q}b_{1}...b_{r}}^{i_{1}...i_{p}a_{1}...a_{r}}\left( {\breve{u}}%
\right) {\breve{\delta}}_{i_{1}}\otimes ...\otimes {\breve{\delta}}%
_{i_{p}}\otimes {\breve{\partial}}_{a_{1}}\otimes ...\otimes {\breve{\partial%
}}_{a_{r}}\otimes {\breve{d}}^{j_{1}}\otimes ...\otimes {\breve{d}}%
^{j_{q}}\otimes {\breve{\delta}}^{b_{1}}...\otimes {\breve{\delta}}^{b_{r}}.
\notag
\end{equation}

We shall respectively use the denotations $\mathcal{X}\left({\breve E}%
\right) $ (or $\mathcal{X} \left( M\right) ),\ \Lambda ^p\left({\breve {%
\mathcal{E}}}\right) $ or\newline
$\left(\Lambda ^p\left( M\right) \right) $ and $\mathcal{F}\left( {\breve E}
\right)$ (or $\mathcal{F}$ $\left( M\right) $) for the module of d--vector
fields on ${\breve {\mathcal{E}}}$ (or\newline
$M$), the exterior algebra of p--forms on ${\breve {\mathcal{E}}}$\ (or $M)$
and the set of real functions on ${\breve {\mathcal{E}}}$ (or $M).$

\subsubsection{Anholonomic frames in hvc--bundles}

\index{anholonomic}

The anholnomic frames adapted to a N--connection in hvc--bundle $%
\tilde{\mathcal{E}}$ are defined by the set of coefficients (\ref{nconho});
having restricted the constructions to a vector (covector) shell we obtain
some generalizations of the formulas for corresponding N(or \v{N}%
)--connection elongation of partial derivatives defined by
(\ref{dder}) (or (\ref{ddercv})) and (\ref{ddif}) (or
(\ref{ddifcv})).

We introduce the adapted partial derivatives (anholonomic N--frames, or
N--bases) in $\tilde{\mathcal{E}}$ by applying the coefficients (\ref{nconho}%
)
\begin{equation*}
{\tilde{\delta}}_\alpha =\frac{{\tilde{\delta}}}{\delta \tilde{u}^\alpha }%
=\left( \delta _i,\delta _{a_1},\delta _{a_2},{\breve{\delta}}^{a_3},...,{%
\breve{\delta}}^{a_{z-1}},\partial _{a_z}\right) ,
\end{equation*}
where 
\begin{eqnarray*}
&&\delta _i=\partial _i-N_i^{~a_1}\partial _{a_1}-N_i^{~a_2}\partial
_{a_2}+N_{ia_3}{\breve{\partial}}^{a_3}-...+N_{ia_{z-1}}{\breve{\partial}}%
^{a_{z-1}}-N_i^{~a_z}\partial _{a_z}, \\
&&\delta _{a_1}=\partial _{a_1}-N_{a_1}^{~a_2}\partial _{a_2}+N_{a_1a_3}{%
\breve{\partial}}^{a_3}-...+N_{a_1a_{z-1}}{\breve{\partial}}%
^{a_{z-1}}-N_{a_1}^{~a_z}\partial _{a_z}, \\
&&\delta _{a_2}=\partial _{a_2}+N_{a_2a_3}{\breve{\partial}}%
^{a_3}-...+N_{a_2a_{z-1}}{\breve{\partial}}^{a_{z-1}}-N_{a_2}^{~a_z}\partial
_{a_z}, \\
&&{\breve{\delta}}^{a_3}={\tilde{\partial}}^{a_3}-N^{a_3a_4}\partial
_{a_4}-...+N_{~a_{z-1}}^{a_3}{\breve{\partial}}^{a_{z-1}}-N^{~a_3a_z}%
\partial _{a_z}, \\
&&................. \\
&&{\breve{\delta}}^{a_{z-1}}={\tilde{\partial}}^{a_{z-1}}-N^{~a_{z-1}a_z}%
\partial _{a_z}, \\
&&\partial _{a_z}=\partial /\partial y^{a_z}.
\end{eqnarray*}
These formulas can be written in the matrix form:
\begin{equation}
{\tilde{\delta}}_{_{\bullet }}=\widehat{\mathbf{N}}(u)\times {\tilde{\partial%
}}_{_{\bullet }}  \label{dderho}
\end{equation}
where
\begin{eqnarray}
{\tilde{\delta}}_{_{\bullet }} &=&\left(
\begin{array}{l}
\delta _i \\
\delta _{a_1} \\
\delta _{a_2} \\
{\breve{\delta}}^{a_3} \\
... \\
{\breve{\delta}}^{a_{z-1}} \\
\partial _{a_z}%
\end{array}
\right) ,\quad {\tilde{\partial}}_{_{\bullet }}=\left(
\begin{array}{l}
\partial _i \\
\partial _{a_1} \\
\partial _{a_2} \\
{\tilde{\partial}}^{a_3} \\
... \\
{\tilde{\partial}}^{a_{z-1}} \\
\partial _{a_z}%
\end{array}
\right) ,\quad  \label{rows}
\end{eqnarray}
and
\begin{eqnarray}
\widehat{\mathbf{N}} &=&\left(
\begin{array}{llllllll}
1 & -N_i^{~a_1} & -N_i^{~a_2} & N_{ia_3} & -N_i^{~a_4} & ... & N_{ia_{z-1}}
& -N_i^{~a_z} \\
0 & 1 & -N_{a_1}^{~a_2} & N_{a_1a_3} & -N_{a_1}^{~a_4} & ... & N_{a_1a_{z-1}}
& -N_{a_1}^{~a_z} \\
0 & 0 & 1 & N_{a_2a_3} & -N_{a_2}^{~a_4} & ... & N_{a_2a_{z-1}} &
-N_{a_2}^{~a_z} \\
0 & 0 & 0 & 1 & -N^{a_3a_4} & ... & N_{~a_{z-1}}^{a_3} & -N^{~a_3a_z} \\
... & ... & ... & ... & ... & ... & ... & ... \\
0 & 0 & 0 & 0 & 0 & ... & 1 & -N^{~a_{z-1}a_z} \\
0 & 0 & 0 & 0 & 0 & ... & 0 & 1%
\end{array}
\right) .  \notag
\end{eqnarray}

The adapted differentials (anholonomic N--coframes, or N--cobases) in $%
\tilde{\mathcal{E}}$ are introduced in the simplest form by using
matrix formalism: The respective dual matrices to (\ref{rows})
\begin{eqnarray*}
{\tilde{\delta}}^{\bullet } &=&\{{\tilde{\delta}}^\alpha \}=\left(
\begin{array}{lllllll}
d^i & \delta ^{a_1} & \delta ^{a_2} & {\breve{\delta}}_{a_3} & ... & {\breve{%
\delta}}_{a_{z-1}} & \delta ^{a_z}%
\end{array}
\right) , \\
{\tilde{d}}^{\bullet } &=&\{{\tilde{\partial}}^\alpha \}=\left(
\begin{array}{lllllll}
d^i & d^{a_1} & d^{a_2} & d_{a_3} & ... & {d}_{a_{z-1}} & d^{a_z}%
\end{array}
\right)
\end{eqnarray*}
are related via a matrix relation
\begin{equation}
{\tilde{\delta}}^{\bullet }={\tilde{d}}^{\bullet }\widehat{\mathbf{M}}
\label{ddifho}
\end{equation}
which defines the formulas for anholonomic N--coframes. The matrix $\widehat{%
\mathbf{M}}$ from (\ref{ddifho}) is the inverse to $\widehat{\mathbf{N}},$
i. e. satisfies the condition
\begin{equation}
\widehat{\mathbf{M}}\times \widehat{\mathbf{N}}=I.  \label{invmatr}
\end{equation}

The\textbf{\ anholonomic coefficients, } $\widetilde{\mathbf{w}}=\{%
\widetilde{w}_{\beta \gamma }^\alpha \left( \widetilde{u}\right) \},$ on
hcv--bundle $\tilde{\mathcal{E}}$ are expressed via coefficients of the
matrix $\widehat{\mathbf{N}}$ and their partial derivatives following the
relations
\begin{equation}
\left[ \widetilde{\delta }_\alpha ,\widetilde{\delta }_\beta \right] =%
\widetilde{\delta }_\alpha \widetilde{\delta }_\beta -\widetilde{\delta }%
_\beta \widetilde{\delta }_\alpha =\widetilde{w}_{\beta \gamma }^\alpha
\left( \widetilde{u}\right) \widetilde{\delta }_\alpha .  \label{anholho}
\end{equation}
We omit the explicit formulas on shells.

A d--tensor formalism can be also developed on the space $\tilde{\mathcal{E}}%
.$ In this case the indices have to be stipulated for every shell
separately, like for v--bundles or cv--bundles.

Let us consider some examples for particular cases of hcv--bundles:

\subsubsection{Anholonomic frames in osculator bundle}

For the osculator bundle $T^2M=Osc^2M$ from subsection (\ref{nconoscs}) the
formulas (\ref{dderho}) and (\ref{ddifho}) are written respectively in the
form
\begin{equation*}
{\tilde{\delta}}_\alpha =\left( \frac \delta {\delta x^i},\frac \delta
{\delta y_{(1)}^i},\frac \partial {\partial y_{(2)}^i}\right) ,
\end{equation*}
where
\begin{eqnarray*}
\frac \delta {\delta x^i} &=&\frac \partial {\partial x^i}-N_{(1)i}^{\qquad
j}\frac \partial {\partial y_{(1)}^i}-N_{(2)i}^{\qquad j}\frac \partial
{\partial y_{(2)}^i}, \\
\frac \delta {\delta y_{(1)}^i} &=&\frac \partial {\partial
y_{(1)}^i}-N_{(2)i}^{\qquad j}\frac \partial {\partial y_{(2)}^j},
\end{eqnarray*}
and
\begin{equation}
{\tilde{\delta}}^\alpha =\left( dx^i,\delta y_{(1)}^i,\delta
y_{(2)}^i\right) ,  \label{ddifosc2}
\end{equation}
where
\begin{eqnarray*}
\delta y_{(1)}^i &=&dy_{(1)}^i+M_{(1)j}^idx^j, \\
\delta y_{(2)}^i &=&dy_{(2)}^i+M_{(1)j}^idy_{(1)}^j+M_{(2)j}^idx^j,
\end{eqnarray*}
with the dual coefficients $M_{(1)j}^i$ and $M_{(2)j}^i$ (see (\ref{invmatr}%
)) expressed via primary coefficients $N_{(1)j}^i$ and $N_{(2)j}^i$ as
\begin{equation*}
M_{(1)j}^i=N_{(1)j}^i,M_{(2)j}^i=N_{(2)j}^{\quad i}+N_{(1)m}^{\quad
i}N_{(1)j}^{\quad m}.
\end{equation*}

\subsubsection{Anholonomic frames in dual osculator bundle}

\label{nconoscds1}

Following the definitions for dual osculator bundle $(T^{*2}M,p^{*2},M)$ in
subsection (\ref{nconoscds}) the formulas (\ref{dderho}) and (\ref{ddifho})
are written respectively in the form
\begin{equation*}
{\tilde{\delta}}_\alpha =\left( \frac \delta {\delta x^i},\frac \delta
{\delta y_{(1)}^i},\frac \partial {\partial p_{(2)}i}\right) ,
\end{equation*}
where
\begin{eqnarray*}
\frac \delta {\delta x^i} &=&\frac \partial {\partial x^i}-N_{(1)i}^{\qquad
j}\frac \partial {\partial y_{(1)}^i}+N_{(2)ij}^{\quad }\frac \partial
{\partial p_{(2)}j}, \\
\frac \delta {\delta y_{(1)}^i} &=&\frac \partial {\partial
y_{(1)}^i}+N_{(2)ij}^{\qquad }\frac \partial {\partial p_{(2)j}},
\end{eqnarray*}
and
\begin{equation}
{\tilde{\delta}}^\alpha =\left( dx^i,\delta y_{(1)}^i,\delta p_{(2)i}\right)
,  \label{ddifosc2d}
\end{equation}
where
\begin{eqnarray*}
\delta y_{(1)}^i &=&dy_{(1)}^i+N_{(1)j}^idx^j, \\
\delta p_{(2)i} &=&dp_{(2)i}-N_{(2)ij}dx^j,
\end{eqnarray*}
with the dual coefficients $M_{(1)j}^i$ and $M_{(2)j}^i$ (see (\ref{invmatr}%
)) were expressed via $N_{(1)j}^i$ and $N_{(2)j}^i$ like in Ref. \cite{mhss}.

\section{Distinguished connections and metrics}

In general, distinguished objects (d--objects)
\index{d--objects } on a v--bundle $\mathcal{E}$ (or cv--bundle ${%
\breve{\mathcal{E}}}$) are introduced as geometric objects with various
group and coordinate transforms coordinated with the N--connection structure
on $\mathcal{E}$ (or ${\breve{\mathcal{E}}}$). For example, a distinguished
connection (in brief, \textbf{d--connection)}
\index{d--connection} $D$ on $\mathcal{E}$ (or ${%
\breve{\mathcal{E}}}$) is defined as a linear connection $D$ on $E$ (or $%
\breve{E}$) conserving under a parallelism the global decomposition (\ref%
{whitney}) 
(or (\ref{whitneyc})) into horizontal and vertical subbundles of $T\mathcal{E%
}$ (or $T{\breve{\mathcal{E}}}).$ A covariant derivation
associated to a d--connection becomes d--covariant. We shall give
necessary formulas for cv--bundles in round brackets.

\subsection{D--connections}

\subsubsection{D--connections in v--bundles (cv--bundles)}

A N--connection in a v--bundle $\mathcal{E}$ (cv--bundle ${\breve{\mathcal{E}%
}}$) induces a corresponding decomposition of d--tensors into sums of
horizontal and vertical parts, for example, for every d--vector $X\in
\mathcal{X\left( E\right) }$ (${\breve X}\in \mathcal{X}\left( {\breve{%
\mathcal{E}}}\right) $ ) and 1--form $A\in \Lambda ^1\left( \mathcal{E}%
\right) $ ($\breve{A}\in \Lambda ^1\left( {\breve{\mathcal{E}}}\right) $) we
have respectively
\begin{eqnarray}
X &=&hX+vX\mathbf{\ \quad }\mbox{and \quad }A=hA+vA,  \label{vdecomp} \\
(\breve{X} &=&h\breve{X}+vX\mathbf{\ \quad }\mbox{and \quad }\breve{A}=h%
\breve{A}+v\breve{A})  \notag
\end{eqnarray}
where
\begin{equation*}
hX = X^i\delta _i,vX=X^a\partial _a \ (h\breve{X} = \breve{X}^i{\tilde \delta%
}_i,v\breve{X}=\breve{X}_a{\breve \partial}^a)
\end{equation*}
and
\begin{equation*}
hA =A_i\delta ^i,vA=A_ad^a \ (h\breve{A} = \breve{A}_i{\breve \delta}^i,v%
\breve{A}= \breve{A}^a{\breve d}_a).
\end{equation*}

In consequence, we can associate to every d--covariant derivation along the
d--vector (\ref{vdecomp}), 
$D_X=X\circ D$ ($D_{\breve{X}}=\breve{X}\circ D$) two new operators of h-
and v--covariant derivations
\begin{eqnarray*}
D_X^{(h)}Y &=&D_{hX}Y\quad \mbox{ and \quad }D_X^{\left( v\right)
}Y=D_{vX}Y,\quad \forall Y{\ \in }\mathcal{X\left( E\right) } \\
(D_{\breve{X}}^{(h)}\breve{Y} &=&D_{h\breve{X}}\breve{Y}\quad
\mbox{ and
\quad }D_{\breve{X}}^{\left( v\right) }\breve{Y}=D_{v\breve{X}}\breve{Y}%
,\quad \forall \breve{Y}{\ \in }\mathcal{X}\left( {\breve{\mathcal{E}}}%
\right) )
\end{eqnarray*}
for which the following conditions hold:
\begin{eqnarray}
D_XY &=&D_X^{(h)}Y{\ +}D_X^{(v)}Y  \label{dcovpr} \\
(D_{\breve{X}}\breve{Y} &=&D_{\breve{X}}^{(h)}\breve{Y}{\ +}D_{\breve{X}%
}^{(v)}\breve{Y}),  \notag
\end{eqnarray}
where
\begin{eqnarray*}
D_X^{(h)}f &=&(hX\mathbf{)}f\mbox{ \quad and\quad }D_X^{(v)}f=(vX\mathbf{)}%
f,\quad X,Y\mathbf{\in }\mathcal{X\left( E\right) },f\in \mathcal{F}\left(
M\right) \\
({\breve{D}}_{\breve{X}}^{(h)}f &=&(h\breve{X}\mathbf{)}f%
\mbox{ \quad
and\quad }{\breve{D}}_{\breve{X}}^{(v)}f=(v\breve{X}\mathbf{)}f,\quad \breve{%
X},\breve{Y}{\in }\mathcal{X}\left( {\breve{\mathcal{E}}}\right) ,f\in
\mathcal{F}\left( M\right) ).
\end{eqnarray*}

The components $\Gamma _{\beta \gamma }^\alpha $ ( ${\breve{\Gamma}}_{\beta
\gamma }^\alpha )$of a d--connection ${\breve{D}}_\alpha =({\breve{\delta}}%
_\alpha \circ D),$ locally adapted to the N---connection structure with
respect to the frames (\ref{dder}) and (\ref{ddif}) ((\ref{ddercv}) and (\ref%
{ddifcv})), are defined by the equations
\begin{equation*}
D_\alpha \delta _\beta =\Gamma _{\alpha \beta }^\gamma \delta _\gamma ~({%
\breve{D}}_\alpha {\breve{\delta}}_\beta ={\breve{\Gamma}}_{\alpha \beta
}^\gamma {\breve{\delta}}_\gamma ~),
\end{equation*}
from which one immediately follows
\begin{equation}
\Gamma _{\alpha \beta }^\gamma \left( u\right) =\left( D_\alpha \delta
_\beta \right) \circ \delta ^\gamma \quad ~({\breve{\Gamma}}_{\alpha \beta
}^\gamma \left( {\breve{u}}\right) =\left( {\breve{D}}_\alpha {\breve{\delta}%
}_\beta \right) \circ {\breve{\delta}}^\gamma ).  \label{gamma}
\end{equation}

The coefficients of operators of h- and v--covariant derivations,
\begin{eqnarray*}
D_k^{(h)} &=&\{L_{jk}^i,L^a_{bk\;}\}\mbox{ and }D_c^{(v)}=%
\{C_{jk}^i,C_{bc}^a\} \\
({\breve{D}}_k^{(h)} &=&\{{\breve{L}}_{jk}^i,{\breve{L}}_{a k}^{~b}\}%
\mbox{
and }{\breve{D}}^{(v)c}=\{{\breve C}_{~j}^{i~c},{\breve C}_a^{~bc}\})
\end{eqnarray*}
(see (\ref{dcovpr})), are introduced as corresponding h- and
v--paramet\-ri\-za\-ti\-ons of (\ref{gamma})
\begin{eqnarray}
L_{jk}^i &=&\left( D_k\delta _j\right) \circ d^i,\quad L_{bk}^a=\left(
D_k\partial _b\right) \circ \delta ^a  \label{hgamma} \\
({\breve{L}}_{jk}^i &=&\left( {\breve{D}}_k{\breve{\delta}}_j\right) \circ
d^i,\quad {\breve{L}}_{a k}^{~b}=\left( {\breve{D}}_k{\breve \partial}%
^b\right) \circ {\breve{\delta}}_a)  \notag
\end{eqnarray}
and
\begin{eqnarray}
C_{jc}^i &=&\left( D_c\delta _j\right) \circ d^i,\quad C_{bc}^a=\left(
D_c\partial _b\right) \circ \delta ^a  \label{vgamma} \\
({\breve C}_{~j}^{i~c} &=&\left( {\breve{D}}^c{\breve{\delta}}_j\right)
\circ d^i,\quad {\breve C}_a^{~bc}=\left( {\breve{D}}^c{\breve \partial}%
^b\right) \circ {\breve{\delta}}_a).  \notag
\end{eqnarray}

A set of components (\ref{hgamma}) and (\ref{vgamma}) \
\begin{equation*}
\Gamma _{\alpha \beta }^\gamma =[L_{jk}^i,L_{bk}^a,C_{jc}^i,C_{bc}^a]~\left(
{\breve{\Gamma}}_{\alpha \beta }^\gamma =[{\breve{L}}_{jk}^i,{\breve{L}}%
_{ak}^{~b},{\breve{C}}_{~j}^{i~c},{\breve{C}}_a^{~bc}]\right)
\end{equation*}
completely defines the local action of a d---connection $D$ in $\mathcal{E}$
(${\breve{D}}$ in ${\breve{\mathcal{E}}).}$

For instance, having taken on $\mathcal{E}$ (${\breve{\mathcal{E}})}$ a
d---tensor field of type $\left(
\begin{array}{cc}
1 & 1 \\
1 & 1%
\end{array}
\right) ,$
\begin{eqnarray*}
\mathbf{t} &=&t_{jb}^{ia}\delta _i\otimes \partial _a\otimes d^j\otimes
\delta ^b, \\
{{\tilde{\mathbf{t}}}} &=&{\breve{t}}_{ja}^{ib}{\breve{\delta}}_i\otimes {%
\breve{\partial}}^a\otimes d^j\otimes {\breve{\delta}}_b,
\end{eqnarray*}
and a d--vector $\mathbf{X}$ ($\mathbf{{\breve{X}}})$ we obtain
\begin{eqnarray*}
D_X\mathbf{t} &=&D_X^{(h)}\mathbf{t+}D_X^{(v)}\mathbf{t=}\left( X^k{\breve{t}%
}_{jb|k}^{ia}+X^ct_{jb\perp c}^{ia}\right) \delta _i\otimes \partial
_a\otimes d^j\otimes \delta ^b, \\
({\breve{D}}_{\breve{X}}{{\tilde{\mathbf{t}}}} &=&\breve{D}_{\breve{X}}^{(h)}%
{\tilde{\mathbf{t}}}+\breve{D}_{\breve{X}}^{(v)}{\tilde{\mathbf{t}}}=\left(
\breve{X}^k{\breve{t}}_{ja|k}^{ib}+\breve{X}_c{\breve{t}}_{ja}^{ib\perp
c}\right) {\breve{\delta}}_i\otimes {\breve{\partial}}^a\otimes d^j\otimes {%
\breve{\delta}}_b)
\end{eqnarray*}
where the h--covariant derivative is written
\begin{eqnarray*}
t_{jb|k}^{ia} &=&\delta
_kt_{jb}^{ia}+L_{hk}^it_{jb}^{ha}+L_{ck}^at_{jb}^{ic}-L_{jk}^ht_{hb}^{ia}-L_{bk}^ct_{jc}^{ia}
\\
({\breve{t}}_{ja|k}^{ib} &=&{\breve{\delta}}_k{\breve{t}}_{ja}^{ib}+{\breve{L%
}}_{hk}^i{\breve{t}}_{ja}^{hb}+{\breve{L}}_{ck}^{~b}{\breve{t}}_{ja}^{ic}-{%
\breve{L}}_{jk}^h{\breve{t}}_{ha}^{ib}-{\breve{L}}_{ck}^{~b}{\breve{t}}%
_{ja}^{ic})
\end{eqnarray*}
and the v-covariant derivative is written
\begin{eqnarray}
t_{jb\perp c}^{ia} &=&\partial
_ct_{jb}^{ia}+C_{hc}^it_{jb}^{ha}+C_{dc}^at_{jb}^{id}-C_{jc}^ht_{hb}^{ia}-C_{bc}^dt_{jd}^{ia}
\label{covder1} \\
({\breve{t}}_{ja}^{ib\perp c} &=&{\breve{\partial}}^c{\breve{t}}_{ja}^{ib}+{%
\breve{C}}_{~j}^{i~c}{\breve{t}}_{ja}^{hb}+{\breve{C}}_a^{~dc}{\breve{t}}%
_{jd}^{ib}-{\breve{C}}_{~j}^{i~c}{\breve{t}}_{ha}^{ib}-{\breve{C}}_d^{~bc}{%
\breve{t}}_{ja}^{id}).  \label{covder2}
\end{eqnarray}
For a scalar function $f\in \mathcal{F}\left(\mathcal{E}\right) $ ( $f\in
\mathcal{F}\left( {\breve{\mathcal{E}}}\right) $) we have
\begin{eqnarray*}
D_k^{(h)} &=&\frac{\delta f}{\delta x^k}=\frac{\partial f}{\partial x^k}%
-N_k^a\frac{\partial f}{\partial y^a}\mbox{ and }D_c^{(v)}f=\frac{\partial f%
}{\partial y^c} \\
(\breve{D}_k^{(h)} &=&\frac{{\breve{\delta}}f}{\delta x^k}=\frac{\partial f}{%
\partial x^k}+N_{ka}\frac{\partial f}{\partial p_a}\mbox{ and }\breve{D}%
^{(v)c}f=\frac{\partial f}{\partial p_c}).
\end{eqnarray*}

\subsubsection{D--connections in hvc--bundles}

\index{d--connection} The theory of connections in higher order anisotropic
vector superbundles and vector bundles was elaborated in Refs. \cite%
{vstr2,vhsp,vbook}. Here we re--formulate that formalism for the case when
some shells of higher order anisotropy could be covector spaces by stating
the general rules of covariant derivation compatible with the N--connection
structure in hvc--bundle $%
\tilde{\mathcal{E}}$ and omit details and cumbersome formulas.

For a hvc--bundle of type ${\tilde{\mathcal{E}}}={\tilde{\mathcal{E}}}%
[v(1),v(2),cv(3),...,cv(z-1),v(z)]$ a d--connection ${\tilde
\Gamma}_{\alpha \beta }^\gamma $ has the next shell decomposition
of components (on induction being on the $p$-th shell, considered
as the base space, which in this case a hvc--bundle, we introduce
in a usual manner, like a vector or covector fiber, the $(p+1)$-th
shell)
\begin{eqnarray*}
{\tilde \Gamma}_{\alpha \beta }^\gamma &=&\{\Gamma _{\alpha _1\beta
_1}^{\gamma
_1}=[L_{j_1k_1}^{i_1},L_{b_1k_1}^{a_1},C_{j_1c_1}^{i_1},C_{b_1c_1}^{a_1}], \\
& & \Gamma _{\alpha _2\beta _2}^{\gamma _2} =
[L_{j_2k_2}^{i_2},L_{b_2k_2}^{a_2},C_{j_2c_2}^{i_2},C_{b_2c_2}^{a_2}], \\
& & {\breve{\Gamma}}_{\alpha _3\beta _3}^{\gamma _3} =[{\breve{L}}%
_{j_3k_3}^{i_3},{\breve{L}}_{a_3k_3}^{~b_3},{\breve{C}}_{~j_3}^{i_3~c_3},{%
\breve{C}}_{a_3}^{~b_3c_3}], \\
&& ...................................., \\
& & {\breve{\Gamma}}_{\alpha _{z-1}\beta _{z-1}}^{\gamma _{z-1}}= [{\breve{L}%
}_{j_{z-1}k_{z-1}}^{i_{z-1}},{\breve{L}}_{a_{z-1}k_{z-1}}^{~b_{z-1}},{\breve{%
C}}_{~j_{z-1}}^{i_{z-1}~c_{z-1}},{\breve{C}}_{a_{z-1}}^{~b_{z-1}c_{z-1}}], \\
& & \Gamma _{\alpha _z\beta _z}^{\gamma _z}
=[L_{j_zk_z}^{i_z},L_{b_zk_z}^{a_z},C_{j_zc_z}^{i_z},C_{b_zc_z}^{a_z}]\}.
\end{eqnarray*}
These coefficients determine the rules of a covariant derivation $\tilde D$
on ${\tilde {\mathcal{E}}}.$

For example, let us consider a d--tensor ${\tilde{\mathbf{t}}}$ of type
\begin{equation*}
\left(
\begin{array}{llllll}
1 & 1_{1} & 1_{2} & {\breve{1}}_{3} & ... & 1_{z} \\
1 & 1_{1} & 1_{2} & {\breve{1}}_{3} & ... & 1_{z}%
\end{array}%
\right)
\end{equation*}%
with corresponding tensor product of components of anholonomic N--frames (%
\ref{dderho}) and (\ref{ddifho})
\begin{equation*}
{\tilde{\mathbf{t}}}={\tilde{t}}_{jb_{1}b_{2}{\breve{a}}_{3}...{\breve{a}}%
_{z-1}b_{z}}^{ia_{1}a_{2}{\breve{b}}_{3}...{\breve{b}}_{z-1}a_{z}}\delta
_{i}\otimes \partial _{a_{1}}\otimes d^{j}\otimes \delta ^{b_{1}}\otimes
\partial _{a_{2}}\otimes \delta ^{b_{2}}\otimes {\breve{\partial}}%
^{a_{3}}\otimes {\breve{\delta}}_{b_{3}}{...}\otimes {\breve{\partial}}%
^{a_{z-1}}\otimes {\breve{\delta}}_{bz-1}\otimes \partial _{a_{z}}\otimes
\delta ^{b_{z}}.
\end{equation*}%
The d--covariant derivation $\tilde{D}$ of ${\tilde{\mathbf{t}}}$ is to be
performed separately for every shall according the rule (\ref{covder1}) if a
shell is defined by a vector subspace, or according the rule (\ref{covder2})
if the shell is defined by a covector subspace.

\subsection{Metric structure}

\subsubsection{D--metrics in v--bundles}

We define a \textbf{metric structure }$\mathbf{G\ }$
\index{metric} in the total space $E$ of a v--bundle $\mathcal{E=}$ $\left(
E,p,M\right) $ over a connected and paracompact base $M$ as a symmetric
covariant tensor field of type $\left( 0,2\right) $,
\begin{equation*}
\mathbf{G} = G_{\alpha \beta } du^{\alpha}\otimes du^\beta
\end{equation*}
being non degenerate and of constant signature on $E.$

Nonlinear connection $\mathbf{N}$ and metric $\mathbf{G}$ structures on $%
\mathcal{E}$ are mutually compatible it there are satisfied the conditions:
\begin{equation}  \label{comp}
\mathbf{G}\left( \delta _i,\partial _a\right) =0,\mbox{or equivalently, }%
G_{ia}\left( u\right) -N_i^b\left( u\right) h_{ab}\left( u\right) =0,
\end{equation}
where $h_{ab}=\mathbf{G}\left( \partial _a,\partial _b\right) $ and $G_{ia}=%
\mathbf{G}\left( \partial _i,\partial _a\right),$ which gives
\begin{equation}  \label{ncon1}
N_i^b\left( u\right) = h^{ab}\left( u\right) G_{ia}\left( u\right)
\end{equation}
( the matrix $h^{ab}$ is inverse to $h_{ab}).$ In consequence one obtains
the following decomposition of metric:
\begin{equation}  \label{metrdec}
\mathbf{G}(X,Y)\mathbf{=hG}(X,Y)+\mathbf{vG}(X,Y),
\end{equation}
where the d--tensor $\mathbf{hG}(X,Y)$ = $ \mathbf{ G}(hX,hY)$ is
of type $\left(
\begin{array}{cc}
0 & 0 \\
2 & 0%
\end{array}
\right) $ and the d--tensor $\mathbf{vG}(X,Y) = \mathbf{G}(vX,vY)$ is of
type $\left(
\begin{array}{cc}
0 & 0 \\
0 & 2%
\end{array}
\right) .$ With respect to anholonomic basis (\ref{dder}) the d--metric (\ref%
{metrdec}) is written
\begin{equation}  \label{dmetric}
\mathbf{G}=g_{\alpha \beta }\left( u\right) \delta ^\alpha \otimes \delta
^\beta =g_{ij}\left( u\right) d^i\otimes d^j+h_{ab}\left( u\right) \delta
^a\otimes \delta ^b,
\end{equation}
where $g_{ij}=\mathbf{G}\left( \delta _i,\delta _j\right) .$

A metric structure of type (\ref{metrdec}) (equivalently, of type (\ref%
{dmetric})) or a metric on $E$ with components satisfying the constraints (\ref%
{comp}), (equivalently (\ref{ncon1})) defines an adapted to the given
N--connection inner (d--scalar) product on the tangent bundle $\mathcal{TE}$.

We shall say that a d--connection $%
\widehat{D}_X$ is compatible with the d-scalar product on $\mathcal{TE\ }$
(i. e. it is a standard d--connection) if
\begin{equation*}
\widehat{D}_X\left( \mathbf{X\cdot Y}\right) =\left( \widehat{D}_X\mathbf{Y}%
\right) \cdot \mathbf{Z+Y\cdot }\left( \widehat{D}_X\mathbf{Z}\right)
,\forall \mathbf{X,Y,Z}\mathbf{\in }\mathcal{X\left( E\right) }.
\end{equation*}
An arbitrary d--connection $D_X$ differs from the standard one $\widehat{D}%
_X $ by an operator $\widehat{P}_X\left( u\right) =\{X^\alpha \widehat{P}%
_{\alpha \beta }^\gamma \left( u\right) \},$ called the deformation d-tensor
with respect to $\widehat{D}_X,$ which is just a d-linear transform of $%
\mathcal{E}_u,$\ $\forall \ u\in \mathcal{E}$. The explicit form of $%
\widehat{P}_X $ can be found by using the corresponding axiom defining
linear connections \cite{lue}
\begin{equation*}
\left( D_X-\widehat{D}_X\right) fZ=f\left( D_X-\widehat{D}_X\right) Z\mathbf{%
,}
\end{equation*}
written with respect to N--elongated bases (\ref{dder}) and (\ref{ddif}).
From the last expression we obtain
\begin{equation*}
\widehat{P}_X\left( u\right) =\left[ (D_X-\widehat{D}_X)\delta _\alpha
\left( u\right) \right] \delta ^\alpha \left( u\right) ,
\end{equation*}
therefore
\begin{equation}
D_XZ\mathbf{\ }=\widehat{D}_XZ\mathbf{\ +}\widehat{P}_XZ.  \label{deft}
\end{equation}

A d--connection $D_X$ is \textbf{metric} (or \textbf{compatible } with
metric $\mathbf{G}$) on $\mathcal{E}$ if
\begin{equation*}
D_X\mathbf{G} =0,\forall X\mathbf{\in }\mathcal{X\left( E\right) }.
\end{equation*}
With respect to anholonomic frames these conditions are written
\begin{equation}  \label{comatib}
D_\alpha g_{\beta\gamma}=0,
\end{equation}
where by $g_{\beta\gamma}$ we denote the coefficients in the block form (\ref%
{dmetric}).

\subsubsection{D--metrics in cv-- and hvc--bundles}

The presented considerations on self--consistent definition of
N--connection, d--connection and metric structures in v--bundles
can re--formulated in a similar fashion for another types of
anisotropic space--times, on cv--bundles and on shells of
hvc--bundles. For simplicity, we give here only the analogous
formulas for the metric d--tensor (\ref{dmetric}):

\begin{itemize}
\item On cv--bundle ${\breve{\mathcal{E}}}$ we write
\begin{equation}
{\breve{\mathbf{G}}}={\breve{g}}_{\alpha \beta }\left( {\breve{u}}\right) {%
\breve{\delta}}^\alpha \otimes {\breve{\delta}}^\beta ={\breve{g}}%
_{ij}\left( {\breve{u}}\right) d^i\otimes d^j+{\breve{h}}^{ab}\left( {%
\breve{u}}\right) {\breve{\delta}}_a\otimes {\breve{\delta}}_b,
\label{dmetricvc}
\end{equation}
where ${\breve{g}}_{ij}={\breve{\mathbf{G}}}\left( {\breve{\delta}}_i,{%
\breve{\delta}}_j\right) $ and ${\breve{h}}^{ab}={\breve{\mathbf{G}}}\left( {%
\breve{\partial}}^a,{\breve{\partial}}^b\right) $ and the N--coframes are
given by formulas (\ref{ddifcv}).

For simplicity, we shall consider that the metricity conditions are
satisfied, ${\breve{D}}_\gamma {\breve{g}}_{\alpha \beta}=0.$

\item On hvc--bundle ${\tilde{\mathcal{E}}}$ we write
\begin{eqnarray}
 & &   \label{dmetrichcv} \\
{\tilde{\mathbf{G}}} &=&{\tilde{g}}_{\alpha \beta }\left( {\tilde{u}}\right)
{\tilde{\delta}}^{\alpha }\otimes {\tilde{\delta}}^{\beta }={\tilde{g}}%
_{ij}\left( {\tilde{u}}\right) d^{i}\otimes d^{j}+{\tilde{h}}%
_{a_{1}b_{1}}\left( {\tilde{u}}\right) {\delta }^{a_{1}}\otimes {\delta }%
^{b_{1}}+{\tilde{h}}_{a_{2}b_{2}}\left( {\tilde{u}}\right) {\delta }%
^{a_{2}}\otimes {\delta }^{b_{2}} \notag \\
&&+{\tilde{h}}^{a_{3}b_{3}}\left( {\tilde{u}}\right) {\breve{\delta}}%
_{a_{3}}\otimes {\breve{\delta}}_{b_{3}}+...+{\tilde{h}}^{a_{z-1}b_{z-1}}%
\left( {\tilde{u}}\right) {\breve{\delta}}_{a_{z-1}}\otimes {\breve{\delta}}%
_{b_{z-1}}+{\tilde{h}}_{a_{z}b_{z}}\left( {\tilde{u}}\right) {\delta }%
^{a_{z}}\otimes {\delta }^{b_{z}},  \notag
\end{eqnarray}%
where ${\tilde{g}}_{ij}={\tilde{\mathbf{G}}}\left( {\tilde{\delta}}_{i},{%
\tilde{\delta}}_{j}\right) $ and ${\tilde{h}}_{a_{1}b_{1}}={\tilde{\mathbf{G}%
}}\left( \partial _{a_{1}},\partial _{b_{1}}\right) ,$ ${\tilde{h}}%
_{a_{2}b_{2}}={\tilde{\mathbf{G}}}\left( \partial _{a_{2}},\partial
_{b_{2}}\right) ,$ ${\tilde{h}}^{a_{3}b_{3}}={\tilde{\mathbf{G}}}\left( {%
\breve{\partial}}^{a_{3}},{\breve{\partial}}^{b_{3}}\right) ,....$ and the
N--coframes are given by formulas (\ref{ddifho}).

The metricity conditions are ${\tilde{D}}_{\gamma }{\tilde{g}}_{\alpha \beta
}=0.$

\item On osculator bundle $T^{2}M=Osc^{2}M$ we have a particular case of (%
\ref{dmetrichcv}) when
\begin{eqnarray}
{\tilde{\mathbf{G}}} &=&{\tilde{g}}_{\alpha \beta }\left( {\tilde{u}}\right)
{\tilde{\delta}}^{\alpha }\otimes {\tilde{\delta}}^{\beta }
\label{dmetrichosc2} \\
&=&{\tilde{g}}_{ij}\left( {\tilde{u}}\right) d^{i}\otimes d^{j}+{\tilde{h}}%
_{ij}\left( {\tilde{u}}\right) {\delta y}_{(1)}^{i}\otimes {\delta y}%
_{(1)}^{i}+{\tilde{h}}_{ij}\left( {\tilde{u}}\right) {\delta y}%
_{(2)}^{i}\otimes {\delta y}_{(2)}^{i}  \notag
\end{eqnarray}%
where the N--coframes are given by (\ref{ddifosc2}).

\item On dual osculator bundle $\left( T^{*2}M,p^{*2},M\right) $ we have
another particular case of (\ref{dmetrichcv}) when
\begin{eqnarray}  \label{dmetrichosc2d}
{\tilde{\mathbf{G}}} &=& {\tilde{g}}_{\alpha \beta }\left( {\tilde{u}}%
\right) {\tilde{\delta}}^\alpha \otimes {\tilde{\delta}}^\beta \\
&=& {\tilde{g}}_{ij}\left( {\tilde{u}}\right) d^i\otimes d^j+{\tilde{h}}%
_{ij}\left( {\tilde{u}}\right) {\delta y}_{(1)}^i\otimes {\delta y}_{(1)}^i+{%
\tilde{h}}^{ij}\left( {\tilde{u}}\right) {\delta p}_i^{(2)}\otimes {\delta p}%
_i^{(2)}  \notag
\end{eqnarray}
where the N--coframes are given by (\ref{ddifosc2d}).
\end{itemize}


\subsection{Some remarkable d--connections}

We emphasize that the geometry of connections in a v--bundle $\mathcal{E}$
is very reach. If a triple of fundamental geometric objects $\left(
N_i^a\left( u\right) ,\Gamma _{\beta \gamma }^\alpha \left( u\right)
,g_{\alpha \beta }\left( u\right) \right) $ is fixed on $\mathcal{E}$, a
multi--connection structure (with corresponding different rules of covariant
derivation, which are, or not, mutually compatible and with the same, or
not, induced d--scalar products in $\mathcal{TE)}$ is defined on this
v--bundle. We can give a priority to a connection structure following some
physical arguments, like the reduction to the Christoffel symbols in the
holonomic case, mutual compatibility between metric and N--connection and
d--connection structures and so on.

In this subsection we enumerate some of the connections and covariant
derivations in v--bundle $\mathcal{E}$, cv--bundle ${\breve{\mathcal{E}}}$
and in some hvc--bundles which can present interest in investigation of
locally anisotropic gravitational and matter field interactions :

\begin{enumerate}
\item Every N--connection in $\mathcal{E}$ with coefficients $N_i^a\left(
x,y\right) $ being differentiable on y--variables, induces a structure of
linear connection $N_{\beta \gamma }^\alpha ,$ where
\begin{equation}
N_{bi}^a=\frac{\partial N_i^a}{\partial y^b}\mbox{ and
}N_{bc}^a\left( x,y\right) =0.  \label{nlinearized}
\end{equation}
For some $Y\left( u\right) =Y^i\left( u\right) \partial _i+Y^a\left(
u\right) \partial _a$ and $B\left( u\right) =B^a\left( u\right) \partial _a$
one introduces a covariant derivation as
\begin{equation*}
D_Y^{(\widetilde{N})}B=\left[ Y^i\left( \frac{\partial B^a}{\partial x^i}%
+N_{bi}^aB^b\right) +Y^b\frac{\partial B^a}{\partial y^b}\right] \frac
\partial {\partial y^a}.
\end{equation*}

\item The d--connection of Berwald type \cite{berw} on v--bundle $\mathcal{E}
$ (cv--bundle ${\breve{\mathcal{E}})}$
\begin{eqnarray}
\Gamma _{\beta \gamma }^{(B)\alpha } &=&\left( L_{jk}^i,\frac{\partial N_k^a%
}{\partial y^b},0,C_{bc}^a\right) ,  \label{berwald} \\
({\breve{\Gamma}}_{\beta \gamma }^{(B)\alpha } &=&\left( \breve{L}_{jk}^i,-%
\frac{\partial \breve{N}_{ka}}{\partial p_b},0,{\breve{C}}_a^{~bc}\right) )
\notag
\end{eqnarray}
where 
\begin{eqnarray}
L_{.jk}^i\left( x,y\right) &=&\frac 12g^{ir}\left( \frac{\delta g_{jk}}{%
\delta x^k}+\frac{\delta g_{kr}}{\delta x^j}-\frac{\delta g_{jk}}{\delta x^r}%
\right) ,  \label{lccoef} \\
C_{.bc}^a\left( x,y\right) &=&\frac 12h^{ad}\left( \frac{\partial h_{bd}}{%
\partial y^c}+\frac{\partial h_{cd}}{\partial y^b}-\frac{\partial h_{bc}}{%
\partial y^d}\right)  \notag \\
(\breve{L}_{.jk}^i\left( x,p\right) &=&\frac 12\breve{g}^{ir}\left( \frac{{%
\breve{\delta}}\breve{g}_{jk}}{\delta x^k}+\frac{{\breve{\delta}}\breve{g}%
_{kr}}{\delta x^j}-\frac{{\breve{\delta}}\breve{g}_{jk}}{\delta x^r}\right) ,
\notag \\
{\breve{C}}_a^{~bc}\left( x,p\right) &=&\frac 12\breve{h}_{ad}\left( \frac{%
\partial \breve{h}^{bd}}{\partial p_c}+\frac{\partial \breve{h}^{cd}}{%
\partial p_b}-\frac{\partial \breve{h}^{bc}}{\partial p_d}\right) ),  \notag
\end{eqnarray}

which is hv---metric, i.e. there are satisfied the conditions $%
D_k^{(B)}g_{ij}=0$ and $D_c^{(B)}h_{ab}=0$ ($\breve{D}_k^{(B)}\breve{g}%
_{ij}=0$ and $\breve{D}^{(B)c}\breve{h}^{ab}=0).$

\item The canonical d--connection $\mathbf{\Gamma ^{(c)}}$ (or $\mathbf{%
\breve{\Gamma}^{(c)})}$ on a v--bundle (or cv--bundle) is associated to a
metric $\mathbf{G}$ (or $\mathbf{\breve{G})}$ of type (\ref{dmetric}) (or (%
\ref{dmetricvc})),
\begin{equation*}
\Gamma _{\beta \gamma }^{(c)\alpha
}=[L_{jk}^{(c)i},L_{bk}^{(c)a},C_{jc}^{(c)i},C_{bc}^{(c)a}]~(\breve{\Gamma}%
_{\beta \gamma }^{(c)\alpha }=[\breve{L}_{jk}^{(c)i},\breve{L}%
_{~a~.k}^{(c).b},\breve{C}_{~j}^{(c)i\ c},{\breve{C}}_a^{(c)~bc}])
\end{equation*}
with coefficients
\begin{eqnarray}
L_{jk}^{(c)i} &=&L_{.jk}^i,C_{bc}^{(c)a}=C_{.bc}^a~(\breve{L}_{jk}^{(c)i}=%
\breve{L}_{.jk}^i,{\breve{C}}_a^{(c)~bc}={\breve{C}}_a^{~bc}),%
\mbox{ (see
(\ref{lccoef})}  \notag \\
L_{bi}^{(c)a} &=&\frac{\partial N_i^a}{\partial y^b}+\frac 12h^{ac}\left(
\frac{\delta h_{bc}}{\delta x^i}-\frac{\partial N_i^d}{\partial y^b}h_{dc}-%
\frac{\partial N_i^d}{\partial y^c}h_{db}\right)  \notag \\
~(\breve{L}_{~a~.i}^{(c).b} &=&-\frac{\partial {\breve{N}}_i^a}{\partial p_b}%
+\frac 12\breve{h}_{ac}\left( \frac{{\breve{\delta}}\breve{h}^{bc}}{\delta
x^i}+\frac{\partial {\breve{N}}_{id}}{\partial p_b}\breve{h}^{dc}+\frac{%
\partial {\breve{N}}_{id}}{\partial p_c}\breve{h}^{db}\right) ),  \notag \\
~C_{jc}^{(c)i} &=&\frac 12g^{ik}\frac{\partial g_{jk}}{\partial y^c}~(\breve{%
C}_{~j}^{(c)i\ c}=\frac 12\breve{g}^{ik}\frac{\partial \breve{g}_{jk}}{%
\partial p_c}).  \label{inters}
\end{eqnarray}
This is a metric d--connection which satisfies conditions
\begin{eqnarray*}
D_k^{(c)}g_{ij} &=&0,D_c^{(c)}g_{ij}=0,D_k^{(c)}h_{ab}=0,D_c^{(c)}h_{ab}=0 \\
(\breve{D}_k^{(c)}\breve{g}_{jk} &=&0,\breve{D}^{(c)c}\breve{g}_{jk}=0,%
\breve{D}_k^{(c)}\breve{h}^{bc}=0,\breve{D}^{(c)c}\breve{h}^{ab}=0).
\end{eqnarray*}
In physical applications we shall use the canonical connection and for
simplicity we shall omit the index $(c).$ The coefficients (\ref{inters})$\,$%
are to be extended to higher order if we are dealing with derivations of
geometrical objects with ''shell'' indices. In this case the fiber indices
are to be stipulated for every type of shell into consideration.

\item We can consider the N--adapted Christoffel d--symbols
\begin{equation}
\widetilde{\Gamma }_{\beta \gamma }^\alpha =\frac 12g^{\alpha \tau }\left(
\delta _\gamma g_{\tau \beta }+\delta _\beta g_{\tau \gamma }-\delta
g_{\beta \gamma }\right) ,  \label{dchrist}
\end{equation}
which have the components of d--connection $\widetilde{\Gamma }_{\beta
\gamma }^\alpha =\left( L_{jk}^i,0,0,C_{bc}^a\right) ,$ with $L_{jk}^i$ and $%
C_{bc}^a$ as in (\ref{lccoef}) 
if $g_{\alpha \beta }$ is taken in the form (\ref{dmetric}). 
\end{enumerate}

Arbitrary linear connections on a v-bundle $\mathcal{E}$ can be also
characterized by theirs deformation tensors (see (\ref{deft})) 
with respect, for instance, to the d--connect\-i\-on (\ref{dchrist}):
\begin{equation*}
\Gamma _{\beta \gamma }^{(B)\alpha }=\widetilde{\Gamma }_{\beta \gamma
}^\alpha +P_{\beta \gamma }^{(B)\alpha },\Gamma _{\beta \gamma }^{(c)\alpha
}=\widetilde{\Gamma }_{\beta \gamma }^\alpha +P_{\beta \gamma }^{(c)\alpha }
\end{equation*}
or, in general,
\begin{equation*}
\Gamma _{\beta \gamma }^\alpha =\widetilde{\Gamma }_{\beta \gamma }^\alpha
+P_{\beta \gamma }^\alpha ,
\end{equation*}
where $P_{\beta \gamma }^{(B)\alpha },P_{\beta \gamma }^{(c)\alpha }$ and $%
P_{\beta \gamma }^\alpha $ are respectively the deformation d-tensors of
d--connect\-i\-ons (\ref{berwald}),\ (\ref{inters}) 
or of a general one. Similar deformation d--tensors can be introduced for
d--connections on cv--bundles and hvc--bundles. We omit explicit formulas.

\subsection{Amost Hermitian anisotropic spaces}

\index{almost Hermitian}

The are possible very interesting particular constructions \cite%
{ma87,ma94,mhss} on t--bundle $TM$ provided with N--connection which defines
a N--adapted frame structure $\delta _{\alpha }=(\delta _{i},%
\dot{\partial}_{i})$ (for the same formulas (\ref{dder}) and (\ref{ddif})
but with identified fiber and base indices). We are using the 'dot' symbol
in order to distinguish the horizontal and vertical operators because on
t--bundles the indices could take the same values both for the base and
fiber objects. This allow us to define an almost complex structure $\mathbf{J%
}=\{J_{\alpha }^{\ \beta }\}$ on $TM$ as follows
\begin{equation}
\mathbf{J}(\delta _{i})=-\dot{\partial}_{i},\ \mathbf{J}(\dot{\partial}%
_{i})=\delta _{i}.  \label{alcomp}
\end{equation}
It is obvious that $\mathbf{J}$ is well--defined and $\mathbf{J}^{2}=-I.$

For d--metrics of type (\ref{dmetric}), on $TM,$ we can consider the case
when\newline
$g_{ij}(x,y)=h_{ab}(x,y),$ i. e.
\begin{equation}  \label{dmetrict}
\mathbf{G}_{(t)}= g_{ij}(x,y)dx^i\otimes dx^j + g_{ij}(x,y)\delta y^i\otimes
\delta y^j,
\end{equation}
where the index $(t)$ denotes that we have geometrical object defined on
tangent space.

An almost complex
\index{almost complex} structure $J_{\alpha}^{\ \beta}$ is compatible with a
d--metric of type (\ref{dmetrict}) and a d--connection $D$ on tangent bundle
$TM$ if the conditions
\begin{equation*}
J_{\alpha}^{\ \beta} J_{\gamma}^{\ \delta} g_{\beta \delta} = g_{\alpha
\gamma} \ \mbox{ and } \ D_{\alpha} J^{\gamma}_{\ \beta}=0
\end{equation*}
are satisfied.

The pair $(\mathbf{G}_{(t)},\mathbf{J})$ is an almost Hermitian structure on
$TM.$

One can introduce an almost sympletic 2--form associated to the almost
Hermitian structure $(\mathbf{G}_{(t)},\mathbf{J}),$
\begin{equation}  \label{hermit}
\theta = g_{ij}(x,y)\delta y^i\wedge dx^j.
\end{equation}

If the 2--form (\ref{hermit}), defined by the coefficients $g_{ij},$ is
closed, we obtain an almost K\"{a}hlerian structure in $TM.$

\begin{definition}
\label{kahlcon} An almost K\"{a}hler metric connection is a linear
connection $D^{(H)}$ on $T{%
\tilde{M}}=TM\setminus \{0\}$ with the properties:

\begin{enumerate}
\item $D^{(H)}$ preserve by parallelism the vertical distribution defined by
the N--connection structure;

\item $D^{(H)}$ is compatible with the almost K\"{a}hler structure $(\mathbf{%
G}_{(t)},\mathbf{J})$, i. e.
\begin{equation*}
D_{X}^{(H)}g=0,\ D_{X}^{(H)}J=0,\ \forall X\in \mathcal{X}\left( T{\tilde{M}}%
\right) .
\end{equation*}
\end{enumerate}
\end{definition}

By straightforward calculation we can prove that a d--connection
\\ $D\Gamma =\left(L_{jk}^i,L_{jk}^i,C_{jc}^i,C_{jc}^i\right) $
with the coefficients defined by
\begin{eqnarray}  \label{kahlerconm}
D^{(H)}_{\delta _i}\delta _j & = & L_{jk}^i \delta _i,\ D^{(H)}_{\delta _i}%
\dot{\partial}_j = L_{jk}^i \dot{\partial}_i, \\
D^{(H)}_{\delta _i}\delta _j & = & C_{jk}^i \delta _i,\ D^{(H)}_{\delta _i}%
\dot{\partial}_j = C_{jk}^i \dot{\partial}_i,  \notag
\end{eqnarray}
where $L_{jk}^i$ and $C_{ab}^e \to C_{jk}^i,$\ on $TM$ are defined by the
formulas (\ref{lccoef}), define a torsionless (see the next section on
torsion structures) metric d--connection which satisfy the compatibility
conditions (\ref{comatib}).

Almost complex structures and almost K\"{a}hler models of Finsler, Lagrange,
Hamilton and Cartan geometries (of first an higher orders) are investigated
in details in Refs. \cite{m1,m2,mhss,vbook}.

\section{ Torsions and Curvatures}

\index{torsion}
\index{curvature}

In this section we outline the basic definitions and formulas for the
torsion and curvature structures in v--bundles and cv--bundles provided with
N--connection structure.

\subsection{N--connection curvature}

\begin{enumerate}
\item The curvature $\mathbf{\Omega }$$\,$ of a nonlinear connection $%
\mathbf{N}$ in a v--bundle $\mathcal{E}$ can be defined in local form as $%
\mathbf{\ }$ \cite{ma87,ma94}:
\begin{equation*}
\mathbf{\Omega }=\frac 12\Omega _{ij}^ad^i\bigwedge d^j\otimes \partial _a,
\end{equation*}
where
\begin{eqnarray}
\Omega _{ij}^a &=&\delta _jN_i^a-\delta _iN_j^a  \label{ncurv} \\
&=&\partial _jN_i^a-\partial _iN_j^a+N_i^bN_{bj}^a-N_j^bN_{bi}^a,  \notag
\end{eqnarray}
$N_{bi}^a$ being that from (\ref{nlinearized}).

\item For the curvature $\mathbf{%
\breve{\Omega}},$ of a nonlinear connection $\mathbf{\breve{N}}$ in a
cv--bundle $\mathcal{\breve{E}}$ we introduce
\begin{equation*}
\mathbf{\breve{\Omega}}\,=\frac 12\breve{\Omega}_{ija}d^i\bigwedge
d^j\otimes \breve{\partial}^a,
\end{equation*}
where
\begin{eqnarray}
\breve{\Omega}_{ija} &=&-\breve{\delta}_j\breve{N}_{ia}+\breve{\delta}_i%
\breve{N}_{ja}  \label{ncurvcv} \\
&=&-\partial _j\breve{N}_{ia}+\partial _i\breve{N}_{ja}+\breve{N}_{ib}\breve{%
N}_{ja}^{\quad b}-\breve{N}_{jb}\breve{N}_{ja}^{\quad b},  \notag \\
\breve{N}_{ja}^{\quad b} &=&\breve{\partial}^b\breve{N}_{ja}=\partial \breve{%
N}_{ja}/\partial p_b.  \notag
\end{eqnarray}

\item Curvatures $\mathbf{\tilde{\Omega}}$$\,$ of different type of
nonlinear connections $\mathbf{\tilde{N}}$ in higher order
an\-isot\-ro\-pic bundles were analyzed for different type of
higher order tangent/dual tangent bundles and higher order
prolongations of generalized Finsler, Lagrange and Hamiloton
spaces in Refs. \cite{m1,m2,mhss} and for higher
order anisotropic superspaces and spinor bundles in Refs. \cite%
{vbook,vsp1,vhsp,vstr2}: For every higher order anisotropy shell we shall
define the coefficients (\ref{ncurv}) or (\ref{ncurvcv}) in dependence of
the fact with type of subfiber we are considering (a vector or covector
fiber).
\end{enumerate}

\subsection{d--Torsions in v- and cv--bundles}

The torsion $\mathbf{T}$ of a d--connection $\mathbf{D\ }$ in v--bundle $%
\mathcal{E}$ (cv--bundle $\mathcal{\breve{E})}$ is defined by the equation
\begin{equation}
\mathbf{T\left( X,Y\right) =XY_{\circ }^{\circ }T\doteq }D_X\mathbf{Y-}D_Y%
\mathbf{X\ -\left[ X,Y\right] .}  \label{torsion}
\end{equation}
One holds the following h- and v--decompositions
\begin{equation*}
\mathbf{T\left( X,Y\right) =T\left( hX,hY\right) +T\left( hX,vY\right)
+T\left( vX,hY\right) +T\left( vX,vY\right) .}
\end{equation*}
We consider the projections:
\begin{equation*}
\mathbf{hT\left( X,Y\right) ,vT\left( hX,hY\right) ,hT\left( hX,hY\right)
,...}
\end{equation*}
and say that, for instance, $\mathbf{hT\left( hX,hY\right) }$ is the
h(hh)--torsion of $\mathbf{D}$ ,\newline
$\mathbf{vT\left( hX,hY\right) }$ is the v(hh)--torsion of $\mathbf{D}$ and
so on.

The torsion (\ref{torsion}) in v-bundle is locally determined by five
d--tensor fields, torsions, defined as
\begin{eqnarray}
T_{jk}^{i} &=&\mathbf{hT}\left( \delta _{k},\delta _{j}\right) \cdot
d^{i},\quad T_{jk}^{a}=\mathbf{vT}\left( \delta _{k},\delta _{j}\right)
\cdot \delta ^{a},  \label{dtorsions} \\
P_{jb}^{i} &=&\mathbf{hT}\left( \partial _{b},\delta _{j}\right) \cdot
d^{i},\quad P_{jb}^{a}=\mathbf{vT}\left( \partial _{b},\delta _{j}\right)
\cdot \delta ^{a},S_{bc}^{a}=\mathbf{vT}\left( \partial _{c},\partial
_{b}\right) \cdot \delta ^{a}.  \notag
\end{eqnarray}%
Using formulas (\ref{dder}), (\ref{ddif}), (\ref{ncurv})
and (\ref{torsion}) 
we can computer \cite{ma87,ma94} in explicit form the components of torsions
(\ref{dtorsions}) 
for a d--connection of type (\ref{hgamma}) and (\ref{vgamma}):
\begin{eqnarray}
T_{.jk}^{i} &=&T_{jk}^{i}=L_{jk}^{i}-L_{kj}^{i},\quad
T_{ja}^{i}=C_{.ja}^{i},T_{aj}^{i}=-C_{ja}^{i},T_{.ja}^{i}=0,\quad
T_{.ib}^{a}=-P_{.bi}^{a}  \label{dtorsc} \\
\qquad T_{.bc}^{a} &=&S_{.bc}^{a}=C_{bc}^{a}-C_{cb}^{a},T_{.ij}^{a}=\delta
_{j}N_{i}^{a}-\delta _{j}N_{j}^{a},\quad T_{.bi}^{a}=P_{.bi}^{a}=\partial
_{b}N_{i}^{a}-L_{.bj}^{a}.  \notag
\end{eqnarray}

Formulas similar to (\ref{dtorsions}) and (\ref{dtorsc}) hold for
cv--bundles:
\begin{eqnarray}
\check{T}_{jk}^{i} &=&\mathbf{hT}\left( \delta _{k},\delta _{j}\right) \cdot
d^{i},\quad \check{T}_{jka}=\mathbf{vT}\left( \delta _{k},\delta _{j}\right)
\cdot \check{\delta}_{a},  \label{torsionsa} \\
\check{P}_{j}^{i\quad b} &=&\mathbf{hT}\left( \check{\partial}^{b},\delta
_{j}\right) \cdot d^{i},\quad \check{P}_{aj}^{\quad b}=\mathbf{vT}\left(
\check{\partial}^{b},\delta _{j}\right) \cdot \check{\delta}_{a},\check{S}%
_{a}^{\quad bc}=\mathbf{vT}\left( \check{\partial}^{c},\check{\partial}%
^{b}\right) \cdot \check{\delta}_{a}.  \notag
\end{eqnarray}%
and
\begin{eqnarray}
\check{T}_{.jk}^{i} &=&\check{T}_{jk}^{i}=L_{jk}^{i}-L_{kj}^{i},\quad
\check{T}_{j}^{ia}=\check{C}_{.j}^{i\ a},\check{T}_{~\ j}^{ia}=-\check{C}%
_{j}^{i~a},~\check{T}_{.j}^{i~a}=0,\check{T}_{a~b}^{~j}=-\check{P}_{a\
b}^{~j}  \label{dtorsca} \\
\qquad \check{T}_{a}^{~bc} &=&\check{S}_{a}^{~bc}=\check{C}_{a}^{~bc}-%
\check{C}_{a}^{~cb},\check{T}_{.ija}=-\delta _{j}\check{N}_{ia}+\delta _{j}%
\check{N}_{ja},\quad \check{T}_{a}^{~bi}=\check{P}_{a}^{~bi}=-\check{\partial%
}^{b}\check{N}_{ia}-\check{L}_{a}^{~bi}.  \notag
\end{eqnarray}

The formulas for torsion can be generalized for hvc--bundles (on every shell
we must write (\ref{dtorsc}) or (\ref{dtorsca}) in dependence of the type of
shell, vector or co-vector one, we are dealing).

\subsection{d--Curvatures in v- and cv--bundles}

\index{d--curvatures}

The curvature $\mathbf{R}$ of a d--connection in v--bundle $\mathcal{E}$ is
defined by the equation
\begin{equation*}
\mathbf{R}\left( X,Y\right) Z = XY_{\bullet }^{\bullet }R\bullet Z =D_XD_Y{\
Z} - D_Y D_X Z - D_{[X,Y]}{Z.}
\end{equation*}
One holds the next properties for the h- and v--decompositions of curvature:
\begin{eqnarray}  \label{curvaturehv}
\mathbf{vR}\left( X,Y\right) hZ &=& 0,\ \mathbf{hR}\left( X,Y\right) vZ=0, \\
\mathbf{R}\left( X,Y\right) Z & = & \mathbf{hR} \left( X,Y\right) hZ +
\mathbf{vR} \left( X,Y\right) vZ.  \notag
\end{eqnarray}

From (\ref{curvaturehv}) 
and the equation $\mathbf{R\left( X,Y\right) =-R\left( Y,X\right) }$ we get
that the curvature of a d--con\-nec\-ti\-on $\mathbf{D}$ in $\mathcal{E}$ is
completely determined by the following six d--tensor fields:
\begin{eqnarray}
R_{h.jk}^{.i} &=& d^i\cdot \mathbf{R}\left( \delta _k,\delta _j\right)
\delta _h,~R_{b.jk}^{.a}=\delta ^a\cdot \mathbf{R}\left( \delta _k,\delta
_j\right) \partial _b,  \label{rps} \\
P_{j.kc}^{.i} &=& d^i\cdot \mathbf{R}\left( \partial _c,\partial _k\right)
\delta _j,~P_{b.kc}^{.a}=\delta ^a\cdot \mathbf{R}\left( \partial
_c,\partial _k\right) \partial _b,  \notag \\
S_{j.bc}^{.i} &=& d^i\cdot \mathbf{R}\left( \partial _c,\partial _b\right)
\delta _j,~S_{b.cd}^{.a}=\delta ^a\cdot \mathbf{R}\left( \partial
_d,\partial _c\right) \partial _b.  \notag
\end{eqnarray}
By a direct computation, using 
(\ref{dder}),(\ref{ddif}),(\ref{hgamma}),(\ref{vgamma}) and (\ref{rps}) we
get:
\begin{eqnarray}
R_{h.jk}^{.i} &=& {\delta}_h L_{.hj}^i - {\delta}_j L_{.hk}^i +L_{.hj}^m
L_{mk}^i-L_{.hk}^m L_{mj}^i+C_{.ha}^i R_{.jk}^a,  \label{dcurvatures} \\
R_{b.jk}^{.a} &=& {\delta}_k L_{.bj}^a - {\delta}_j L_{.bk}^a +L_{.bj}^c
L_{.ck}^a - L_{.bk}^c L_{.cj}^a + C_{.bc}^a R_{.jk}^c,  \notag \\
P_{j.ka}^{.i} &=& {\partial}_a L_{.jk}^i - \left({\delta}_k C_{.ja}^i +
L_{.lk}^i C_{.ja}^l - L_{.jk}^l C_{.la}^i - L_{.ak}^c C_{.jc}^i\right)
+C_{.jb}^i P_{.ka}^b,  \notag \\
P_{b.ka}^{.c} &=& {\partial}_a L_{.bk}^c -\left( {\delta}_k C_{.ba}^c +
L_{.dk}^{c} C_{.ba}^d - L_{.bk}^d C_{.da}^c - L_{.ak}^d C_{.bd}^c \right)
+C_{.bd}^c P_{.ka}^d,  \notag \\
S_{j.bc}^{.i} &=& {\partial}_c C_{.jb}^i - {\partial}_b C_{.jc}^i +
C_{.jb}^h C_{.hc}^i - C_{.jc}^h C_{hb}^i,  \notag \\
S_{b.cd}^{.a} & = & {\partial}_d C_{.bc}^a - {\partial}_c C_{.bd}^a +
C_{.bc}^e C_{.ed}^a - C_{.bd}^eC_{.ec}^a.  \notag
\end{eqnarray}

We note that d--torsions (\ref{dtorsc}) and d--curvatures (\ref{dcurvatures}%
) 
are computed in explicit form by particular cases of d--connections (\ref%
{berwald}), (\ref{inters}) and (\ref{dchrist}).

For cv--bundles we have
\begin{eqnarray}
\check{R}_{h.jk}^{.i} &=&d^{i}\cdot \mathbf{R}\left( \delta _{k},\delta
_{j}\right) \delta _{h},~\check{R}_{\ a.jk}^{b}=%
\check{\delta}_{a}\cdot \mathbf{R}\left( \delta _{k},\delta _{j}\right)
\check{\partial}^{b},  \label{rpsa} \\
\check{P}_{j.k}^{.i\quad c} &=&d^{i}\cdot \mathbf{R}\left( \check{\partial}%
^{c},\partial _{k}\right) \delta _{j},~\check{P}_{\quad a.k}^{b\quad c}=%
\check{\delta}_{a}\cdot \mathbf{R}\left( \check{\partial}^{c},\partial
_{k}\right) \check{\partial}^{b},  \notag \\
\check{S}_{j.}^{.ibc} &=&d^{i}\cdot \mathbf{R}\left( \check{\partial}^{c},%
\check{\partial}^{b}\right) \delta _{j},~\check{S}_{.a}^{b.cd}=\check{\delta}%
_{a}\cdot \mathbf{R}\left( \check{\partial}^{d},\check{\partial}^{c}\right)
\check{\partial}^{b}.  \notag
\end{eqnarray}%
and
\begin{eqnarray}
\check{R}_{h.jk}^{.i} &=&{\check{\delta}}_{h}L_{.hj}^{i}-{\check{\delta}}%
_{j}L_{.hk}^{i}+L_{.hj}^{m}L_{mk}^{i}-L_{.hk}^{m}L_{mj}^{i}+C_{.h}^{i~a}%
\check{R}_{.ajk},  \label{dcurvaturesa} \\
\check{R}_{.ajk}^{b.} &=&{\check{\delta}}_{k}\check{L}_{a.j}^{~b}-{\check{%
\delta}}_{j}\check{L}_{~b~k}^{a}+\check{L}_{cj}^{~b}\check{L}_{.ak}^{~c}-%
\check{L}_{ck}^{b}\check{L}_{a.j}^{~c}+\check{C}_{a}^{~bc}\check{R}_{c.jk},
\notag \\
\check{P}_{j.k}^{.i~a} &=&{\check{\partial}}^{a}L_{.jk}^{i}-\left( {\check{%
\delta}}_{k}\check{C}_{.j}^{i~a}+L_{.lk}^{i}\check{C}_{.j}^{l~a}-L_{.jk}^{l}%
\check{C}_{.l}^{i~a}-\check{L}_{ck}^{~a}\check{C}_{.j}^{i~c}\right) +%
\check{C}_{.j}^{i~b}\check{P}_{bk}^{~\quad a},  \notag \\
\check{P}_{ck}^{b~a} &=&{\check{\partial}}^{a}\check{L}_{c.k}^{~b}-({\check{%
\delta}}_{k}\check{C}_{c.}^{~ba}+\check{L}_{c.k}^{bd}\check{C}_{d}^{\ ba}-%
\check{L}_{d.k}^{\quad b}\check{C}_{c.}^{\ ad})-\check{L}_{dk}^{\quad a}%
\check{C}_{c.}^{\ bd})+\check{C}_{c.}^{\ bd}\check{P}_{d.k}^{\quad a},
\notag \\
\check{S}_{j.}^{.ibc} &=&{\check{\partial}}^{c}\check{C}_{.j}^{i\ b}-{\check{%
\partial}}^{b}\check{C}_{.j}^{i\ c}+\check{C}_{.j}^{h\ b}\check{C}_{.h}^{i\
c}-\check{C}_{.j}^{h\ c}\check{C}_{h}^{i\ b},  \notag \\
\check{S}_{\ a.}^{b\ cd} &=&{\check{\partial}}^{d}\check{C}_{a.}^{\ bc}-{%
\check{\partial}}^{c}\check{C}_{a.}^{\ bd}+\check{C}_{e.}^{\ bc}\check{C}%
_{a.}^{\ ed}-\check{C}_{e.}^{\ bd}\check{C}_{.a}^{\ ec}.  \notag
\end{eqnarray}

The formulas for curvature can be also generalized for hvc--bundles (on
every shell we must write (\ref{dtorsc}) or (\ref{torsionsa}) in dependence
of the type of shell, vector or co-vector one, we are dealing).

\section{Generalizations of Finsler Spaces}

We outline the basic definitions and formulas for Finsler, Lagrange and
generalized Lagrange spaces (constructed on tangent bundle) and for Cartan,
Hamilton and generalized Hamilton spaces (constructed on cotangent bundle).
The original results are given in details in monographs \cite{ma87,ma94,mhss}

\subsection{Finsler Spaces}

\index{Finsler space}

The Finsler geometry is modeled on tangent bundle $TM.$

\begin{definition}
A Finsler space (manifold) is a pair $F^{n}=\left( M,F(x,y)\right) $ \ where
$M$ is a real $n$--dimensional differentiable manifold and $F:TM\rightarrow
\mathcal{R}$ \ is a scalar function which satisfy the following conditions:

\begin{enumerate}
\item $F$ is a differentiable function on the manifold $%
\widetilde{TM}$ $=TM\backslash \{0\}$ and $F$ is continous on the null
section of the projection $\pi :TM\rightarrow M;$

\item $F$ is a positive function, homogeneous on the fibers of the $TM,$ i.
e. $F(x,\lambda y)=\lambda F(x,y),\lambda \in \mathcal{R};$

\item The Hessian of $F^{2}$ with elements
\begin{equation}
g_{ij}^{(F)}(x,y)=\frac{1}{2}\frac{\partial ^{2}F^{2}}{\partial
y^{i}\partial y^{j}}  \label{finm}
\end{equation}
is positively defined on $\widetilde{TM}.$
\end{enumerate}
\end{definition}

The function $F(x,y)$ and $g_{ij}(x,y)$ are called respectively the
fundamental function and the fundamental (or metric) tensor of the Finsler
space $F.$

One considers ''anisotropic'' (depending on directions $y^{i})$ Christoffel
symbols, for simplicity we write $g_{ij}^{(F)}=g_{ij},$
\begin{equation*}
\gamma _{~jk}^{i}(x,y)=\frac{1}{2}g^{ir}\left( \frac{\partial g_{rk}}{%
\partial x^{j}}+\frac{\partial g_{jr}}{\partial x^{k}}-\frac{\partial g_{jk}%
}{\partial x^{r}}\right) ,
\end{equation*}
which are used for definition of the Cartan N--connection,
\begin{equation}
N_{(c)~j}^{i}=\frac{1}{2}\frac{\partial }{\partial y^{j}}\left[ \gamma
_{~nk}^{i}(x,y)y^{n}y^{k}\right] .  \label{ncartan}
\end{equation}
This N--connection can be used for definition of an almost complex structure
like in (\ref{alcomp}) and to define on $TM$ a d--metric
\begin{equation}
\mathbf{G}_{(F)}=g_{ij}(x,y)dx^{i}\otimes dx^{j}+g_{ij}(x,y)\delta
y^{i}\otimes \delta y^{j},  \label{dmfin}
\end{equation}
with $g_{ij}(x,y)$ taken as (\ref{finm}).

Using the Cartan N--connection (\ref{ncartan}) and Finsler metric tensor (%
\ref{finm}) (or, equivalently, the d--metric (\ref{dmfin})) we can introduce
the canonical d--connection
\begin{equation*}
D\Gamma \left( N_{(c)}\right) =\Gamma _{(c)\beta \gamma }^{\alpha }=\left(
L_{(c)~jk}^{i},C_{(c)~jk}^{i}\right)
\end{equation*}
with the coefficients computed like in (\ref{kahlerconm}) and (\ref{lccoef})
with $h_{ab}\rightarrow g_{ij}.$ The d--connection $D\Gamma \left(
N_{(c)}\right) $ has the unique property that it is torsionless and
satisfies the metricity conditions both for the horizontal and vertical
components, i. e. $D_{\alpha }g_{\beta \gamma }=0.$

The d--curvatures
\begin{equation*}
\check{R}_{h.jk}^{.i}=\{\check{R}_{h.jk}^{.i},\check{P}_{j.k}^{.i\quad
l},S_{(c)j.kl}^{.i}\}
\end{equation*}
on a Finsler space provided with Cartan N--connection and Finsler metric
structures are computed following the formulas (\ref{dcurvatures}) when the $%
a,b,c...$ indices are identified with $i,j,k,...$ indices. It should be
emphasized that in this case all values $g_{ij,}\Gamma _{(c)\beta \gamma
}^{\alpha }$ and $R_{(c)\beta .\gamma \delta }^{.\alpha }$ are defined by a
fundamental function $F\left( x,y\right) .$

In general, we can consider that a Finsler space is provided with a metric $%
g_{ij}=\partial ^{2}F^{2}/2\partial y^{i}\partial y^{j},$ but the
N--connection and d--connection are be defined in a different manner, even
not be determined by $F.$

\subsection{Lagrange and Generalized Lagrange Spaces}

\index{Lagrange space}

The notion of Finsler spaces was extended by J. Kern \cite{ker} and R. Miron %
\cite{mironlg}. It is widely developed in monographs \cite{ma87,ma94} and
exteded to superspaces by S. Vacaru \cite{vlasg,vstr2,vbook}.

The idea of extension was to consider instead of the homogeneous fundamental
function $F(x,y)$ in a Finsler space a more general one, a Lagrangian $%
L\left( x,y\right) $, defined as a differentiable mapping $L:(x,y)\in
TM\rightarrow L(x,y)\in \mathcal{R},$ of class $C^{\infty }$ on manifold $%
\widetilde{TM}$ and continous on the null section $0:M\rightarrow TM$ of the
projection $\pi :TM\rightarrow M.$ A Lagrangian is regular if it is
differentiable and the Hessian
\begin{equation}
g_{ij}^{(L)}(x,y)=\frac{1}{2}\frac{\partial ^{2}L^{2}}{\partial
y^{i}\partial y^{j}}  \label{lagm}
\end{equation}
is of rank $n$ on $M.$

\begin{definition}
A Lagrange space is a pair $L^{n}=\left( M,L(x,y)\right) $ where $M$ is a
smooth real $n$--dimensional manifold provided with regular Lagrangian \ $%
L(x,y)$ structure $L:TM\rightarrow \mathcal{R}$ $\ $for which $g_{ij}(x,y)$
from (\ref{lagm}) has a constant signature over the manifold $\widetilde{TM}%
. $
\end{definition}

The fundamental Lagrange function $L(x,y)$ defines a canonical
N--con\-nec\-ti\-on
\begin{equation*}
N_{(cL)~j}^{i}=\frac{1}{2}\frac{\partial }{\partial y^{j}}\left[
g^{ik}\left( \frac{\partial ^{2}L^{2}}{\partial y^{k}\partial y^{h}}y^{h}-%
\frac{\partial L}{\partial x^{k}}\right) \right]
\end{equation*}
as well a d-metric
\begin{equation}
\mathbf{G}_{(L)}=g_{ij}(x,y)dx^{i}\otimes dx^{j}+g_{ij}(x,y)\delta
y^{i}\otimes \delta y^{j},  \label{dmlag}
\end{equation}
with $g_{ij}(x,y)$ taken as (\ref{lagm}). As well we can introduce an almost
K\"{a}hlerian structure and an almost Hermitian model of $L^{n},$ denoted as
$H^{2n}$ as in the case of Finsler spaces but with a proper fundamental
Lagange function and metric tensor $g_{ij}.$ The canonical metric
d--connection $D\Gamma \left( N_{(cL)}\right) =\Gamma _{(cL)\beta \gamma
}^{\alpha }=\left( L_{(cL)~jk}^{i},C_{(cL)~jk}^{i}\right) $ is to computed
by the same formulas (\ref{kahlerconm}) and (\ref{lccoef}) with $%
h_{ab}\rightarrow g_{ij}^{(L)},$ for $N_{(cL)~j}^{i}.$ The d--torsions (\ref%
{dtorsc}) and d--curvatures (\ref{dcurvatures}) are defined, in this case,
by $L_{(cL)~jk}^{i}$ and $C_{(cL)~jk}^{i}.$ We also note that instead of $%
N_{(cL)~j}^{i}$ and $\Gamma _{(cL)\beta \gamma }^{\alpha }$ one can consider
on a $L^{n}$--space arbitrary N--connections $N_{~j}^{i},$ d--connections $%
\Gamma _{\beta \gamma }^{\alpha }$ which are not defined only by $L(x,y)$
and $g_{ij}^{(L)}$ but can be metric, or non--metric with respect to the
Lagrange metric.

The next step of generalization is to consider an arbitrary metric $%
g_{ij}\left( x,y\right) $ on $TM$ instead of (\ref{lagm}) which is the
second derivative of ''anisotropic'' coordinates $y^{i}$ of a Lagrangian %
\cite{mironlg,mironlgg}.

\begin{definition}
A generalized Lagrange space is a pair $GL^{n}=\left( M,g_{ij}(x,y)\right) $
where $g_{ij}(x,y)$ is a covariant, symmetric d--tensor field, of rank $n$
and of constant signature on $\widetilde{TM}.$
\end{definition}

\bigskip One can consider different classes of N-- and d--connections on $%
TM, $ which are compatible (metric) or non compatible with
(\ref{dmlag}) for arbitrary $g_{ij}(x,y).$ We can apply all
formulas for d--connections, N--curvatures, d--torsions and
d--curvatures as in a v--bundle $\mathcal{E},$ but reconsidering
them on $TM,$ by changing \ $h_{ab}\rightarrow g_{ij}(x,y)$ and
$N_{i}^{a}\rightarrow N_{~i}^{k}.$

\subsection{Cartan Spaces}

\index{Cartan space}

The theory of Cartan spaces (see, for instance, \cite{run,kaw1}) \ was
formulated in a new fashion in R. Miron's works \cite{mironc1,mironc2} by
considering them as duals to the Finsler spaces (see details and references
in \cite{mhss}). Roughly, a Cartan space is constructed on a cotangent
bundle $T^{\ast }M$ like a Finsler space on the corresponding tangent bundle
$TM.$

Consider a real smooth manifold $M,$ the cotangent bundle $\left( T^{\ast
}M,\pi ^{\ast },M\right) $ and the manifold $%
\widetilde{T^{\ast }M}=T^{\ast }M\backslash \{0\}.$

\begin{definition}
A Cartan space is a pair $C^{n}=\left( M,K(x,p)\right) $ \ such that $%
K:T^{\ast }M\rightarrow \mathcal{R}$ is a scalar function which satisfy the
following conditions:

\begin{enumerate}
\item $K$ is a differentiable function on the manifold $\widetilde{T^{\ast }M%
}$ $=T^{\ast }M\backslash \{0\}$ and continuous on the null
section of the projection $\pi ^{\ast }:T^{\ast }M\rightarrow M;$

\item $K$ is a positive function, homogeneous on the fibers of the $T^{\ast
}M,$ i. e. $K(x,\lambda p)=\lambda F(x,p),\lambda \in \mathcal{R};$

\item The Hessian of $K^{2}$ with elements
\begin{equation}
\check{g}_{(K)}^{ij}(x,p)=\frac{1}{2}\frac{\partial ^{2}K^{2}}{\partial
p_{i}\partial p_{j}}  \label{carm}
\end{equation}
is positively defined on $\widetilde{T^{\ast }M}.$
\end{enumerate}
\end{definition}

The function $K(x,y)$ and $\check{g}^{ij}(x,p)$ are called \ respectively
the fundamental function and the fundamental (or metric) tensor of the
Cartan space $C^{n}.$ We use symbols like $"\check{g}"$ as to emphasize that
the geometrical objects are defined on a dual space.

One considers ''anisotropic'' (depending on directions, momenta, $p_{i})$
\newline
Christoffel symbols, for simplicity, we write the inverse to (\ref{carm}) as $%
g_{ij}^{(K)}=\check{g}_{ij},$
\begin{equation*}
\check{\gamma}_{~jk}^{i}(x,p)=\frac{1}{2}\check{g}^{ir}\left( \frac{\partial
\check{g}_{rk}}{\partial x^{j}}+\frac{\partial \check{g}_{jr}}{\partial x^{k}%
}-\frac{\partial \check{g}_{jk}}{\partial x^{r}}\right) ,
\end{equation*}
which are used for definition of the canonical N--connection,
\begin{equation}
\check{N}_{ij}=\check{\gamma}_{~ij}^{k}p_{k}-\frac{1}{2}\gamma
_{~nl}^{k}p_{k}p^{l}{\breve{\partial}}^{n}\check{g}_{ij},~{\breve{\partial}}%
^{n}=\frac{\partial }{\partial p_{n}}.  \label{nccartan}
\end{equation}
This N--connection can be used for definition of an almost complex structure
like in (\ref{alcomp}) and to define on $T^{\ast }M$ a d--metric
\begin{equation}
\mathbf{\check{G}}_{(k)}=\check{g}_{ij}(x,p)dx^{i}\otimes dx^{j}+\check{g}%
^{ij}(x,p)\delta p_{i}\otimes \delta p_{j},  \label{dmcar}
\end{equation}
with $\check{g}^{ij}(x,p)$ taken as (\ref{carm}).

Using the canonical N--connection (\ref{nccartan}) and Finsler metric tensor
(\ref{carm}) (or, equivalently, the d--metric (\ref{dmcar}) we can introduce
the canonical d--connection
\begin{equation*}
D\check{\Gamma}\left( \check{N}_{(k)}\right) =\check{\Gamma}_{(k)\beta
\gamma }^{\alpha }=\left( \check{H}_{(k)~jk}^{i},\check{C}_{(k)~i}^{\quad
jk}\right)
\end{equation*}%
with the coefficients \ are computed
\begin{equation*}
\check{H}_{(k)~jk}^{i}=\frac{1}{2}\check{g}^{ir}\left( \check{\delta}_{j}%
\check{g}_{rk}+\check{\delta}_{k}\check{g}_{jr}-\check{\delta}_{r}\check{g}%
_{jk}\right) ,\check{C}_{(k)~i}^{\quad jk}=\check{g}_{is}{\breve{\partial}}%
^{s}\check{g}^{jk},
\end{equation*}%
The d--connection $D\check{\Gamma}\left( \check{N}_{(k)}\right) $ has the
unique property that it is torsionless and satisfies the metricity
conditions both for the horizontal and vertical components, i. e. $\check{D}%
_{\alpha }\check{g}_{\beta \gamma }=0.$

The d--curvatures
\begin{equation*}
\check{R}_{(k)\beta .\gamma \delta }^{.\alpha
}=\{R_{(k)h.jk}^{.i},P_{(k)j.km}^{.i},\check{S}_{j.}^{.ikl}\}
\end{equation*}
on a Finsler space provided with Cartan N--connection and Finsler metric
structures are computed following the formulas (\ref{dcurvaturesa}) when the
$a,b,c...$ indices are identified with $i,j,k,...$ indices. It should be
emphasized that in this case all values $\check{g}_{ij,}\check{\Gamma}%
_{(k)\beta \gamma }^{\alpha }$ and $\check{R}_{(k)\beta .\gamma \delta
}^{.\alpha }$ are defined by a fundamental function $K\left( x,p\right) .$

In general, we can consider that a Cartan space is provided with a metric $%
\check{g}^{ij}=\partial ^{2}K^{2}/2\partial p_{i}\partial p_{j},$ but the
N--connection and d--connection could be defined in a different manner, even
not be determined by $K.$

\subsection{ Generalized Hamilton and Hamilton Spaces}

\index{Hamilton space}

The geometry of Hamilton spaces was defined and investigated by R. Miron in
the papers \cite{mironh1,mironh2,mironh3} (see details and references in %
\cite{mhss}). It was developed on the cotangent bundel as a dual geometry to
the geometry of Lagrange spaces. \ Here we start with the definition of
generalized Hamilton spaces and then consider the particular case.

\begin{definition}
A generalized Hamilton space is a pair\newline
$GH^{n}=\left( M,%
\check{g}^{ij}(x,p)\right) $ where $M$ is a real $n$--dimensional manifold
and $\check{g}^{ij}(x,p)$ is a contravariant, symmetric, nondegenerate of
rank $n$ and of constant signature on $\widetilde{T^{\ast }M}.$
\end{definition}

\bigskip The value $\check{g}^{ij}(x,p)$ is called the fundamental (or
metric) tensor of the space $GH^{n}.$ One can define such values for every
paracompact manifold $M.$ In general, a N--connection on $GH^{n}$ is not
determined by $\check{g}^{ij}.$ Therefore we can consider arbitrary
coefficients $\check{N}_{ij}\left( x,p\right) $ and define on $T^{\ast }M$ a
d--metric like (\ref{dmetricvc})
\begin{equation}
{\breve{\mathbf{G}}}={\breve{g}}_{\alpha \beta }\left( {\breve{u}}\right) {%
\breve{\delta}}^{\alpha }\otimes {\breve{\delta}}^{\beta }={\breve{g}}%
_{ij}\left( {\breve{u}}\right) d^{i}\otimes d^{j}+{\check{g}}^{ij}\left( {%
\breve{u}}\right) {\breve{\delta}}_{i}\otimes {\breve{\delta}}_{j},
\label{dmghs}
\end{equation}
This N--coefficients $\check{N}_{ij}\left( x,p\right) $ and d--metric
structure (\ref{dmghs}) allow to define an almost K\"{a}hler model of
generalized Hamilton spaces and to define canonical d--connections,
d--torsions and d-curvatures (see respectively the formulas (\ref{lccoef}), (%
\ref{inters}), (\ref{dtorsca}) and (\ref{dcurvatures}) with the fiber
coefficients redefined for the cotangent bundle $T^{\ast }M$ ).

A generalized Hamilton space $GH^{n}=\left( M,\check{g}^{ij}(x,p)\right) $
is called reducible to a Hamilton one if there exists a Hamilton function $%
H\left( x,p\right) $ on $T^{\ast }M$ such that
\begin{equation}
\check{g}^{ij}(x,p)=\frac{1}{2}\frac{\partial ^{2}H}{\partial p_{i}\partial
p_{j}}.  \label{hsm}
\end{equation}

\begin{definition}
A Hamilton space is a pair $H^{n}=\left( M,H(x,p)\right) $ \ such that $%
H:T^{\ast }M\rightarrow \mathcal{R}$ is a scalar function which satisfy the
following conditions:

\begin{enumerate}
\item $H$ is a differentiable function on the manifold $\widetilde{T^{\ast }M%
}$ $=T^{\ast }M\backslash \{0\}$ and continuous on the null
section of the projection $\pi ^{\ast }:T^{\ast }M\rightarrow M;$

\item The Hessian of $H$ with elements (\ref{hsm}) is positively defined on $%
\widetilde{T^{\ast }M}$ and $\check{g}^{ij}(x,p)$ is nondegenerate matrix of
rank $n$ and of constant signature.
\end{enumerate}
\end{definition}

For Hamilton spaces the canonical N--connection (defined by $H$ and its
Hessian) exists,
\begin{equation*}
\check{N}_{ij}=\frac{1}{4}\{\check{g}_{ij},H\}-\frac{1}{2}\left( \check{g}%
_{ik}\frac{\partial ^{2}H}{\partial p_{k}\partial x^{j}}+\check{g}_{jk}\frac{%
\partial ^{2}H}{\partial p_{k}\partial x^{i}}\right) ,
\end{equation*}
where the Poisson brackets, for arbitrary functions $f$ and $g$ on $T^{\ast
}M,$ act as
\begin{equation*}
\{f,g\}=\frac{\partial f}{\partial p_{i}}\frac{\partial g}{\partial x^{i}}-%
\frac{\partial g}{\partial p_{i}}\frac{\partial p}{\partial x^{i}}.
\end{equation*}

The canonical d--connection $D\check{\Gamma}\left( \check{N}_{(c)}\right) =%
\check{\Gamma}_{(c)\beta \gamma }^{\alpha }=\left( \check{H}_{(c)~jk}^{i},%
\check{C}_{(c)~i}^{\quad jk}\right) $is defined by the coefficients
\begin{equation*}
\check{H}_{(c)~jk}^{i}=\frac{1}{2}\check{g}^{is}\left( \check{\delta}_{j}%
\check{g}_{sk}+\check{\delta}_{k}\check{g}_{js}-\check{\delta}_{s}\check{g}%
_{jk}\right) ,\check{C}_{(c)~i}^{\quad jk}=-\frac{1}{2}\check{g}_{is}\check{%
\partial}^{j}\check{g}^{sk}.
\end{equation*}%
In result we can compute the d--torsions and d--curvatures like on
cv--bundle \ or on Cartan spaces. On Hamilton spaces all such
objects are defined by the Hamilton function $H(x,p)$ and indices
have to be reconsidered for co--fibers of the co-tangent bundle.

\section{Gravity on Vector Bundles}

The components of the Ricci d--tensor
\index{Ricci d--tensor}
\begin{equation*}
R_{\alpha \beta }=R_{\alpha .\beta \tau }^{.\tau }
\end{equation*}%
with respect to a locally adapted frame (\ref{ddif}) are as follows:
\begin{eqnarray}
R_{ij} &=&R_{i.jk}^{.k},\quad
R_{ia}=-^{2}P_{ia}=-P_{i.ka}^{.k},R_{ai}=^{1}P_{ai}=P_{a.ib}^{.b},\quad
R_{ab}=S_{a.bc}^{.c}.  \label{2.33} \\
&&  \notag
\end{eqnarray}%
We point out that because, in general, $^{1}P_{ai}\neq ~^{2}P_{ia}$ the
Ricci d-tensor is non symmetric.

Having defined a d-metric of type in $\mathcal{E}$ we can introduce the
scalar curvature of d--connection $\mathbf{D}$:
\begin{equation}
{%
\overleftarrow{R}}=G^{\alpha \beta }R_{\alpha \beta }=R+S,  \label{2.34}
\end{equation}
where $R=g^{ij}R_{ij}$ and $S=h^{ab}S_{ab}.$

For our further considerations it will be also useful to use an alternative
way of definition torsion (\ref{torsion}) and curvature (\ref{curvaturehv})
by using the commutator
\begin{equation*}
\Delta _{\alpha \beta }\doteq \bigtriangledown _\alpha \bigtriangledown
_\beta -\bigtriangledown _\beta \bigtriangledown _\alpha =2\bigtriangledown
_{[\alpha }\bigtriangledown _{\beta ]}.
\end{equation*}
For components of d--torsion we have
\begin{equation*}
\Delta _{\alpha \beta }f=T_{.\alpha \beta }^\gamma \bigtriangledown _\gamma f
\end{equation*}
for every scalar function $f\,\,$ on $\mathcal{E}$. Curvature can be
introduced as an operator acting on arbitrary d-vector $V^\delta :$
\begin{equation}
(\Delta _{\alpha \beta }-T_{.\alpha \beta }^\gamma \bigtriangledown _\gamma
)V^\delta =R_{~\gamma .\alpha \beta }^{.\delta }V^\gamma  \label{2.37}
\end{equation}
(we note that in this section we shall follow conventions of Miron and
Anastasiei \cite{ma87,ma94} on d-tensors; we can obtain corresponding
Penrose and Rindler abstract index formulas \cite{penr1,penr2} just for a
trivial N-connection structure and by changing denotations for components of
torsion and curvature in this manner:\ $T_{.\alpha \beta }^\gamma
\rightarrow T_{\alpha \beta }^{\quad \gamma }$ and $R_{~\gamma .\alpha \beta
}^{.\delta }\rightarrow R_{\alpha \beta \gamma }^{\qquad \delta }).$

Here we also note that torsion and curvature of a d-connection on $\mathcal{E%
}$ satisfy generalized for locally anisotropic spaces Ricci and Bianchi
identities \cite{ma87,ma94} which in terms of components (\ref{2.37}) are
written respectively as
\begin{equation}
R_{~[\gamma .\alpha \beta ]}^{.\delta }+\bigtriangledown _{[\alpha
}T_{.\beta \gamma ]}^\delta +T_{.[\alpha \beta }^\nu T_{.\gamma ]\nu
}^\delta =0  \label{2.38}
\end{equation}
and
\begin{equation}
\bigtriangledown _{[\alpha }R_{|\nu |\beta \gamma ]}^{\cdot \sigma
}+T_{\cdot [\alpha \beta }^\delta R_{|\nu |.\gamma ]\delta }^{\cdot \sigma
}=0.  \label{2.39}
\end{equation}
Identities (\ref{2.38}) and (\ref{2.39}) can be proved similarly as in \cite%
{penr1} by taking into account that indices play a distinguished character.

We can also consider a la-generalization of the so-called conformal Weyl
tensor (see, for instance, \cite{penr1}) which can be written as a d-tensor
in this form:
\begin{eqnarray}
C_{\quad \alpha \beta }^{\gamma \delta } &=&R_{\quad \alpha \beta }^{\gamma
\delta }-\frac{4}{n+m-2}R_{\quad \lbrack \alpha }^{[\gamma }~\delta _{\quad
\beta ]}^{\delta ]}+\frac{2}{(n+m-1)(n+m-2)}{\overleftarrow{R}~\delta
_{\quad \lbrack \alpha }^{[\gamma }~\delta _{\quad \beta ]}^{\delta ]}.}
  \notag \\
&& \label{2.40}
\end{eqnarray}%
This object is conformally invariant on locally anisotropic spaces provided
with d-connection generated by d-metric structures.

The Einstein equations and conservation laws on v-bundles
provided with N--con\-nec\-ti\-on structures are studied in
detail in \cite{ma87,ma94,ana86,ana87}. In Ref. \cite{vg} we
proved that the locally anisotropic gravity can be formulated in
a gauge like manner and analyzed the conditions when the Einstein
locally anisotropic gravitational field equations are equivalent
to a corresponding form of Yang-Mills equations. In this
subsection we write the locally anisotropic gravitational field
equations in a form more convenient for theirs equivalent
reformulation in locally anisotropic spinor variables.

We define d-tensor $\Phi _{\alpha \beta }$ as to satisfy conditions
\begin{equation}
-2\Phi _{\alpha \beta }\doteq R_{\alpha \beta }-\frac 1{n+m}\overleftarrow{R}%
g_{\alpha \beta }  \label{2.41}
\end{equation}
which is the torsionless part of the Ricci tensor for locally isotropic
spaces \cite{penr1,penr2}, i.e. $\Phi _\alpha ^{~~\alpha }\doteq 0$.\ The
Einstein equations on locally anisotropic spaces
\begin{equation}
\overleftarrow{G}_{\alpha \beta }+\lambda g_{\alpha \beta }= \kappa
E_{\alpha \beta },  \label{2.42}
\end{equation}
where
\begin{equation}
\overleftarrow{G}_{\alpha \beta }=R_{\alpha \beta }-\frac 12\overleftarrow{R}%
g_{\alpha \beta }  \label{2.43}
\end{equation}
is the Einstein d-tensor, $\lambda $ and $\kappa $ are correspondingly the
cosmological and gravitational constants and by $E_{\alpha \beta }$ is
denoted the locally anisotropic energy-momentum d-tensor \cite{ma87,ma94},
can be rewritten in equivalent form:
\begin{equation}
\Phi _{\alpha \beta }=-\frac \kappa 2(E_{\alpha \beta }-\frac 1{n+m}E_\tau
^{~\tau }~g_{\alpha \beta }).  \label{2.44}
\end{equation}
Because the locally anisotropic spaces generally have nonzero
torsions we shall add to (\ref{2.44}) (equivalently to
(\ref{2.42})) a system of algebraic d-field equations with the
source $S_{~\beta \gamma }^\alpha $ being the locally anisotropic
spin density of matter (if we consider a variant of locally
anisotropic Einstein--Cartan theory):
\begin{equation}
T_{~\alpha \beta }^\gamma +2\delta _{~[\alpha }^\gamma T_{~\beta ]\delta
}^\delta =\kappa S_{~\alpha \beta .}^\gamma  \label{2.45}
\end{equation}
From (\ref{2.38}) and (\ref{2.45}) one follows the conservation law of
locally anisotropic spin matter:
\begin{equation*}
\bigtriangledown _\gamma S_{~\alpha \beta }^\gamma -T_{~\delta \gamma
}^\delta S_{~\alpha \beta }^\gamma =E_{\beta \alpha }-E_{\alpha \beta }.
\end{equation*}

Finally, in this section, we remark that all presented geometric
constructions contain those elaborated for generalized Lagrange spaces \cite%
{ma87,ma94} (for which a tangent bundle $TM$ is considered instead of a
v-bundle $\mathcal{E}$ ). We also note that the Lagrange (Finsler) geometry
is characterized by a metric with components parametized as $g_{ij}=\frac 12%
\frac{\partial ^2\mathcal{L}}{\partial y^i\partial y^j}$ $\left(
g_{ij}=\frac 12\frac{\partial ^2\Lambda ^2}{\partial y^i\partial y^j}\right)
$ and $h_{ij}=g_{ij},$ where $\mathcal{L=L}$ $(x,y)$ $\left( \Lambda
=\Lambda \left( x,y\right) \right) $ is a Lagrangian $\left(
\mbox{Finsler
metric}\right) $ on $TM$ (see details in \cite{ma87,ma94,mat,bej}).


\chapter{ Clifford Ha--Structures}

\index{Clifford ha--sructure}

The theory of anisotropic spinors formulated in the Part II is extended for
higher order anisotropic (ha) spaces. In brief, such spinors will be called
ha--spinors which are defined as some Clifford ha--structures defined with
respect to a distinguished quadratic form (\ref{dmetrichcv}) on a
hvc--bundle. For simplicity, the bulk of formulas will be given with respect
to higher order vector bundles. To rewrite such formulas for hvc--bundles is
to consider for the ''dual'' shells of higher order anisotropy some dual
vector spaces and associated dual spinors.

\section{Distinguished Clifford Algebras}

The typical fiber of dv--bundle $\xi _d\ ,\ \pi _d:\ HE\oplus V_1E\oplus
...\oplus V_zE\rightarrow E$ is a d-vector space, $\mathcal{F}=h\mathcal{F}%
\oplus v_1\mathcal{F\oplus }...\oplus v_z\mathcal{F},$ split into horizontal
$h\mathcal{F}$ and verticals $v_p\mathcal{F},p=1,...,z$ subspaces, with a
bilinear quadratic form $G(g,h)$ induced by a hvc--bundle metric (\ref%
{dmetrichcv}). Clifford algebras (see, for example, Refs. \cite%
{kar,tur,penr2}) formulated for d-vector spaces will be called Clifford
d-algebras \cite{vjmp,viasm1,vdeb}. We shall consider the main properties of
Clifford d--algebras. The proof of theorems will be based on the technique
developed in Ref. \cite{kar} correspondingly adapted to the distinguished
character of spaces in consideration.

Let $k$ be a number field (for our purposes $k=\mathcal{R}$ or $k=\mathcal{C}%
,\mathcal{R}$ and $\mathcal{C},$ are, respectively real and complex number
fields) and define $\mathcal{F},$ as a d-vector space on $k$ provided with
nondegenerate symmetric quadratic form (metric)\ $G.$ Let $C$ be an algebra
on $k$ (not necessarily commutative) and $j\ :\ \mathcal{F}$ $\rightarrow C$
a homomorphism of underlying vector spaces such that $j(u)^2=\;G(u)\cdot 1\
(1$ is the unity in algebra $C$ and d-vector $u\in \mathcal{F}).$ We are
interested in definition of the pair $\left( C,j\right) $ satisfying the
next universitality conditions. For every $k$-algebra $A$ and arbitrary
homomorphism $\varphi :\mathcal{F}\rightarrow A$ of the underlying d-vector
spaces, such that $\left( \varphi (u)\right) ^2\rightarrow G\left( u\right)
\cdot 1,$ there is a unique homomorphism of algebras $\psi \ :\ C\rightarrow
A$ transforming the diagram 1 into a commutative one.

The algebra solving this problem will be denoted as $C\left( \mathcal{F}%
,A\right) $ [equivalently as $C\left( G\right) $ or $C\left( \mathcal{F}%
\right) ]$ and called as Clifford d--algebra associated with pair
$\left( \mathcal{F},G\right).$

\begin{theorem}
The above-presented diagram has a unique solution $\left( C,j\right) $ up to
isomorphism.
\end{theorem}

\textbf{Proof:} (We adapt for d-algebras that of Ref. \cite{kar}, p. 127 and
extend for higher order anisotropies a similar proof presented in the Part
II). For a universal problem the uniqueness is obvious if we prove the
existence of solution $C\left( G\right) $ . To do this we use tensor algebra
$\mathcal{L}^{(F)}=\oplus \mathcal{L}_{qs}^{pr}\left( \mathcal{F}\right) $ =$%
\oplus _{i=0}^\infty T^i\left( \mathcal{F}\right) ,$ where $T^0\left(
\mathcal{F}\right) =k$ and $T^i\left( \mathcal{F}\right) =k$ and $T^i\left(
\mathcal{F}\right) =\mathcal{F}\otimes ...\otimes \mathcal{F}$ for $i>0.$
Let $I\left( G\right) $ be the bilateral ideal generated by elements of form
$\epsilon \left( u\right) =u\otimes u-G\left( u\right) \cdot 1$ where $u\in
\mathcal{F}$ and $1$ is the unity element of algebra $\mathcal{L}\left(
\mathcal{F}\right) .$ Every element from $I\left( G\right) $ can be written
as $\sum\nolimits_i\lambda _i\epsilon \left( u_i\right) \mu _i,$ where $%
\lambda _i,\mu _i\in \mathcal{L}(\mathcal{F})$ and $u_i\in \mathcal{F}.$ Let
$C\left( G\right) $ =$\mathcal{L}(\mathcal{F})/I\left( G\right) $ and define
$j:\mathcal{F}\rightarrow C\left( G\right) $ as the composition of
monomorphism $i:{\mathcal{F}\rightarrow L}^1(\mathcal{F})\subset \mathcal{L}(%
\mathcal{F})$ and projection $p:\mathcal{L}\left( \mathcal{F}\right)
\rightarrow C\left( G\right) .$ In this case pair $\left( C\left( G\right)
,j\right) $ is the solution of our problem. From the general properties of
tensor algebras the homomorphism $\varphi :\mathcal{F}\rightarrow A$ can be
extended to $\mathcal{L}(\mathcal{F})$ , i.e., the diagram 2 is commutative,
where $\rho $ is a monomorphism of algebras. Because $\left( \varphi \left(
u\right) \right) ^2=G\left( u\right) \cdot 1,$ then $\rho $ vanishes on
ideal $I\left( G\right) $ and in this case the necessary homomorphism $\tau $
is defined. As a consequence of uniqueness of $\rho ,$ the homomorphism $%
\tau $ is unique.

Tensor d--algebra $\mathcal{L}(\mathcal{F)}$ can be considered as a $%
\mathcal{Z}/2$ graded algebra. Really, let us in\-tro\-duce $\mathcal{L}%
^{(0)}(\mathcal{F})=\sum_{i=1}^\infty T^{2i}\left( \mathcal{F}\right) $ and $%
\mathcal{L}^{(1)}(\mathcal{F})=\sum_{i=1}^\infty T^{2i+1}\left( \mathcal{F}%
\right) .$ Setting $I^{(\alpha )}\left( G\right) =I\left( G\right) \cap
\mathcal{L}^{(\alpha )}(\mathcal{F}).$ Define $C^{(\alpha )}\left( G\right) $
as $p\left( \mathcal{L}^{(\alpha )}(\mathcal{F})\right) ,$ where $p:\mathcal{%
L}\left( \mathcal{F}\right) \rightarrow C\left( G\right) $ is the canonical
projection. Then $C\left( G\right) =C^{(0)}\left( G\right) \oplus
C^{(1)}\left( G\right) $ and in consequence we obtain that the Clifford
d--algebra is $\mathcal{Z}/2$ graded.

It is obvious that Clifford d-algebra
\index{Clifford d--algebra} functorially depends on pair $\left( \mathcal{F}%
,G\right) .$ If $f:\mathcal{F}\rightarrow\mathcal{F}^{\prime }$ is a
homomorphism of k-vector spaces, such that $G^{\prime }\left( f(u)\right)
=G\left( u\right) ,$ where $G$ and $G^{\prime }$ are, respectively, metrics
on $\mathcal{F}$ and $\mathcal{F}^{\prime },$ then $f$ induces an
homomorphism of d-algebras
\begin{equation*}
C\left( f\right) :C\left( G\right) \rightarrow C\left( G^{\prime }\right)
\end{equation*}
with identities $C\left( \varphi \cdot f\right) =C\left( \varphi \right)
C\left( f\right) $ and $C\left( Id_{\mathcal{F}}\right) =Id_{C(\mathcal{F)}%
}. $

If $\mathcal{A}^{\alpha}$ and $\mathcal{B}^{\beta}$ are $\mathcal{Z}/2$%
--graded d--algebras, then their graded tensorial product $\mathcal{A}%
^\alpha \otimes \mathcal{B}^\beta $ is defined as a d-algebra for k-vector
d-space $\mathcal{A}^\alpha \otimes \mathcal{B}^\beta $ with the graded
product induced as $\left( a\otimes b\right) \left( c\otimes d\right)
=\left( -1\right) ^{\alpha \beta }ac\otimes bd,$ where $b\in \mathcal{B}%
^\alpha $ and $c\in \mathcal{A}^\alpha \quad \left( \alpha ,\beta
=0,1\right) .$

Now we re--formulate for d--algebras the Chevalley theorem \cite{chev}:

\begin{theorem}
The Clifford d-algebra
\begin{equation*}
C\left( h\mathcal{F}\oplus v_1\mathcal{F}\oplus ...\oplus v_z\mathcal{F}%
,g+h_1+...+h_z\right)
\end{equation*}
is naturally isomorphic to $C(g)\otimes C\left( h_1\right) \otimes
...\otimes C\left( h_z\right) .$
\end{theorem}

\textbf{Proof. }Let $n:h\mathcal{F}\rightarrow C\left( g\right) $ and $%
n_{(p)}^{\prime }:v_{(p)}\mathcal{F}\rightarrow C\left( h_{(p)}\right) $ be
canonical maps and map
\begin{equation*}
m:h\mathcal{F}\oplus v_{1}\mathcal{F}\oplus ...\oplus v_{z}\mathcal{F}%
\rightarrow C(g)\otimes C\left( h_{1}\right) \otimes ...\otimes C\left(
h_{z}\right)
\end{equation*}%
is defined as
\begin{equation}
m(x,y_{(1)},...,y_{(z)})=n(x)\otimes 1\otimes ...\otimes 1+1\otimes
n^{\prime }(y_{(1)})\otimes ...\otimes 1+1\otimes ...\otimes 1\otimes
n^{\prime }(y_{(z)}),  \notag
\end{equation}%
$x\in h\mathcal{F},y_{(1)}\in v_{(1)}\mathcal{F},...,y_{(z)}\in v_{(z)}%
\mathcal{F}.$ We have
\begin{eqnarray*}
\left( m(x,y_{(1)},...,y_{(z)})\right) ^{2} &=&\left[ \left( n\left(
x\right) \right) ^{2}+\left( n^{\prime }\left( y_{(1)}\right) \right)
^{2}+...+\left( n^{\prime }\left( y_{(z)}\right) \right) ^{2}\right] \cdot 1
\\
&=&[g\left( x\right) +h\left( y_{(1)}\right) +...+h\left( y_{(z)}\right) ].
\end{eqnarray*}%
\ Taking into account the universality property of Clifford d-algebras we
conclude that $m_{1}+...+m_{z}$ induces the homomorphism
\begin{eqnarray}
\zeta :C\left( h\mathcal{F}\oplus v_{1}\mathcal{F}\oplus ...\oplus v_{z}%
\mathcal{F},g+h_{1}+...+h_{z}\right) \rightarrow \notag \\ C\left( h\mathcal{F}%
,g\right)
\widehat{\otimes }C\left( v_{1}\mathcal{F},h_{1}\right) \widehat{\otimes }%
...C\left( v_{z}\mathcal{F},h_{z}\right) .  \notag
\end{eqnarray}%
We also can define a homomorphism
\begin{eqnarray}
\upsilon :C\left( h\mathcal{F},g\right) \widehat{\otimes }C\left( v_{1}%
\mathcal{F},h_{(1)}\right) \widehat{\otimes }...\widehat{\otimes }C\left(
v_{z}\mathcal{F},h_{(z)}\right) \rightarrow \notag \\
 C\left( h\mathcal{F}\oplus v_{1} \mathcal{F}\oplus ...
 \oplus v_{z}\mathcal{F},g+h_{(1)}+...+h_{(z)}\right)
\notag
\end{eqnarray}%
by using formula $\upsilon \left( x\otimes y_{(1)}\otimes
...\otimes y_{(z)}\right) =\delta \left( x\right) \delta
_{(1)}^{\prime }\left( y_{(1)}\right) ...\delta _{(z)}^{\prime
}\left( y_{(z)}\right) ,$ where homomorphysms $\delta $ and
$\delta _{(1)}^{\prime },...,\delta _{(z)}^{\prime }$ are,
respectively, induced by embedding of $h\mathcal{F}$ and
$v_{1}\mathcal{F}$ into $h\mathcal{F}\oplus v_{1}\mathcal{F}\oplus
...\oplus v_{z}\mathcal{F}:$%
\begin{eqnarray*}
\delta &:&C\left( h\mathcal{F},g\right) \rightarrow C\left( h\mathcal{F}%
\oplus v_{1}\mathcal{F}\oplus ...\oplus v_{z}\mathcal{F}%
,g+h_{(1)}+...+h_{(z)}\right) , \\
\delta _{(1)}^{\prime } &:&C\left( v_{1}\mathcal{F},h_{(1)}\right)
\rightarrow C\left( h\mathcal{F}\oplus v_{1}\mathcal{F}\oplus ...\oplus v_{z}%
\mathcal{F},g+h_{(1)}+...+h_{(z)}\right) , \\
&&................................... \\
\delta _{(z)}^{\prime } &:&C\left( v_{z}\mathcal{F},h_{(z)}\right)
\rightarrow C\left( h\mathcal{F}\oplus v_{1}\mathcal{F}\oplus ...\oplus v_{z}%
\mathcal{F},g+h_{(1)}+...+h_{(z)}\right) .
\end{eqnarray*}

Superposition of homomorphisms $\zeta $ and $\upsilon $ lead to
identities
\begin{eqnarray}
\upsilon \zeta &=&Id_{C\left( h\mathcal{F},g\right) \widehat{\otimes }%
C\left( v_{1}\mathcal{F},h_{(1)}\right) \widehat{\otimes }...\widehat{%
\otimes }C\left( v_{z}\mathcal{F},h_{(z)}\right) },\zeta \upsilon
=Id_{C\left( h\mathcal{F},g\right) \widehat{\otimes }C\left( v_{1}\mathcal{F}%
,h_{(1)}\right) \widehat{\otimes }...\widehat{\otimes }C\left( v_{z}\mathcal{%
F},h_{(z)}\right) }.  \label{2.37a} \\
&&  \notag
\end{eqnarray}%
Really, d-algebra $C\left( h\mathcal{F}\oplus v_{1}\mathcal{F}\oplus
...\oplus v_{z}\mathcal{F},g+h_{(1)}+...+h_{(z)}\right) $ is generated by
elements of type $m(x,y_{(1)},...y_{(z)}).$ Calculating
\begin{eqnarray*}
&&\upsilon \zeta \left( m\left( x,y_{(1)},...y_{(z)}\right) \right)
=\upsilon (n\left( x\right) \otimes 1\otimes ...\otimes 1+1\otimes
n_{(1)}^{\prime }\left( y_{(1)}\right) \otimes ...\otimes 1 \\
&&+...+1\otimes ....\otimes n_{(z)}^{\prime }\left( y_{(z)}\right) )=\delta
\left( n\left( x\right) \right) \delta \left( n_{(1)}^{\prime }\left(
y_{(1)}\right) \right) ...\delta \left( n_{(z)}^{\prime }\left(
y_{(z)}\right) \right) \\
&=&m\left( x,0,...,0\right) +m(0,y_{(1)},...,0)+...+m(0,0,...,y_{(z)}) \\
&=&m\left( x,y_{(1)},...,y_{(z)}\right) ,
\end{eqnarray*}%
we prove the first identity in (\ref{2.37a}).

On the other hand, d-algebra
\begin{equation*}
C\left( h\mathcal{F},g\right) \widehat{\otimes }C\left( v_1\mathcal{F}%
,h_{(1)}\right) \widehat{\otimes }...\widehat{\otimes }C\left( v_z\mathcal{F}%
,h_{(z)}\right)
\end{equation*}
is generated by elements of type
\begin{equation*}
n\left( x\right) \otimes 1\otimes ...\otimes ,1\otimes n_{(1)}^{\prime
}\left( y_{(1)}\right) \otimes ...\otimes 1,...1\otimes ....\otimes
n_{(z)}^{\prime }\left( y_{(z)}\right) ,
\end{equation*}
we prove the second identity in (\ref{2.37a}).

Following from the above--mentioned properties of homomorphisms $\zeta $ and
$\upsilon $ we can assert that the natural isomorphism is explicitly
constructed.$\Box $

In consequence of the presented in this section Theorems we conclude that
all operations with Clifford d-algebras can be reduced to calculations for $%
C\left( h\mathcal{F},g\right) $ and\\ $C\left( v_{(p)}\mathcal{F}%
,h_{(p)}\right) $ which are usual Clifford algebras of dimension $2^n$ and,
respectively, $2^{m_p}$ \cite{kar,ati}.

Of special interest is the case when $k=\mathcal{R}$ and $\mathcal{F}$ is
isomorphic to vector space $\mathcal{R}^{p+q,a+b}$ provided with quadratic
form
\begin{equation*}
-x_1^2-...-x_p^2+x_{p+q}^2-y_1^2-...-y_a^2+...+y_{a+b}^2.
\end{equation*}
In this case, the Clifford algebra, denoted as $\left(
C^{p,q},C^{a,b}\right) ,\,$ is generated by symbols $%
e_1^{(x)},e_2^{(x)},...,e_{p+q}^{(x)},e_1^{(y)},e_2^{(y)},...,e_{a+b}^{(y)}$
satisfying properties
\begin{eqnarray*}
\left( e_i\right) ^2 &=&-1~\left( 1\leq i\leq p\right) ,\left( e_j\right)
^2=-1~\left( 1\leq j\leq a\right) , \\
\left( e_k\right) ^2 &=&1~(p+1\leq k\leq p+q), \\
\left( e_j\right) ^2 &=&1~(n+1\leq s\leq a+b),~e_ie_j=-e_je_i,~i\neq j.\,
\end{eqnarray*}

Explicit calculations of $C^{p,q}$ and $C^{a,b}$ are possible by using
isomorphisms \cite{kar,penr2}
\begin{eqnarray*}
C^{p+n,q+n} &\simeq &C^{p,q}\otimes M_2\left( \mathcal{R}\right) \otimes
...\otimes M_2\left( \mathcal{R}\right) \\
&\cong &C^{p,q}\otimes M_{2^n}\left( \mathcal{R}\right) \cong M_{2^n}\left(
C^{p,q}\right) ,
\end{eqnarray*}
where $M_s\left( A\right) $ denotes the ring of quadratic matrices of order $%
s$ with coefficients in ring $A.$ Here we write the simplest isomorphisms $%
C^{1,0}\simeq \mathcal{C},C^{0,1}\simeq \mathcal{R}\oplus \mathcal{R}$ and $%
C^{2,0}=\mathcal{H},$ where by $\mathcal{H}$ is denoted the body of
quaternions.

Now, we emphasize that higher order Lagrange and Finsler spaces, denoted $%
H^{2n}$--spaces, admit locally a structure of Clifford algebra on complex
vector spaces. Really, by using almost \ Hermitian structure $J_\alpha
^{\quad \beta }$ and considering complex space $\mathcal{C}^n$ with
nondegenarate quadratic form $\sum_{a=1}^n\left| z_a\right| ^2,~z_a\in
\mathcal{C}^2$ induced locally by metric (\ref{dmetrichcv}) (rewritten in
complex coordinates $z_a=x_a+iy_a)$ we define Clifford algebra $%
\overleftarrow{C}^n=\underbrace{\overleftarrow{C}^1\otimes ...\otimes
\overleftarrow{C}^1}_n,$ where $\overleftarrow{C}^1=\mathcal{C\otimes }_R%
\mathcal{C=C\oplus C}$ or in consequence, $\overleftarrow{C}^n\simeq
C^{n,0}\otimes _{\mathcal{R}}\mathcal{C}\approx C^{0,n}\otimes _{\mathcal{R}}%
\mathcal{C}.$ Explicit calculations lead to isomorphisms
\begin{equation*}
\overleftarrow{C}^2=C^{0,2}\otimes _{\mathcal{R}}\mathcal{C}\approx
M_2\left( \mathcal{R}\right) \otimes _{\mathcal{R}}\mathcal{C}\approx
M_2\left( \overleftarrow{C}^n\right) ,~C^{2p}\approx M_{2^p}\left( \mathcal{C%
}\right)
\end{equation*}
and
\begin{equation*}
\overleftarrow{C}^{2p+1}\approx M_{2^p}\left( \mathcal{C}\right) \oplus
M_{2^p}\left( \mathcal{C}\right) ,
\end{equation*}
which show that complex Clifford algebras, defined locally for $H^{2n}$%
-spaces, have periodicity 2 on $p.$

Considerations presented in the proof of theorem 2.2 show that map $j:%
\mathcal{F}\rightarrow C\left( \mathcal{F}\right) $ is monomorphic, so we
can identify space $\mathcal{F}$ with its image in $C\left( \mathcal{F}%
,G\right) ,$ denoted as $u\rightarrow \overline{u},$ if $u\in C^{(0)}\left(
\mathcal{F},G\right) ~\left( u\in C^{(1)}\left( \mathcal{F},G\right) \right)
;$ then $u=\overline{u}$ ( respectively, $\overline{u}=-u).$

\begin{definition}
\index{Clifford d--group} The set of elements $u\in C\left( G\right) ^{*},$
where $C\left( G\right) ^{*}$ denotes the multiplicative group of invertible
elements of $C\left( \mathcal{F},G\right) $ satisfying $%
\overline{u}\mathcal{F}u^{-1}\in \mathcal{F},$ is called the twisted
Clifford d-group, denoted as $\widetilde{\Gamma }\left( \mathcal{F}\right) .$
\end{definition}

Let $\widetilde{\rho }:\widetilde{\Gamma }\left( \mathcal{F}\right)
\rightarrow GL\left( \mathcal{F}\right) $ be the homorphism given by $%
u\rightarrow \rho \widetilde{u},$ where $\widetilde{\rho }_u\left( w\right) =%
\overline{u}wu^{-1}.$ We can verify that $\ker \widetilde{\rho }=\mathcal{R}%
^{*}$is a subgroup in $\widetilde{\Gamma }\left( \mathcal{F}\right) .$

The canonical map $j:\mathcal{F}\rightarrow C\left( \mathcal{F}\right) $ can
be interpreted as the linear map $\mathcal{F}\rightarrow C\left( \mathcal{F}%
\right) ^0$ satisfying the universal property of Clifford d-algebras. This
leads to a homomorphism of algebras, $C\left( \mathcal{F}\right) \rightarrow
C\left( \mathcal{F}\right) ^t,$ considered by an anti-involution of $C\left(
\mathcal{F}\right) $ and denoted as $u\rightarrow ~^tu.$ More exactly, if $%
u_1...u_n\in \mathcal{F,}$ then $t_u=u_n...u_1$ and $^t\overline{u}=%
\overline{^tu}=\left( -1\right) ^nu_n...u_1.$

\begin{definition}
\index{spinor norm} The spinor norm of arbitrary $u\in C\left( \mathcal{F}%
\right) $ is defined as\newline
$S\left( u\right) =~^t%
\overline{u}\cdot u\in C\left( \mathcal{F}\right) .$
\end{definition}

It is obvious that if $u,u^{\prime },u^{\prime \prime }\in
\widetilde{\Gamma }\left( \mathcal{F}\right) ,$ then
$S(u,u^{\prime })=S\left( u\right) S\left( u^{\prime }\right) $
and \newline $S\left( uu^{\prime }u^{\prime \prime }\right)
=S\left( u\right) S\left( u^{\prime }\right) S\left( u^{\prime
\prime }\right) .$ For $u,u^{\prime }\in \mathcal{F} S\left(
u\right) =-G\left( u\right) $ and  \\ $S\left( u,u^{\prime
}\right) =S\left( u\right) S\left( u^{\prime }\right) =S\left(
uu^{\prime }\right) .$

Let us introduce the orthogonal group $O\left( G\right) \subset GL\left(
G\right) $ defined by metric $G$ on $\mathcal{F}$ and denote sets
\begin{equation*}
SO\left( G\right) =\{u\in O\left( G\right) ,\det \left| u\right|
=1\},~Pin\left( G\right) =\{u\in \widetilde{\Gamma }\left( \mathcal{F}%
\right) ,S\left( u\right) =1\}
\end{equation*}
and $Spin\left( G\right) =Pin\left( G\right) \cap C^0\left( \mathcal{F}%
\right) .$ For ${\mathcal{F}\cong \mathcal{R}}^{n+m}$ we write $Spin\left(
n_E\right) .$ By straightforward calculations (see similar considerations in
Ref. \cite{kar}) we can verify the exactness of these sequences:
\begin{eqnarray*}
1 &\rightarrow &\mathcal{Z}/2\rightarrow Pin\left( G\right) \rightarrow
O\left( G\right) \rightarrow 1, \\
1 &\rightarrow &\mathcal{Z}/2\rightarrow Spin\left( G\right) \rightarrow
SO\left( G\right) \rightarrow 0, \\
1 &\rightarrow &\mathcal{Z}/2\rightarrow Spin\left( n_E\right) \rightarrow
SO\left( n_E\right) \rightarrow 1.
\end{eqnarray*}
We conclude this subsection by emphasizing that the spinor norm was defined
with respect to a quadratic form induced by a metric in dv--bundle $\mathcal{E%
}^{<z>}$. This approach differs from those presented in Refs. \cite{asa88}
and \cite{ono}.

\section{Clifford Ha--Bundles}

We shall consider two variants of generalization of spinor constructions
defined for d-vector spaces to the case of distinguished vector bundle
spaces enabled with the structure of N-connection. The first is to use the
extension to the category of vector bundles. The second is to define the
Clifford fibration associated with compatible linear d-connection and metric
$G$ on a dv--bundle. We shall analyze both variants.

\subsection{Clifford d--module structure in dv--bundles}

\index{Clifford d--module}

Because functor $\mathcal{F}\to C(\mathcal{F})$ is smooth we can extend it
to the category of vector bundles of type
\begin{equation*}
\xi ^{<z>}=\{\pi _d:HE^{<z>}\oplus V_1E^{<z>}\oplus ...\oplus
V_zE^{<z>}\rightarrow E^{<z>}\}.
\end{equation*}
Recall that by $\mathcal{F}$ we denote the typical fiber of such bundles.
For $\xi ^{<z>}$ we obtain a bundle of algebras, denoted as $C\left( \xi
^{<z>}\right) ,\,$ such that $C\left( \xi ^{<z>}\right) _u=C\left( \mathcal{F%
}_u\right) .$ Multiplication in every fiber defines a continuous
map
\begin{equation*}
C\left( \xi ^{<z>}\right) \times C\left( \xi ^{<z>}\right)
\rightarrow C\left( \xi ^{<z>}\right).
\end{equation*}
If $\xi ^{<z>}$ is a distinguished vector bundle on number field $k$%
,\thinspace \thinspace the structure of the\\ $C\left( \xi ^{<z>}\right) $%
-module, the d-module, the d-module, on $\xi ^{<z>}$ is given by
the continuous map\\ $C\left( \xi ^{<z>}\right) \times _E\xi
^{<z>}\rightarrow \xi
^{<z>}$ with every fiber $\mathcal{F}_u$ provided with the structure of the $%
C\left( \mathcal{F}_u\right) -$module, correlated with its $k$-module
structure, Because $\mathcal{F}\subset C\left( \mathcal{F}\right) ,$ we have
a fiber to fiber map $\mathcal{F}\times _E\xi ^{<z>}\rightarrow \xi ^{<z>},$
inducing on every fiber the map $\mathcal{F}_u\times _E\xi
_{(u)}^{<z>}\rightarrow \xi _{(u)}^{<z>}$ ($\mathcal{R}$-linear on the first
factor and $k$-linear on the second one ). Inversely, every such bilinear
map defines on $\xi ^{<z>}$ the structure of the $C\left( \xi ^{<z>}\right) $%
-module by virtue of universal properties of Clifford d--algebras.
Equivalently, the above--mentioned bilinear map defines a morphism of
v--bundles
\begin{equation*}
m:\xi ^{<z>}\rightarrow HOM\left( \xi ^{<z>},\xi ^{<z>}\right) \quad
[HOM\left( \xi ^{<z>},\xi ^{<z>}\right)
\end{equation*}
denotes the bundles of homomorphisms] when $\left( m\left( u\right) \right)
^2=G\left( u\right) $ on every point.

Vector bundles $\xi ^{<z>}$ provided with $C\left( \xi ^{<z>}\right) $%
--structures are objects of the category with morphisms being morphisms of
dv-bundles, which induce on every point $u\in \xi ^{<z>}$ morphisms of $%
C\left( \mathcal{F}_u\right) -$modules. This is a Banach category contained
in the category of finite-dimensional d-vector spaces on filed $k.$

Let us denote by $H^s\left( \mathcal{E}^{<z>},GL_{n_E}\left( \mathcal{R}%
\right) \right) ,\,$ where $n_E=n+m_1+...+m_z,\,$ the s-dimensional
cohomology group of the algebraic sheaf of germs of continuous maps of
dv-bundle $\mathcal{E}^{<z>}$ with group $GL_{n_E}\left( \mathcal{R}\right) $
the group of automorphisms of $\mathcal{R}^{n_E}\,$ (for the language of
algebraic topology see, for example, Refs. \cite{kar} and \cite{god}). We
shall also use the group $SL_{n_E}\left( \mathcal{R}\right) =\{A\subset
GL_{n_E}\left( \mathcal{R}\right) ,\det A=1\}.\,$ Here we point out that
cohomologies $H^s(M,Gr)$ characterize the class of a principal bundle $\pi
:P\rightarrow M $ on $M$ with structural group $Gr.$ Taking into account
that we deal with bundles distinguished by an N-connection we introduce into
consideration cohomologies $H^s\left( \mathcal{E}^{<z>},GL_{n_E}\left(
\mathcal{R}\right) \right) $ as distinguished classes (d-classes) of bundles
$\mathcal{E}^{<z>}$ provided with a global N-connection structure.

For a real vector bundle $\xi ^{<z>}$ on compact base $\mathcal{E}^{<z>}$ we
can define the orientation on $\xi ^{<z>}$ as an element $\alpha _d\in
H^1\left( \mathcal{E}^{<z>},GL_{n_E}\left( \mathcal{R}\right) \right) $
whose image on map
\begin{equation*}
H^1\left( \mathcal{E}^{<z>},SL_{n_E}\left( \mathcal{R}\right) \right)
\rightarrow H^1\left( \mathcal{E}^{<z>},GL_{n_E}\left( \mathcal{R}\right)
\right)
\end{equation*}
is the d-class of bundle $\mathcal{E}^{<z>}.$

\begin{definition}
The spinor structure on $\xi ^{<z>}$ is defined as an element\newline
$\beta _d\in H^1\left( \mathcal{E}^{<z>},Spin\left( n_E\right) \right) $
whose image in the composition
\begin{equation*}
H^1\left( \mathcal{E}^{<z>},Spin\left( n_E\right) \right) \rightarrow
H^1\left( \mathcal{E}^{<z>},SO\left( n_E\right) \right) \rightarrow
H^1\left( \mathcal{E}^{<z>},GL_{n_E}\left( \mathcal{R}\right) \right)
\end{equation*}
is the d-class of $\mathcal{E}^{<z>}.$
\end{definition}

The above definition of spinor structures can be re--formulated in terms of
principal bundles. Let $\xi ^{<z>}$ be a real vector bundle of rank n+m on a
compact base $\mathcal{E}^{<z>}.$ If there is a principal bundle $P_d$ with
structural group $SO( n_E ) $ or $Spin( n_E ) ],$ this bundle $\xi ^{<z>}$
can be provided with orientation (or spinor) structure. The bundle $P_d$ is
associated with element\newline
$\alpha _d\in H^1\left(\mathcal{E}^{<z>},SO(n_{<z>})\right) $ [or $\beta
_d\in H^1\left( \mathcal{E}^{<z>}, Spin\left( n_E\right) \right) .$

We remark that a real bundle is oriented if and only if its first
Stiefel--Whitney d--class vanishes,
\begin{equation*}
w_1\left( \xi _d\right) \in H^1\left( \xi ,\mathcal{Z}/2\right) =0,
\end{equation*}
where $H^1\left( \mathcal{E}^{<z>},\mathcal{Z}/2\right) $ is the first group
of Chech cohomology with coefficients in $\mathcal{Z}/2,$ Considering the
second Stiefel--Whitney class $w_2\left( \xi ^{<z>}\right) \in H^2\left(
\mathcal{E}^{<z>},\mathcal{Z}/2\right) $ it is well known that vector bundle
$\xi ^{<z>}$ admits the spinor structure if and only if $w_2\left( \xi
^{<z>}\right) =0.$ Finally, we emphasize that taking into account that base
space $\mathcal{E}^{<z>}$ is also a v-bundle, $p:E^{<z>}\rightarrow M,$ we
have to make explicit calculations in order to express cohomologies $%
H^s\left( \mathcal{E}^{<z>},GL_{n+m}\right) \,$ and $H^s\left( \mathcal{E}%
^{<z>},SO\left( n+m\right) \right) $ through cohomologies
\begin{equation*}
H^s\left( M,GL_n\right) ,H^s\left( M,SO\left( m_1\right) \right)
,...H^s\left( M,SO\left( m_z\right) \right) ,
\end{equation*}
which depends on global topological structures of spaces $M$ and $\mathcal{E}%
^{<z>}$ $.$ For general bundle and base spaces this requires a cumbersome
cohomological calculus.

\subsection{Clifford fibration}

\index{Clifford fibration}

Another way of defining the spinor structure is to use Clifford fibrations.
Consider the principal bundle with the structural group $Gr$ being a
subgroup of orthogonal group $O\left( G\right) ,$ where $G$ is a quadratic
nondegenerate form ) defined on the base (also being a bundle space) space $%
\mathcal{E}^{<z>}.$ The fibration associated to principal fibration $P\left(
\mathcal{E}^{<z>},Gr\right) $ with a typical fiber having Clifford algebra $%
C\left( G\right) $ is, by definition, the Clifford fibration $PC\left(
\mathcal{E}^{<z>},Gr\right) .$ We can always define a metric on the Clifford
fibration if every fiber is isometric to $PC\left( \mathcal{E}%
^{<z>},G\right) $ (this result is proved for arbitrary quadratic forms $G$
on pseudo--Riemannian bases \cite{tur}). If, additionally, $Gr\subset
SO\left( G\right) $ a global section can be defined on $PC\left( G\right) .$

Let $\mathcal{P}\left( \mathcal{E}^{<z>},Gr\right) $ be the set of principal
bundles with differentiable base $\mathcal{E}^{<z>}$ and structural group $%
Gr.$ If $g:Gr\rightarrow Gr^{\prime }$ is an homomorphism of Lie groups and $%
P( \mathcal{E}^{<z>},Gr)$ $\subset \mathcal{P}\left( \mathcal{E}%
^{<z>},Gr\right) $ (for simplicity in this subsection we shall denote
mentioned bundles and sets of bundles as $P,P^{\prime }$ and respectively, $%
\mathcal{P},\mathcal{P}^{\prime }),$ we can always construct a principal
bundle with the property that there is an homomorphism $f:P^{\prime
}\rightarrow P$ of principal bundles which can be projected to the identity
map of $\mathcal{E}^{<z>}$ and corresponds to isomorphism $g:Gr\rightarrow
Gr^{\prime }.$ If the inverse statement also holds, the bundle $P^{\prime }$
is called as the extension of $P$ associated to $g$ and $f$ is called the
extension homomorphism denoted as $%
\widetilde{g.}$

Now we can define distinguished spinor structures on bundle spaces
\index{spinor structure}.

\begin{definition}
Let $P\in \mathcal{P}\left( \mathcal{E}^{<z>},O\left( G\right) \right) $ be
a principal bundle. A distinguished spinor structure of $P,$ equivalently a
ds-structure of $\mathcal{E}^{<z>}$ is an extension $%
\widetilde{P}$ of $P$ associated to homomorphism $h:PinG\rightarrow O\left(
G\right) $ where $O\left( G\right) $ is the group of orthogonal rotations,
generated by metric $G,$ in bundle $\mathcal{E}^{<z>}.$
\end{definition}

So, if $\widetilde{P}$ is a spinor structure of the space $\mathcal{E}%
^{<z>}, $ then $\widetilde{P}\in \mathcal{P}\left( \mathcal{E}%
^{<z>},PinG\right).$

The definition of spinor structures on varieties was given in
Ref.\cite{cru1}. In Refs. \cite{cru2} and \cite{cru2} it is
proved that a necessary and sufficient condition for a space time
to be orientable is to admit a global field of orthonormalized
frames. We mention that spinor structures can be also defined on
varieties modeled on Banach spaces \cite{ana77}. As we have shown
similar constructions are possible for the cases when space time
has the structure of a v-bundle with an N-connection.

\begin{definition}
A special distinguished spinor structure, ds-structure,
\index{ds--structure} of principal bundle $P=P\left( \mathcal{E}%
^{<z>},SO\left( G\right) \right) $ is a principal bundle\newline
$%
\widetilde{P}=\widetilde{P}\left( \mathcal{E}^{<z>},SpinG\right) $ for which
a homomorphism of principal bundles $\widetilde{p}:\widetilde{P}\rightarrow
P,$ projected on the identity map of $\mathcal{E}^{<z>}$ and corresponding
to representation
\begin{equation*}
R:SpinG\rightarrow SO\left( G\right) ,
\end{equation*}
is defined.
\end{definition}

In the case when the base space variety is oriented, there is a natural
bijection between tangent spinor structures with a common base. For special
ds--structures we can define, as for any spinor structure, the concepts of
spin tensors, spinor connections, and spinor covariant derivations (see
Refs. \cite{viasm1,vdeb,vsp1}).

\section{Almost Complex Spinor Structures}

\index{complex spinors}

Almost complex structures are an important characteristic of $H^{2n}$-spaces
and of osculator bundles $Osc^{k=2k_1}(M),$ where $k_1=1,2,...$ . For
simplicity in this subsection we restrict our analysis to the case of $%
H^{2n} $-spaces. We can rewrite the almost Hermitian metric \cite{ma87,ma94}%
, $H^{2n}$-metric in complex form \cite{vjmp}:

\begin{equation}
G=H_{ab}\left( z,\xi \right) dz^a\otimes dz^b,  \label{2.38a}
\end{equation}
where
\begin{equation*}
z^a=x^a+iy^a,~%
\overline{z^a}=x^a-iy^a,~H_{ab}\left( z,\overline{z}\right) =g_{ab}\left(
x,y\right) \mid _{y=y\left( z,\overline{z}\right) }^{x=x\left( z,\overline{z}%
\right) },
\end{equation*}
and define almost complex spinor structures. For given metric (\ref{2.38a})
on $H^{2n}$-space there is always a principal bundle $P^U$ with unitary
structural group U(n) which allows us to transform $H^{2n}$-space into
v-bundle $\xi ^U\approx P^U\times _{U\left( n\right) }\mathcal{R}^{2n}.$
This statement will be proved after we introduce complex spinor structures
on oriented real vector bundles \cite{kar}.

Let us consider momentarily $k=\mathcal{C}$ and introduce into consideration
[instead of the group $Spin(n)]$ the group $Spin^c\times _{\mathcal{Z}%
/2}U\left( 1\right) $ being the factor group of the product $Spin(n)\times
U\left( 1\right) $ with the respect to equivalence
\begin{equation*}
\left( y,z\right) \sim \left( -y,-a\right) ,\quad y\in Spin(m).
\end{equation*}
This way we define the short exact sequence
\begin{equation}
1\rightarrow U\left( 1\right) \rightarrow Spin^c\left( n\right) \overset{S^c}%
{\to }SO\left( n\right) \rightarrow 1,  \label{2.39a}
\end{equation}
where $\rho ^c\left( y,a\right) =\rho ^c\left( y\right) .$ If $\lambda $ is
oriented , real, and rank $n,$ $\gamma $-bundle $\pi :E_\lambda \rightarrow
M^n,$ with base $M^n,$ the complex spinor structure, spin structure, on $%
\lambda $ is given by the principal bundle $P$ with structural group $%
Spin^c\left( m\right) $ and isomorphism $\lambda \approx P\times
_{Spin^c\left( n\right) }\mathcal{R}^n$ (see (\ref{2.39a})). For such
bundles the categorial equivalence can be defined as
\begin{equation}
\epsilon ^c:\mathcal{E}_{\mathcal{C}}^T\left( M^n\right) \rightarrow
\mathcal{E}_{\mathcal{C}}^\lambda \left( M^n\right) ,  \label{2.40a}
\end{equation}
where $\epsilon ^c\left( E^c\right) =P\bigtriangleup _{Spin^c\left( n\right)
}E^c$ is the category of trivial complex bundles on $M^n,\mathcal{E}_{%
\mathcal{C}}^\lambda \left( M^n\right) $ is the category of
complex v-bundles on $M^n$ with action of Clifford bundle\\
$C\left( \lambda \right) ,P\bigtriangleup _{Spin^c(n)}$ and $E^c$
is the factor space of the bundle product $P\times _ME^c$ with
respect to the equivalence $\left( p,e\right)
\sim \left( p\widehat{g}^{-1},\widehat{g}e\right) ,p\in P,e\in E^c,$ where $%
\widehat{g}\in Spin^c\left( n\right) $ acts on $E$ by via the imbedding $%
Spin\left( n\right) \subset C^{0,n}$ and the natural action $U\left(
1\right) \subset \mathcal{C}$ on complex v-bundle $\xi ^c,E^c=tot\xi ^c,$
for bundle $\pi ^c:E^c\rightarrow M^n.$

Now we return to the bundle $\xi =\mathcal{E}^{<1>}.$ A real v-bundle (not
being a spinor bundle) admits a complex spinor structure if and only if
there exist a homomorphism $\sigma :U\left( n\right) \rightarrow
Spin^c\left( 2n\right) $ making the diagram 3 commutative. The explicit
construction of $\sigma $ for arbitrary $\gamma $-bundle is given in Refs. %
\cite{kar} and \cite{ati}. For $H^{2n}$-spaces it is obvious that a diagram
similar to (\ref{2.40a}) can be defined for the tangent bundle $TM^n,$ which
directly points to the possibility of defining the $^cSpin$-structure on $%
H^{2n}$--spaces.

Let $\lambda $ be a complex, rank\thinspace $n,$ spinor bundle with
\begin{equation}
\tau :Spin^c\left( n\right) \times _{\mathcal{Z}/2}U\left( 1\right)
\rightarrow U\left( 1\right)  \label{2.41a}
\end{equation}
the homomorphism defined by formula $\tau \left( \lambda ,\delta \right)
=\delta ^2.$ For $P_s$ being the principal bundle with fiber $Spin^c\left(
n\right) $ we introduce the complex linear bundle $L\left( \lambda ^c\right)
=P_S\times _{Spin^c(n)}\mathcal{C}$ defined as the factor space of $%
P_S\times \mathcal{C}$ on equivalence relation

\begin{equation*}
\left( pt,z\right) \sim \left( p,l\left( t\right) ^{-1}z\right) ,
\end{equation*}
where $t\in Spin^c\left( n\right) .$ This linear bundle is associated to
complex spinor structure on $\lambda ^c.$

If $\lambda ^c$ and $\lambda ^{c^{\prime }}$ are complex spinor bundles, the
Whitney sum $\lambda ^c\oplus \lambda ^{c^{\prime }}$ is naturally provided
with the structure of the complex spinor bundle. This follows from the
holomorphism
\begin{equation}
\omega ^{\prime }:Spin^c\left( n\right) \times Spin^c\left( n^{\prime
}\right) \rightarrow Spin^c\left( n+n^{\prime }\right) ,  \label{2.42a}
\end{equation}
given by formula $\left[ \left( \beta ,z\right) ,\left( \beta
^{\prime },z^{\prime }\right) \right] \rightarrow \left[ \omega
\left( \beta ,\beta ^{\prime }\right) ,zz^{\prime }\right] ,$
where $\omega $ is the homomorphism making the diagram 4
commutative. Here, $z,z^{\prime }\in U\left( 1\right) .$ It is
obvious that $L\left( \lambda ^c\oplus \lambda ^{c^{\prime
}}\right) $ is isomorphic to $L\left( \lambda ^c\right) \otimes
L\left( \lambda ^{c^{\prime }}\right).$

We conclude this subsection by formulating our main result on complex spinor
structures for $H^{2n}$-spaces:

\begin{theorem}
Let $\lambda ^c$ be a complex spinor bundle of rank $n$ and $H^{2n}$-space
considered as a real vector bundle $\lambda ^c\oplus \lambda ^{c^{\prime }}$
provided with almost complex structure $J_{\quad \beta }^\alpha ;$
multiplication on $i$ is given by $\left(
\begin{array}{cc}
0 & -\delta _j^i \\
\delta _j^i & 0%
\end{array}
\right) $. Then, the diagram 5 is commutative up to isomorphisms $\epsilon
^c $ and $\widetilde{\epsilon }^c$ defined as in (\ref{2.40a}), $\mathcal{H}$
is functor $E^c\rightarrow E^c\otimes L\left( \lambda ^c\right) $ and $%
\mathcal{E}_{\mathcal{C}}^{0,2n}\left( M^n\right) $ is defined by functor $%
\mathcal{E}_{\mathcal{C}}\left( M^n\right) \rightarrow \mathcal{E}_{\mathcal{%
C}}^{0,2n}\left( M^n\right) $ given as correspondence $E^c\rightarrow
\Lambda \left( \mathcal{C}^n\right) \otimes E^c$ (which is a categorial
equivalence), $\Lambda \left( \mathcal{C}^n\right) $ is the exterior algebra
on $\mathcal{C}^n.$ $W$ is the real bundle $\lambda ^c\oplus \lambda
^{c^{\prime }}$ provided with complex structure.
\end{theorem}

\textbf{Proof: }We use composition of homomorphisms
\begin{equation*}
\mu :Spin^c\left( 2n\right) \overset{\pi }{\to }SO\left( n\right) \overset{r}%
{\to }U\left( n\right) \overset{\sigma }{\to }Spin^c\left( 2n\right) \times
_{\mathcal{Z}/2}U\left( 1\right) ,
\end{equation*}
commutative diagram 6 and introduce composition of homomorphisms
\begin{equation*}
\mu :Spin^c\left( n\right) \overset{\Delta }{\to }Spin^c\left( n\right)
\times Spin^c\left( n\right) \overset{{\omega }^c}{\to }Spin^c\left(
n\right) ,
\end{equation*}
where $\Delta $ is the diagonal homomorphism and $\omega ^c$ is defined as
in (\ref{2.42a}). Using homomorphisms (\ref{2.41a}) and ((\ref{2.42a})) we
obtain formula $\mu \left( t\right) =\mu \left( t\right) r\left( t\right) .$

Now consider bundle $P\times _{Spin^c\left( n\right) }Spin^c\left( 2n\right)
$ as the principal $Spin^c\left( 2n\right) $-bundle, associated to $M\oplus
M $ being the factor space of the product $P\times Spin^c\left( 2n\right) $
on the equivalence relation $\left( p,t,h\right) \sim \left( p,\mu \left(
t\right) ^{-1}h\right) .$ In this case the categorial equivalence (\ref%
{2.40a}) can be rewritten as
\begin{equation*}
\epsilon ^c\left( E^c\right) =P\times _{Spin^c\left( n\right) }Spin^c\left(
2n\right) \Delta _{Spin^c\left( 2n\right) }E^c
\end{equation*}
and seen as factor space of $P\times Spin^c\left( 2n\right) \times _ME^c$ on
equivalence relation
\begin{equation*}
\left( pt,h,e\right) \sim \left( p,\mu \left( t\right) ^{-1}h,e\right) %
\mbox{and}\left( p,h_1,h_2,e\right) \sim \left( p,h_1,h_2^{-1}e\right)
\end{equation*}
(projections of elements $p$ and $e$ coincides on base $M).$ Every element
of $\epsilon ^c\left( E^c\right) $ can be represented as $P\Delta
_{Spin^c\left( n\right) }E^c,$ i.e., as a factor space $P\Delta E^c$ on
equivalence relation $\left( pt,e\right) \sim \left( p,\mu ^c\left( t\right)
,e\right) ,$ when $t\in Spin^c\left( n\right) .$ The complex line bundle $%
L\left( \lambda ^c\right) $ can be interpreted as the factor space of\newline
$P\times _{Spin^c\left( n\right) }\mathcal{C}$ on equivalence relation $%
\left( pt,\delta \right) \sim \left( p,r\left( t\right) ^{-1}\delta \right)
. $

Putting $\left( p,e\right) \otimes \left( p,\delta \right) \left( p,\delta
e\right) $ we introduce morphism
\begin{equation*}
\epsilon ^c\left( E\right) \times L\left( \lambda ^c\right) \rightarrow
\epsilon ^c\left( \lambda ^c\right)
\end{equation*}
with properties
\begin{eqnarray*}
\left( pt,e\right) \otimes \left( pt,\delta \right) &\rightarrow &\left(
pt,\delta e\right) =\left( p,\mu ^c\left( t\right) ^{-1}\delta e\right) , \\
\left( p,\mu ^c\left( t\right) ^{-1}e\right) \otimes \left( p,l\left(
t\right) ^{-1}e\right) &\rightarrow &\left( p,\mu ^c\left( t\right) r\left(
t\right) ^{-1}\delta e\right)
\end{eqnarray*}
pointing to the fact that we have defined the isomorphism correctly and that
it is an isomorphism on every fiber. $\Box $

\chapter{Spinors and Ha--Spaces}

\section{ D--Spinor Techniques}

The purpose of this section is to show how a corresponding abstract spinor
technique entailing notational and calculations advantages can be developed
for arbitrary splits of dimensions of a d-vector space $\mathcal{F}=h%
\mathcal{F}\oplus v_1\mathcal{F}\oplus ...\oplus v_z\mathcal{F}$, where $%
\dim h\mathcal{F}=n$ and $\dim v_p\mathcal{F}=m_p.$ For convenience we shall
also present some necessary coordinate expressions.

The problem of a rigorous definition of spinors on locally
anisotropic spaces (d--spinors) was posed and solved
\cite{vjmp,viasm1,vsp1} in the framework of the formalism of
Clifford and spinor structures on v-bundles provided with
compatible nonlinear and distinguished connections and metric. We
introduced
d-spinors as corresponding objects of the Clifford d-algebra $\mathcal{C}%
\left( \mathcal{F},G\right) $, defined for a d-vector space $\mathcal{F}$ in
a standard manner (see, for instance, \cite{kar}) and proved that operations
with $\mathcal{C}\left( \mathcal{F},G\right) \ $ can be reduced to
calculations for $\mathcal{C}\left( h\mathcal{F},g\right) ,\mathcal{C}\left(
v_1\mathcal{F},h_1\right) ,...$ and $\mathcal{C}\left( v_z\mathcal{F}%
,h_z\right) ,$ which are usual Clifford algebras of respective dimensions $%
2^n,2^{m_1},...$ and $2^{m_z}$ (if it is necessary we can use quadratic
forms $g$ and $h_p$ correspondingly induced on $h\mathcal{F}$ and $v_p%
\mathcal{F}$ by a metric $\mathbf{G}$ (\ref{dmetrichcv})). Considering the
orthogonal subgroup $O\mathbf{\left( G\right) }\subset GL\mathbf{\left(
G\right) }$ defined by a metric $\mathbf{G}$ we can define the d-spinor norm
and parametrize d-spinors by ordered pairs of elements of Clifford algebras $%
\mathcal{C}\left( h\mathcal{F},g\right) $ and $\mathcal{C}\left( v_p\mathcal{%
F},h_p\right) ,p=1,2,...z.$ We emphasize that the splitting of a Clifford
d-algebra associated to a dv-bundle $\mathcal{E}^{<z>}$ is a straightforward
consequence of the global decomposition defining a N-connection structure in
$\mathcal{E}^{<z>}$.

In this subsection we shall omit detailed proofs which in most cases are
mechanical but rather tedious. We can apply the methods developed in \cite%
{pen,penr1,penr2,lue} in a straightforward manner on h- and v-subbundles in
order to verify the correctness of affirmations.

\subsection{Clifford d--algebra, d--spinors and d--twistors}

\index{d--spinor}
\index{d--twistor}

In order to relate the succeeding constructions with Clifford d-algebras %
\cite{vjmp,viasm1} we consider a la-frame decomposition of the metric (\ref%
{dmetrichcv}):
\begin{equation*}
G_{<\alpha ><\beta >}\left( u\right) =l_{<\alpha >}^{<%
\widehat{\alpha }>}\left( u\right) l_{<\beta >}^{<\widehat{\beta }>}\left(
u\right) G_{<\widehat{\alpha }><\widehat{\beta }>},
\end{equation*}
where the frame d-vectors and constant metric matrices are distinguished as

\begin{eqnarray*}
l_{<\alpha >}^{<\widehat{\alpha }>}\left( u\right) &=&\left(
\begin{array}{cccc}
l_j^{\widehat{j}}\left( u\right) & 0 & ... & 0 \\
0 & l_{a_1}^{\widehat{a}_1}\left( u\right) & ... & 0 \\
... & ... & ... & ... \\
0 & 0 & .. & l_{a_z}^{\widehat{a}_z}\left( u\right)%
\end{array}
\right) , \\
G_{<\widehat{\alpha }><\widehat{\beta }>} &=&\left(
\begin{array}{cccc}
g_{\widehat{i}\widehat{j}} & 0 & ... & 0 \\
0 & h_{\widehat{a}_1\widehat{b}_1} & ... & 0 \\
... & ... & ... & ... \\
0 & 0 & 0 & h_{\widehat{a}_z\widehat{b}_z}%
\end{array}
\right) ,
\end{eqnarray*}
$g_{\widehat{i}\widehat{j}}$ and $h_{\widehat{a}_1\widehat{b}_1},...,h_{%
\widehat{a}_z\widehat{b}_z}$ are diagonal matrices with $g_{\widehat{i}%
\widehat{i}}=$ $h_{\widehat{a}_1\widehat{a}_1}=...=h_{\widehat{a}_z\widehat{b%
}_z}=\pm 1.$

To generate Clifford d-algebras we start with matrix equations
\begin{equation}
\sigma _{<\widehat{\alpha }>}\sigma _{<\widehat{\beta }>}+\sigma _{<\widehat{%
\beta }>}\sigma _{<\widehat{\alpha }>}=-G_{<\widehat{\alpha }><\widehat{%
\beta }>}I,  \label{2.43a}
\end{equation}
where $I$ is the identity matrix, matrices $\sigma _{<\widehat{\alpha }%
>}\,(\sigma $-objects) act on a d-vector space $\mathcal{F}=h\mathcal{F}%
\oplus v_1\mathcal{F}\oplus ...\oplus v_z\mathcal{F}$ and theirs components
are distinguished as
\begin{equation}
\sigma _{<\widehat{\alpha }>}\,=\left\{ (\sigma _{<\widehat{\alpha }>})_{%
\underline{\beta }}^{\cdot \underline{\gamma }}=\left(
\begin{array}{cccc}
(\sigma _{\widehat{i}})_{\underline{j}}^{\cdot \underline{k}} & 0 & ... & 0
\\
0 & (\sigma _{\widehat{a}_1})_{\underline{b}_1}^{\cdot \underline{c}_1} & ...
& 0 \\
... & ... & ... & ... \\
0 & 0 & ... & (\sigma _{\widehat{a}_z})_{\underline{b}_z}^{\cdot \underline{c%
}_z}%
\end{array}
\right) \right\} ,  \label{2.44a}
\end{equation}
indices \underline{$\beta $},\underline{$\gamma $},... refer to spin spaces
of type $\mathcal{S}=S_{(h)}\oplus S_{(v_1)}\oplus ...\oplus S_{(v_z)}$ and
underlined Latin indices \underline{$j$},$\underline{k},...$ and $\underline{%
b}_1,\underline{c}_1,...,\underline{b}_z,\underline{c}_z...$ refer
respectively to h-spin space $\mathcal{S}_{(h)}$ and v$_p$-spin space $%
\mathcal{S}_{(v_p)},(p=1,2,...,z)\ $which are correspondingly associated to
a h- and v$_p$-decomposition of a dv-bundle $\mathcal{E}^{<z>}.$ The
irreducible algebra of matrices $\sigma _{<\widehat{\alpha }>}$ of minimal
dimension $N\times N,$ where $N=N_{(n)}+N_{(m_1)}+...+N_{(m_z)},$ $\dim
\mathcal{S}_{(h)}$=$N_{(n)}$ and $\dim \mathcal{S}_{(v_p)}$=$N_{(m_p)},$ has
these dimensions
\begin{eqnarray*}
{N_{(n)}} &=&{\left\{
\begin{array}{rl}
{\ 2^{(n-1)/2},} & n=2k+1 \\
{2^{n/2},\ } & n=2k;%
\end{array}
\right. }\quad , \\
\quad {N}_{(m_p)}{} &=&{}\left|
\begin{array}{cc}
2^{(m_p-1)/2}, & m_p=2k_p+1 \\
2^{m_p}, & m_p=2k_p%
\end{array}
\right| ,
\end{eqnarray*}
where $k=1,2,...,k_p=1,2,....$

The Clifford d-algebra is generated by sums on $n+1$ elements of form
\begin{equation}
A_1I+B^{\widehat{i}}\sigma _{\widehat{i}}+C^{\widehat{i}\widehat{j}}\sigma _{%
\widehat{i}\widehat{j}}+D^{\widehat{i}\widehat{j}\widehat{k}}\sigma _{%
\widehat{i}\widehat{j}\widehat{k}}+...  \label{2.45a}
\end{equation}
and sums of $m_p+1$ elements of form
\begin{equation*}
A_{2(p)}I+B^{\widehat{a}_p}\sigma _{\widehat{a}_p}+C^{\widehat{a}_p\widehat{b%
}_p}\sigma _{\widehat{a}_p\widehat{b}_p}+D^{\widehat{a}_p\widehat{b}_p%
\widehat{c}_p}\sigma _{\widehat{a}_p\widehat{b}_p\widehat{c}_p}+...
\end{equation*}
with antisymmetric coefficients%
\begin{equation*}
C^{\widehat{i}\widehat{j}}=C^{[\widehat{i}\widehat{j}]},C^{\widehat{a}_p%
\widehat{b}_p}=C^{[\widehat{a}_p\widehat{b}_p]},D^{\widehat{i}\widehat{j}%
\widehat{k}}=D^{[\widehat{i}\widehat{j}\widehat{k}]},D^{\widehat{a}_p%
\widehat{b}_p\widehat{c}_p}=D^{[\widehat{a}_p\widehat{b}_p\widehat{c}_p]},...
\end{equation*}
and matrices
\begin{equation*}
\sigma _{\widehat{i}\widehat{j}}=\sigma _{[\widehat{i}}\sigma _{\widehat{j}%
]},\sigma _{\widehat{a}_p\widehat{b}_p}=\sigma _{[\widehat{a}_p}\sigma _{%
\widehat{b}_p]},\sigma _{\widehat{i}\widehat{j}\widehat{k}}=\sigma _{[%
\widehat{i}}\sigma _{\widehat{j}}\sigma _{\widehat{k}]},....
\end{equation*}
Really, we have 2$^{n+1}$ coefficients $\left( A_1,C^{\widehat{i}\widehat{j}%
},D^{\widehat{i}\widehat{j}\widehat{k}},...\right) $ and 2$^{m_p+1}$
coefficients \\ $(A_{2(p)},C^{\widehat{a}_p\widehat{b}_p},D^{\widehat{a}_p%
\widehat{b}_p\widehat{c}_p},...)$ of the Clifford algebra on $\mathcal{F}$.

For simplicity, we shall present the necessary geometric constructions only
for h-spin spaces $\mathcal{S}_{(h)}$ of dimension $N_{(n)}.$ Considerations
for a v-spin space $\mathcal{S}_{(v)}$ are similar but with proper
characteristics for a dimension $N_{(m)}.$

In order to define the scalar (spinor) product on $\mathcal{S}_{(h)}$ we
introduce into consideration this finite sum (because of a finite number of
elements $\sigma _{\lbrack \widehat{i}\widehat{j}...\widehat{k}]}$):
\begin{eqnarray}
&&  \notag \\
&^{(\pm )}E_{\underline{k}\underline{m}}^{\underline{i}\underline{j}}=\delta
_{\underline{k}}^{\underline{i}}\delta _{\underline{m}}^{\underline{j}}+%
\frac{2}{1!}(\sigma _{\widehat{i}})_{\underline{k}}^{.\underline{i}}(\sigma
^{\widehat{i}})_{\underline{m}}^{.\underline{j}}+\frac{2^{2}}{2!}(\sigma _{%
\widehat{i}\widehat{j}})_{\underline{k}}^{.\underline{i}}(\sigma ^{\widehat{i%
}\widehat{j}})_{\underline{m}}^{.\underline{j}}&+\frac{2^{3}}{3!}(\sigma _{%
\widehat{i}\widehat{j}\widehat{k}})_{\underline{k}}^{.\underline{i}}(\sigma
^{\widehat{i}\widehat{j}\widehat{k}})_{\underline{m}}^{.\underline{j}}+...
\label{2.46a}
\end{eqnarray}%
which can be factorized as
\begin{equation}
^{(\pm )}E_{\underline{k}\underline{m}}^{\underline{i}\underline{j}}=N_{(n)}{%
\ }^{(\pm )}\epsilon _{\underline{k}\underline{m}}{\ }^{(\pm )}\epsilon ^{%
\underline{i}\underline{j}}\mbox{ for }n=2k  \label{2.47a}
\end{equation}%
and
\begin{eqnarray}
^{(+)}E_{\underline{k}\underline{m}}^{\underline{i}\underline{j}}
&=&2N_{(n)}\epsilon _{\underline{k}\underline{m}}\epsilon ^{\underline{i}%
\underline{j}},{\ }^{(-)}E_{\underline{k}\underline{m}}^{\underline{i}%
\underline{j}}=0\mbox{ for }n=3(mod4),  \label{2.48a} \\
^{(+)}E_{\underline{k}\underline{m}}^{\underline{i}\underline{j}} &=&0,{\ }%
^{(-)}E_{\underline{k}\underline{m}}^{\underline{i}\underline{j}%
}=2N_{(n)}\epsilon _{\underline{k}\underline{m}}\epsilon ^{\underline{i}%
\underline{j}}\mbox{ for }n=1(mod4).  \notag
\end{eqnarray}

Antisymmetry of $\sigma _{\widehat{i}\widehat{j}\widehat{k}...}$ and the
construction of the objects (\ref{2.45a})--(\ref{2.48a}) define the
properties of $\epsilon $-objects $^{(\pm )}\epsilon _{\underline{k}%
\underline{m}}$ and $\epsilon _{\underline{k}\underline{m}}$ which have an
eight-fold periodicity on $n$ (see details in \cite{penr2} and, with respect
to locally anisotropic spaces, \cite{vjmp}).

For even values of $n$ it is possible the decomposition of every h-spin
space $\mathcal{S}_{(h)}$ into irreducible h-spin spaces $\mathbf{S}_{(h)}$
and $\mathbf{S}_{(h)}^{\prime }$ (one considers splitting of h-indices, for
instance, \underline{$l$}$=L\oplus L^{\prime },\underline{m}=M\oplus
M^{\prime },...;$ for v$_p$-indices we shall write $\underline{a}%
_p=A_p\oplus A_p^{\prime },\underline{b}_p=B_p\oplus B_p^{\prime },...)$ and
defines new $\epsilon $-objects
\begin{equation}
\epsilon ^{\underline{l}\underline{m}}=\frac 12\left( ^{(+)}\epsilon ^{%
\underline{l}\underline{m}}+^{(-)}\epsilon ^{\underline{l}\underline{m}%
}\right) \mbox{ and }\widetilde{\epsilon }^{\underline{l}\underline{m}%
}=\frac 12\left( ^{(+)}\epsilon ^{\underline{l}\underline{m}}-^{(-)}\epsilon
^{\underline{l}\underline{m}}\right)  \label{2.49a}
\end{equation}
We shall omit similar formulas for $\epsilon $-objects with lower indices.

In general, the spinor $\epsilon $-objects should be defined for
every shell of an\-iso\-tro\-py  where instead of dimension $n$
we shall consider the dimensions $m_p$, $1\leq p\leq z,$ of
shells.

We define a d--spinor space $\mathcal{S}_{(n,m_{1})}\ $ as a
direct sum of a horizontal and a vertical spinor spaces, for
instance,
\begin{eqnarray*}
\mathcal{S}_{(8k,8k^{\prime })} &=&\mathbf{S}_{\circ }\oplus \mathbf{S}%
_{\circ }^{\prime }\oplus \mathbf{S}_{|\circ }\oplus \mathbf{S}_{|\circ
}^{\prime },\mathcal{S}_{(8k,8k^{\prime }+1)}\ =\mathbf{S}_{\circ }\oplus
\mathbf{S}_{\circ }^{\prime }\oplus \mathcal{S}_{|\circ }^{(-)},..., \\
\mathcal{S}_{(8k+4,8k^{\prime }+5)} &=&\mathbf{S}_{\triangle }\oplus \mathbf{%
S}_{\triangle }^{\prime }\oplus \mathcal{S}_{|\triangle }^{(-)},...
\end{eqnarray*}%
The scalar product on a $\mathcal{S}_{(n,m_{1})}\ $ is induced by
(corresponding to fixed values of $n$ and $m_{1}$ ) $\epsilon
$-objects  considered for h- and v$_{1}$-components. We present
also an example
for $\mathcal{S}_{(n,m_{1}+...+m_{z})}:$%
\begin{equation}
\mathcal{S}_{(8k+4,8k_{(1)}+5,...,8k_{(p)}+4,...8k_{(z)})}=[\mathbf{S}%
_{\triangle }\oplus \mathbf{S}_{\triangle }^{\prime }\oplus \mathcal{S}%
_{|(1)\triangle }^{(-)}\oplus ...\oplus \mathbf{S}_{|(p)\triangle }\oplus
\mathbf{S}_{|(p)\triangle }^{\prime }\oplus ...\oplus \mathbf{S}_{|(z)\circ
}\oplus \mathbf{S}_{|(z)\circ }^{\prime }.  \notag
\end{equation}

Having introduced d-spinors for dimensions $\left( n,m_1+...+m_z\right) $ we
can write out the generalization for ha--spaces of twistor equations \cite%
{penr1} by using the distinguished $\sigma $-objects (\ref{2.44a}):
\begin{equation}
(\sigma _{(<\widehat{\alpha }>})_{|\underline{\beta }|}^{..\underline{\gamma
}}\quad \frac{\delta \omega ^{\underline{\beta }}}{\delta u^{<\widehat{\beta
}>)}}=\frac 1{n+m_1+...+m_z}\quad G_{<\widehat{\alpha }><\widehat{\beta }%
>}(\sigma ^{\widehat{\epsilon }})_{\underline{\beta }}^{..\underline{\gamma }%
}\quad \frac{\delta \omega ^{\underline{\beta }}}{\delta u^{\widehat{%
\epsilon }}},  \label{2.56a}
\end{equation}
where $\left| \underline{\beta }\right| $ denotes that we do not consider
symmetrization on this index. The general solution of (\ref{2.56a}) on the
d-vector space $\mathcal{F}$ looks like as
\begin{equation}
\omega ^{\underline{\beta }}=\Omega ^{\underline{\beta }}+u^{<\widehat{%
\alpha }>}(\sigma _{<\widehat{\alpha }>})_{\underline{\epsilon }}^{..%
\underline{\beta }}\Pi ^{\underline{\epsilon }},  \label{2.57a}
\end{equation}
where $\Omega ^{\underline{\beta }}$ and $\Pi ^{\underline{\epsilon }}$ are
constant d-spinors. For fixed values of dimensions $n$ and $m=m_1+...m_z$ we
mast analyze the reduced and irreducible components of h- and v$_p$-parts of
equations (\ref{2.56a}) and their solutions (\ref{2.57a}) in order to find
the symmetry properties of a d-twistor $\mathbf{Z^\alpha \ }$ defined as a
pair of d-spinors
\begin{equation*}
\mathbf{Z}^\alpha =(\omega ^{\underline{\alpha }},\pi _{\underline{\beta }%
}^{\prime }),
\end{equation*}
where $\pi _{\underline{\beta }^{\prime }}=\pi _{\underline{\beta }^{\prime
}}^{(0)}\in {\widetilde{\mathcal{S}}}_{(n,m_1,...,m_z)}$ is a constant dual
d-spinor. The problem of definition of spinors and twistors on ha-spaces was
firstly considered in \cite{vdeb} (see also \cite{v87}) in connection with
the possibility to extend the equations (\ref{2.57a}) and theirs solutions (%
\ref{2.58a}), by using nearly autoparallel maps, on curved, locally
isotropic or anisotropic, spaces. We note that the definition of twistors
have been extended to higher order anisotropic spaces with trivial N-- and
d--connections.

\subsection{ Mutual transforms of d-tensors and d-spinors}

The spinor algebra for spaces of higher dimensions can not be considered as
a real alternative to the tensor algebra as for locally isotropic spaces of
dimensions $n=3,4$ \cite{pen,penr1,penr2}. The same holds true for
ha--spaces and we emphasize that it is not quite convenient to perform a
spinor calculus for dimensions $n,m>>4$. Nevertheless, the concept of
spinors is important for every type of spaces, we can deeply understand the
fundamental properties of geometical objects on ha--spaces, and we shall
consider in this subsection some questions concerning transforms of d-tensor
objects into d-spinor ones.

\subsection{ Transformation of d-tensors into d-spinors}

In order to pass from d-tensors to d-spinors we must use $\sigma $-objects (%
\ref{2.44a}) written in reduced or irreduced form \quad (in
dependence of fixed values of dimensions $n$ and $m$):

\begin{eqnarray}
&&(\sigma _{<\widehat{\alpha }>})_{\underline{\beta }}^{\cdot \underline{%
\gamma }},~(\sigma ^{<\widehat{\alpha }>})^{\underline{\beta }\underline{%
\gamma }},~(\sigma ^{<\widehat{\alpha }>})_{\underline{\beta }\underline{%
\gamma }},...,(\sigma _{<\widehat{a}>})^{\underline{b}\underline{c}%
},...,(\sigma _{\widehat{i}})_{\underline{j}\underline{k}},...,(\sigma _{<%
\widehat{a}>})^{AA^{\prime }},...,(\sigma ^{\widehat{i}})_{II^{\prime }},....
\label{2.58a} \\
&&  \notag
\end{eqnarray}%
It is obvious that contracting with corresponding $\sigma $-objects (\ref%
{2.58a}) we can introduce instead of d-tensors indices the d-spinor ones,
for instance,
\begin{equation*}
\omega ^{\underline{\beta }\underline{\gamma }}=(\sigma ^{<\widehat{\alpha }%
>})^{\underline{\beta }\underline{\gamma }}\omega _{<\widehat{\alpha }%
>},\quad \omega _{AB^{\prime }}=(\sigma ^{<\widehat{a}>})_{AB^{\prime
}}\omega _{<\widehat{a}>},\quad ...,\zeta _{\cdot \underline{j}}^{\underline{%
i}}=(\sigma ^{\widehat{k}})_{\cdot \underline{j}}^{\underline{i}}\zeta _{%
\widehat{k}},....
\end{equation*}%
For d-tensors containing groups of antisymmetric indices there is a more
simple procedure of theirs transforming into d-spinors because the objects
\begin{equation}
(\sigma _{\widehat{\alpha }\widehat{\beta }...\widehat{\gamma }})^{%
\underline{\delta }\underline{\nu }},\quad (\sigma ^{\widehat{a}\widehat{b}%
...\widehat{c}})^{\underline{d}\underline{e}},\quad ...,(\sigma ^{\widehat{i}%
\widehat{j}...\widehat{k}})_{II^{\prime }},\quad ...  \label{2.59a}
\end{equation}%
can be used for sets of such indices into pairs of d-spinor indices. Let us
enumerate some properties of $\sigma $-objects of type (\ref{2.59a}) (for
simplicity we consider only h-components having q indices $\widehat{i},%
\widehat{j},\widehat{k},...$ taking values from 1 to $n;$ the properties of v%
$_{p}$-components can be written in a similar manner with respect to indices
$\widehat{a}_{p},\widehat{b}_{p},\widehat{c}_{p}...$ taking values from 1 to
$m$):
\begin{equation}
(\sigma _{\widehat{i}...\widehat{j}})^{\underline{k}\underline{l}}%
\mbox{
 is\ }\left\{ \
\begin{array}{c}
\mbox{symmetric on }\underline{k},\underline{l}\mbox{ for
}n-2q\equiv 1,7~(mod~8); \\
\mbox{antisymmetric on }\underline{k},\underline{l}\mbox{ for
}n-2q\equiv 3,5~(mod~8)%
\end{array}%
\right\}  \label{2.60a}
\end{equation}%
for odd values of $n,$ and an object
\begin{equation*}
(\sigma _{\widehat{i}...\widehat{j}})^{IJ}~\left( (\sigma _{\widehat{i}...%
\widehat{j}})^{I^{\prime }J^{\prime }}\right)
\end{equation*}%
\begin{equation}
\mbox{ is\ }\left\{
\begin{array}{c}
\mbox{symmetric on }I,J~(I^{\prime },J^{\prime })\mbox{ for
}n-2q\equiv 0~(mod~8); \\
\mbox{antisymmetric on }I,J~(I^{\prime },J^{\prime })\mbox{ for
}n-2q\equiv 4~(mod~8)%
\end{array}%
\right\}  \label{2.61a}
\end{equation}%
or
\begin{equation}
(\sigma _{\widehat{i}...\widehat{j}})^{IJ^{\prime }}=\pm (\sigma _{\widehat{i%
}...\widehat{j}})^{J^{\prime }I}\{%
\begin{array}{c}
n+2q\equiv 6(mod8); \\
n+2q\equiv 2(mod8),%
\end{array}
\label{2.62a}
\end{equation}%
with vanishing of the rest of reduced components of the d-tensor $(\sigma _{%
\widehat{i}...\widehat{j}})^{\underline{k}\underline{l}}$ with prime/
unprime sets of indices.

\subsection{ Fundamental d--spinors}

We can transform every d--spinor $\xi ^{\underline{\alpha }}=\left( \xi ^{%
\underline{i}},\xi ^{\underline{a}_1},...,\xi ^{\underline{a}_z}\right) $
into a corresponding d-tensor. For simplicity, we consider this construction
only for a h-component $\xi ^{\underline{i}}$ on a h-space being of
dimension $n$. The values
\begin{equation}
\xi ^{\underline{\alpha }}\xi ^{\underline{\beta }}(\sigma ^{\widehat{i}...%
\widehat{j}})_{\underline{\alpha }\underline{\beta }}\quad \left( n%
\mbox{ is
odd}\right)  \label{2.63a}
\end{equation}
or
\begin{equation}
\xi ^I\xi ^J(\sigma ^{\widehat{i}...\widehat{j}})_{IJ}~\left( \mbox{or }\xi
^{I^{\prime }}\xi ^{J^{\prime }}(\sigma ^{\widehat{i}...\widehat{j}%
})_{I^{\prime }J^{\prime }}\right) ~\left( n\mbox{ is even}\right)
\label{2.64a}
\end{equation}
with a different number of indices $\widehat{i}...\widehat{j},$ taken
together, defines the h-spinor
\index{h--spinor} $\xi ^{%
\underline{i}}\,$ to an accuracy to the sign. We emphasize that it is
necessary to choose only those h-components of d-tensors (\ref{2.63a}) (or (%
\ref{2.64a})) which are symmetric on pairs of indices $\underline{\alpha }%
\underline{\beta }$ (or $IJ\,$ (or $I^{\prime }J^{\prime }$ )) and the
number $q$ of indices $\widehat{i}...\widehat{j}$ satisfies the condition
(as a respective consequence of the properties (\ref{2.60a}) and/ or (\ref%
{2.61a}), (\ref{2.62a}))
\begin{equation}
n-2q\equiv 0,1,7~(mod~8).  \label{2.65a}
\end{equation}
Of special interest is the case when
\begin{equation}
q=\frac 12\left( n\pm 1\right) ~\left( n\mbox{ is odd}\right)  \label{2.66a}
\end{equation}
or
\begin{equation}
q=\frac 12n~\left( n\mbox{ is even}\right) .  \label{2.67a}
\end{equation}
If all expressions (\ref{2.63a}) and/or (\ref{2.64a}) are zero for all
values of $q\,$ with the exception of one or two ones defined by the
conditions (\ref{2.65a}), (\ref{2.66a}) (or (\ref{2.67a})), the value $\xi ^{%
\widehat{i}}$ (or $\xi ^I$ ($\xi ^{I^{\prime }}))$ is called a fundamental
h-spinor. Defining in a similar manner the fundamental v-spinors
\index{v--spinor} we can introduce fundamental d-spinors as pairs of
fundamental h- and v-spinors. Here we remark that a h(v$_p$)-spinor $\xi ^{%
\widehat{i}}~(\xi ^{\widehat{a}_p})\,$ (we can also consider reduced
components) is always a fundamental one for $n(m)<7,$ which is a consequence
of (\ref{2.67a})).

\section{ Differential Geometry of Ha--Spinors}

This subsection is devoted to the differential geometry of d--spinors in
higher order anisotropic spaces. We shall use denotations of type
\begin{equation*}
v^{<\alpha >}=(v^i,v^{<a>})\in \sigma ^{<\alpha >}=(\sigma ^i,\sigma ^{<a>})
\end{equation*}
and
\begin{equation*}
\zeta ^{\underline{\alpha }_p}=(\zeta ^{\underline{i}_p},\zeta ^{\underline{a%
}_p})\in \sigma ^{\alpha _p}=(\sigma ^{i_p},\sigma ^{a_p})\,
\end{equation*}
for, respectively, elements of modules of d-vector and irreduced d-spinor
fields (see details in \cite{vjmp}). D-tensors and d-spinor tensors
(irreduced or reduced) will be interpreted as elements of corresponding $%
\mathcal{\sigma }$--modules, for instance,
\begin{equation*}
q_{~<\beta >...}^{<\alpha >}\in \mathcal{\sigma ^{<\alpha >}\mathbf{/}}%
^{\prime };[-0\mathcal{_{<\beta >}},\psi _{~\underline{\beta }_p\quad ...}^{%
\underline{\alpha }_p\quad \underline{\gamma }_p}\in \mathcal{\sigma }_{~%
\underline{\beta _p}\quad ...}^{\underline{\alpha }_p\quad \underline{\gamma
}_p}~,\xi _{\quad J_pK_p^{\prime }N_p^{\prime }}^{I_pI_p^{\prime }}\in
\mathcal{\sigma }_{\quad J_pK_p^{\prime }N_p^{\prime }}^{I_pI_p^{\prime
}}~,...
\end{equation*}

We can establish a correspondence between the higher order anisotropic
adapted to the N--connection metric $g_{\alpha \beta }$ (\ref{dmetrichcv})
and d-spinor metric $\epsilon _{\underline{\alpha }\underline{\beta }}$ ($%
\epsilon $-objects  for both h- and v$_{p}$-subspaces of $%
\mathcal{E}^{<z>}$) of a ha--space $\mathcal{E}^{<z>}$ by using the relation
\begin{equation}
g_{<\alpha ><\beta >}=-\frac{1}{N(n)+N(m_{1})+...+N(m_{z})}\times ((\sigma
_{(<\alpha >}(u))^{\underline{\alpha }\underline{\beta }}(\sigma _{<\beta
>)}(u))^{\underline{\delta }\underline{\gamma }})\epsilon _{\underline{%
\alpha }\underline{\gamma }}\epsilon _{\underline{\beta }\underline{\delta }%
},  \label{2.68a}
\end{equation}%
where
\begin{equation}
(\sigma _{<\alpha >}(u))^{\underline{\nu }\underline{\gamma }}=l_{<\alpha
>}^{<\widehat{\alpha }>}(u)(\sigma _{<\widehat{\alpha }>})^{<\underline{\nu }%
><\underline{\gamma }>},  \label{2.69a}
\end{equation}%
which is a consequence of formulas (\ref{2.43a})--(\ref{2.49a}). In brief we
can write (\ref{2.68a}) as
\begin{equation}
g_{<\alpha ><\beta >}=\epsilon _{\underline{\alpha }_{1}\underline{\alpha }%
_{2}}\epsilon _{\underline{\beta }_{1}\underline{\beta }_{2}}  \label{2.70a}
\end{equation}%
if the $\sigma $-objects are considered as a fixed structure, whereas $%
\epsilon $-objects are treated as caring the metric ''dynamics ''
, on higher order anisotropic space. This variant is used, for
instance, in the so-called 2-spinor geometry \cite{penr1,penr2}
and should be preferred if we have to make explicit the algebraic
symmetry properties of d-spinor objects by using metric
decomposition (\ref{2.70a}). An alternative way is to consider as
fixed the algebraic structure of $\epsilon $-objects and to use
variable components of $\sigma $-objects of type (\ref{2.69a}) for
developing a variational d-spinor approach to gravitational and
matter field interactions on ha-spaces (the spinor Ashtekar
variables \cite{ash} are introduced in this manner).

We note that a d--spinor metric
\begin{equation*}
\epsilon _{\underline{\nu }\underline{\tau }}=\left(
\begin{array}{cccc}
\epsilon _{\underline{i}\underline{j}} & 0 & ... & 0 \\
0 & \epsilon _{\underline{a}_1\underline{b}_1} & ... & 0 \\
... & ... & ... & ... \\
0 & 0 & ... & \epsilon _{\underline{a}_z\underline{b}_z}%
\end{array}
\right)
\end{equation*}
on the d-spinor space $\mathcal{S}=(\mathcal{S}_{(h)},\mathcal{S}%
_{(v_1)},...,\mathcal{S}_{(v_z)})$ can have symmetric or antisymmetric h (v$%
_p$) -components $\epsilon _{\underline{i}\underline{j}}$ ($\epsilon _{%
\underline{a}_p\underline{b}_p})$ , see $\epsilon $-objects. For
simplicity, in order to avoid cumbersome calculations connected
with
eight-fold periodicity on dimensions $n$ and $m_p$ of a ha-space $\mathcal{E}%
^{<z>},$ we shall develop a general d-spinor formalism only by using
irreduced spinor spaces $\mathcal{S}_{(h)}$ and $\mathcal{S}_{(v_p)}.$

\subsection{ D-covariant derivation on ha--spaces}

Let $\mathcal{E}^{<z>}$ be a ha-space. We define the action on a d-spinor of
a d-covariant operator
\begin{eqnarray*}
\bigtriangledown _{<\alpha >} &=&\left( \bigtriangledown
_{i},\bigtriangledown _{<a>}\right) =(\sigma _{<\alpha >})^{\underline{%
\alpha }_{1}\underline{\alpha }_{2}}\bigtriangledown _{^{\underline{\alpha }%
_{1}\underline{\alpha }_{2}}}=\left( (\sigma _{i})^{\underline{i}_{1}%
\underline{i}_{2}}\bigtriangledown _{^{\underline{i}_{1}\underline{i}%
_{2}}},~(\sigma _{<a>})^{\underline{a}_{1}\underline{a}_{2}}\bigtriangledown
_{^{\underline{a}_{1}\underline{a}_{2}}}\right) \\
&=&((\sigma _{i})^{\underline{i}_{1}\underline{i}_{2}}\bigtriangledown _{^{%
\underline{i}_{1}\underline{i}_{2}}},~(\sigma _{a_{1}})^{\underline{a}_{1}%
\underline{a}_{2}}\bigtriangledown _{(1)^{\underline{a}_{1}\underline{a}%
_{2}}},...,(\sigma _{a_{p}})^{\underline{a}_{1}\underline{a}%
_{2}}\bigtriangledown _{(p)^{\underline{a}_{1}\underline{a}%
_{2}}},...,(\sigma _{a_{z}})^{\underline{a}_{1}\underline{a}%
_{2}}\bigtriangledown _{(z)^{\underline{a}_{1}\underline{a}_{2}}})
\end{eqnarray*}%
(in brief, we shall write
\begin{equation*}
\bigtriangledown _{<\alpha >}=\bigtriangledown _{^{\underline{\alpha }_{1}%
\underline{\alpha }_{2}}}=\left( \bigtriangledown _{^{\underline{i}_{1}%
\underline{i}_{2}}},~\bigtriangledown _{(1)^{\underline{a}_{1}\underline{a}%
_{2}}},...,\bigtriangledown _{(p)^{\underline{a}_{1}\underline{a}%
_{2}}},...,\bigtriangledown _{(z)^{\underline{a}_{1}\underline{a}%
_{2}}}\right) )
\end{equation*}%
as maps
\begin{equation*}
\bigtriangledown _{{\underline{\alpha }}_{1}{\underline{\alpha }}_{2}}\ :\
\mathcal{\sigma }^{\underline{\beta }}\rightarrow \sigma _{<\alpha >}^{%
\underline{\beta }}=\sigma _{{\underline{\alpha }}_{1}{\underline{\alpha }}%
_{2}}^{\underline{\beta }}=
\end{equation*}%
\begin{equation*}
\left( \sigma _{i}^{\underline{\beta }}=\sigma _{{\underline{i}}_{1}{%
\underline{i}}_{2}}^{\underline{\beta }},\sigma _{(1)a_{1}}^{\underline{%
\beta }}=\sigma _{(1){\underline{\alpha }}_{1}{\underline{\alpha }}_{2}}^{%
\underline{\beta }},...,\sigma _{(p)a_{p}}^{\underline{\beta }}=\sigma _{(p){%
\underline{\alpha }}_{1}{\underline{\alpha }}_{2}}^{\underline{\beta }%
},...,\sigma _{(z)a_{z}}^{\underline{\beta }}=\sigma _{(z){\underline{\alpha
}}_{1}{\underline{\alpha }}_{2}}^{\underline{\beta }}\right)
\end{equation*}%
satisfying conditions
\begin{equation*}
\bigtriangledown _{<\alpha >}(\xi ^{\underline{\beta }}+\eta ^{\underline{%
\beta }})=\bigtriangledown _{<\alpha >}\xi ^{\underline{\beta }%
}+\bigtriangledown _{<\alpha >}\eta ^{\underline{\beta }},
\end{equation*}%
and
\begin{equation*}
\bigtriangledown _{<\alpha >}(f\xi ^{\underline{\beta }})=f\bigtriangledown
_{<\alpha >}\xi ^{\underline{\beta }}+\xi ^{\underline{\beta }%
}\bigtriangledown _{<\alpha >}f
\end{equation*}%
for every $\xi ^{\underline{\beta }},\eta ^{\underline{\beta }}\in \mathcal{%
\sigma ^{\underline{\beta }}}$ and $f$ being a scalar field on $\mathcal{E}%
^{<z>}.\mathcal{\ }$ It is also required that one holds the Leibnitz rule
\begin{equation*}
(\bigtriangledown _{<\alpha >}\zeta _{\underline{\beta }})\eta ^{\underline{%
\beta }}=\bigtriangledown _{<\alpha >}(\zeta _{\underline{\beta }}\eta ^{%
\underline{\beta }})-\zeta _{\underline{\beta }}\bigtriangledown _{<\alpha
>}\eta ^{\underline{\beta }}
\end{equation*}%
and that $\bigtriangledown _{<\alpha >}\,$ is a real operator, i.e. it
commuters with the operation of complex conjugation:
\begin{equation*}
\overline{\bigtriangledown _{<\alpha >}\psi _{\underline{\alpha }\underline{%
\beta }\underline{\gamma }...}}=\bigtriangledown _{<\alpha >}(\overline{\psi
}_{\underline{\alpha }\underline{\beta }\underline{\gamma }...}).
\end{equation*}

Let now analyze the question on uniqueness of action on
d--spinors of an operator $\bigtriangledown _{<\alpha >}$
satisfying necessary conditions . Denoting by $\bigtriangledown
_{<\alpha >}^{(1)}$ and $\bigtriangledown _{<\alpha >}$ two such
d-covariant operators we consider the map
\begin{equation}
(\bigtriangledown _{<\alpha >}^{(1)}-\bigtriangledown _{<\alpha >}):\mathcal{%
\sigma ^{\underline{\beta }}\rightarrow \sigma _{\underline{\alpha }_{1}%
\underline{\alpha }_{2}}^{\underline{\beta }}}.  \label{2.71a}
\end{equation}%
Because the action on a scalar $f$ of both operators $\bigtriangledown
_{\alpha }^{(1)}$ and $\bigtriangledown _{\alpha }$ must be identical, i.e.
\begin{equation*}
\bigtriangledown _{<\alpha >}^{(1)}f=\bigtriangledown _{<\alpha >}f,
\end{equation*}%
the action (\ref{2.71a}) on $f=\omega _{\underline{\beta }}\xi ^{\underline{%
\beta }}$ must be written as
\begin{equation*}
(\bigtriangledown _{<\alpha >}^{(1)}-\bigtriangledown _{<\alpha >})(\omega _{%
\underline{\beta }}\xi ^{\underline{\beta }})=0.
\end{equation*}%
In consequence we conclude that there is an element $\Theta _{\underline{%
\alpha }_{1}\underline{\alpha }_{2}\underline{\beta }}^{\quad \quad
\underline{\gamma }}\in \mathcal{\sigma }_{\underline{\alpha }_{1}\underline{%
\alpha }_{2}\underline{\beta }}^{\quad \quad \underline{\gamma }}$ for which
\begin{eqnarray}
&&  \notag \\
\bigtriangledown _{\underline{\alpha }_{1}\underline{\alpha }_{2}}^{(1)}\xi
^{\underline{\gamma }} &=&\bigtriangledown _{\underline{\alpha }_{1}%
\underline{\alpha }_{2}}\xi ^{\underline{\gamma }}+\Theta _{\underline{%
\alpha }_{1}\underline{\alpha }_{2}\underline{\beta }}^{\quad \quad
\underline{\gamma }}\xi ^{\underline{\beta }},\bigtriangledown _{\underline{%
\alpha }_{1}\underline{\alpha }_{2}}^{(1)}\omega _{\underline{\beta }%
}=\bigtriangledown _{\underline{\alpha }_{1}\underline{\alpha }_{2}}\omega _{%
\underline{\beta }}-\Theta _{\underline{\alpha }_{1}\underline{\alpha }_{2}%
\underline{\beta }}^{\quad \quad \underline{\gamma }}\omega _{\underline{%
\gamma }}~.  \label{2.72a}
\end{eqnarray}%
The action of the operator (\ref{2.71a}) on a d-vector $v^{<\beta >}=v^{%
\underline{\beta }_{1}\underline{\beta }_{2}}$ can be written by using
formula (\ref{2.72a}) for both indices $\underline{\beta }_{1}$ and $%
\underline{\beta }_{2}$ :
\begin{eqnarray*}
(\bigtriangledown _{<\alpha >}^{(1)}-\bigtriangledown _{<\alpha >})v^{%
\underline{\beta }_{1}\underline{\beta }_{2}} &=&\Theta _{<\alpha >%
\underline{\gamma }}^{\quad \underline{\beta }_{1}}v^{\underline{\gamma }%
\underline{\beta }_{2}}+\Theta _{<\alpha >\underline{\gamma }}^{\quad
\underline{\beta }_{2}}v^{\underline{\beta }_{1}\underline{\gamma }} \\
&=&(\Theta _{<\alpha >\underline{\gamma }_{1}}^{\quad \underline{\beta }%
_{1}}\delta _{\underline{\gamma }_{2}}^{\quad \underline{\beta }_{2}}+\Theta
_{<\alpha >\underline{\gamma }_{1}}^{\quad \underline{\beta }_{2}}\delta _{%
\underline{\gamma }_{2}}^{\quad \underline{\beta }_{1}})v^{\underline{\gamma
}_{1}\underline{\gamma }_{2}}=Q_{\ <\alpha ><\gamma >}^{<\beta >}v^{<\gamma
>},
\end{eqnarray*}%
where
\begin{equation}
Q_{\ <\alpha ><\gamma >}^{<\beta >}=Q_{\qquad \underline{\alpha }_{1}%
\underline{\alpha }_{2}~\underline{\gamma }_{1}\underline{\gamma }_{2}}^{%
\underline{\beta }_{1}\underline{\beta }_{2}}=\Theta _{<\alpha >\underline{%
\gamma }_{1}}^{\quad \underline{\beta }_{1}}\delta _{\underline{\gamma }%
_{2}}^{\quad \underline{\beta }_{2}}+\Theta _{<\alpha >\underline{\gamma }%
_{1}}^{\quad \underline{\beta }_{2}}\delta _{\underline{\gamma }_{2}}^{\quad
\underline{\beta }_{1}}.  \label{2.73a}
\end{equation}%
The d-commutator $\bigtriangledown _{\lbrack <\alpha >}\bigtriangledown
_{<\beta >]}$ defines the d-torsion. So, applying operators\\ $%
\bigtriangledown _{\lbrack <\alpha >}^{(1)}\bigtriangledown _{<\beta
>]}^{(1)}$ and $\bigtriangledown _{\lbrack <\alpha >}\bigtriangledown
_{<\beta >]}$ on $f=\omega _{\underline{\beta }}\xi ^{\underline{\beta }}$
we can write
\begin{equation*}
T_{\quad <\alpha ><\beta >}^{(1)<\gamma >}-T_{~<\alpha ><\beta >}^{<\gamma
>}=Q_{~<\beta ><\alpha >}^{<\gamma >}-Q_{~<\alpha ><\beta >}^{<\gamma >}
\end{equation*}%
with $Q_{~<\alpha ><\beta >}^{<\gamma >}$ from (\ref{2.73a}).

The action of operator $\bigtriangledown _{<\alpha >}^{(1)}$ on d-spinor
tensors of type $\chi _{\underline{\alpha }_1\underline{\alpha }_2\underline{%
\alpha }_3...}^{\qquad \quad \underline{\beta }_1\underline{\beta }_2...}$
must be constructed by using formula (\ref{2.72a}) for every upper index $%
\underline{\beta }_1\underline{\beta }_2...$ and formula (\ref{2.73a}) for
every lower index $\underline{\alpha }_1\underline{\alpha }_2\underline{%
\alpha }_3...$ .

\subsection{Infeld--van der Waerden co\-ef\-fi\-ci\-ents}

Let
\begin{equation*}
\delta _{\underline{\mathbf{\alpha }}}^{\quad \underline{\alpha }}=\left(
\delta _{\underline{\mathbf{1}}}^{\quad \underline{i}},\delta _{\underline{%
\mathbf{2}}}^{\quad \underline{i}},...,\delta _{\underline{\mathbf{N(n)}}%
}^{\quad \underline{i}},\delta _{\underline{\mathbf{1}}}^{\quad \underline{a}%
},\delta _{\underline{\mathbf{2}}}^{\quad \underline{a}},...,\delta _{%
\underline{\mathbf{N(m)}}}^{\quad \underline{i}}\right)
\end{equation*}%
be a d--spinor basis. The dual to it basis is denoted as
\begin{equation*}
\delta _{\underline{\alpha }}^{\quad \underline{\mathbf{\alpha }}}=\left(
\delta _{\underline{i}}^{\quad \underline{\mathbf{1}}},\delta _{\underline{i}%
}^{\quad \underline{\mathbf{2}}},...,\delta _{\underline{i}}^{\quad
\underline{\mathbf{N(n)}}},\delta _{\underline{i}}^{\quad \underline{\mathbf{%
1}}},\delta _{\underline{i}}^{\quad \underline{\mathbf{2}}},...,\delta _{%
\underline{i}}^{\quad \underline{\mathbf{N(m)}}}\right) .
\end{equation*}%
A d-spinor $\kappa ^{\underline{\alpha }}\in \mathcal{\sigma }$ $^{%
\underline{\alpha }}$ has components $\kappa ^{\underline{\mathbf{\alpha }}%
}=\kappa ^{\underline{\alpha }}\delta _{\underline{\alpha }}^{\quad
\underline{\mathbf{\alpha }}}.$ Taking into account that
\begin{equation*}
\delta _{\underline{\mathbf{\alpha }}}^{\quad \underline{\alpha }}\delta _{%
\underline{\mathbf{\beta }}}^{\quad \underline{\beta }}\bigtriangledown _{%
\underline{\alpha }\underline{\beta }}=\bigtriangledown _{\underline{\mathbf{%
\alpha }}\underline{\mathbf{\beta }}},
\end{equation*}%
we write out the components $\bigtriangledown _{\underline{\alpha }%
\underline{\beta }}$ $\kappa ^{\underline{\gamma }}$ as
\begin{eqnarray}
\delta _{\underline{\mathbf{\alpha }}}^{\quad \underline{\alpha }}~\delta _{%
\underline{\mathbf{\beta }}}^{\quad \underline{\beta }}~\delta _{\underline{%
\gamma }}^{\quad \underline{\mathbf{\gamma }}}~\bigtriangledown _{\underline{%
\alpha }\underline{\beta }}\kappa ^{\underline{\gamma }} &=& \notag \\
\delta _{%
\underline{\mathbf{\epsilon }}}^{\quad \underline{\tau }}~\delta _{%
\underline{\tau }}^{\quad \underline{\mathbf{\gamma }}}~\bigtriangledown _{%
\underline{\mathbf{\alpha }}\underline{\mathbf{\beta }}}\kappa ^{\underline{%
\mathbf{\epsilon }}}+\kappa ^{\underline{\mathbf{\epsilon }}}~\delta _{%
\underline{\epsilon }}^{\quad \underline{\mathbf{\gamma }}}~\bigtriangledown
_{\underline{\mathbf{\alpha }}\underline{\mathbf{\beta }}}\delta _{%
\underline{\mathbf{\epsilon }}}^{\quad \underline{\epsilon }%
} &=& \bigtriangledown _{\underline{\mathbf{\alpha }}\underline{\mathbf{\beta }}%
}\kappa ^{\underline{\mathbf{\gamma }}}+\kappa ^{\underline{\mathbf{\epsilon
}}}\gamma _{~\underline{\mathbf{\alpha }}\underline{\mathbf{\beta }}%
\underline{\mathbf{\epsilon }}}^{\underline{\mathbf{\gamma }}},
\label{2.74a}
\end{eqnarray}%
where the coordinate components of the d--spinor connection $\gamma _{~%
\underline{\mathbf{\alpha }}\underline{\mathbf{\beta }}\underline{\mathbf{%
\epsilon }}}^{\underline{\mathbf{\gamma }}}$ are defined as
\begin{equation}
\gamma _{~\underline{\mathbf{\alpha }}\underline{\mathbf{\beta }}\underline{%
\mathbf{\epsilon }}}^{\underline{\mathbf{\gamma }}}\doteq \delta _{%
\underline{\tau }}^{\quad \underline{\mathbf{\gamma }}}~\bigtriangledown _{%
\underline{\mathbf{\alpha }}\underline{\mathbf{\beta }}}\delta _{\underline{%
\mathbf{\epsilon }}}^{\quad \underline{\tau }}.  \label{2.75a}
\end{equation}%
We call the Infeld - van der Waerden d-symbols a set of $\sigma $-objects ($%
\sigma _{\mathbf{\alpha }})^{\underline{\mathbf{\alpha }}\underline{\mathbf{%
\beta }}}$ parametri\-zed with respect to a coordinate d-spinor basis.
Defining
\begin{equation*}
\bigtriangledown _{<\mathbf{\alpha >}}=(\sigma _{<\mathbf{\alpha >}})^{%
\underline{\mathbf{\alpha }}\underline{\mathbf{\beta }}}~\bigtriangledown _{%
\underline{\mathbf{\alpha }}\underline{\mathbf{\beta }}},
\end{equation*}%
introducing denotations
\begin{equation*}
\gamma ^{\underline{\mathbf{\gamma }}}{}_{<\mathbf{\alpha >\underline{\tau }}%
}\doteq \gamma ^{\underline{\mathbf{\gamma }}}{}_{\mathbf{\underline{\alpha }%
\underline{\beta }\underline{\tau }}}~(\sigma _{<\mathbf{\alpha >}})^{%
\underline{\mathbf{\alpha }}\underline{\mathbf{\beta }}}
\end{equation*}%
and using properties (\ref{2.74a}) we can write relations
\begin{eqnarray}
l_{<\mathbf{\alpha >}}^{<\alpha >}~\delta _{\underline{\beta }}^{\quad
\underline{\mathbf{\beta }}}~\bigtriangledown _{<\alpha >}\kappa ^{%
\underline{\beta }} &=&\bigtriangledown _{<\mathbf{\alpha >}}\kappa ^{%
\underline{\mathbf{\beta }}}+\kappa ^{\underline{\mathbf{\delta }}}~\gamma
_{~<\mathbf{\alpha >}\underline{\mathbf{\delta }}}^{\underline{\mathbf{\beta
}}},  \label{2.76a} \\
l_{<\mathbf{\alpha >}}^{<\alpha >}~\delta _{\underline{\mathbf{\beta }}%
}^{\quad \underline{\beta }}~\bigtriangledown _{<\alpha >}~\mu _{\underline{%
\beta }} &=&\bigtriangledown _{<\mathbf{\alpha >}}~\mu _{\underline{\mathbf{%
\beta }}}-\mu _{\underline{\mathbf{\delta }}}\gamma _{~<\mathbf{\alpha >}%
\underline{\mathbf{\beta }}}^{\underline{\mathbf{\delta }}}.  \notag
\end{eqnarray}%
for d-covariant derivations $~\bigtriangledown _{\underline{\alpha }}\kappa
^{\underline{\beta }}$ and $\bigtriangledown _{\underline{\alpha }}~\mu _{%
\underline{\beta }}.$

We can consider expressions similar to (\ref{2.76a}) for values having both
types of d-spinor and d-tensor indices, for instance,
\begin{eqnarray}
& & l_{<\mathbf{\alpha >}}^{<\alpha >}~l_{<\gamma >}^{<\mathbf{\gamma >}%
}~\delta _{\underline{\mathbf{\delta }}}^{\quad \underline{\delta }%
}~\bigtriangledown _{<\alpha >}\theta _{\underline{\delta }}^{~<\gamma >}=
\notag \\
& & \bigtriangledown _{<\mathbf{\alpha >}}\theta _{\underline{\mathbf{\delta
}}}^{~<\mathbf{\gamma >}}-\theta _{\underline{\mathbf{\epsilon }}}^{~<%
\mathbf{\gamma >}}\gamma _{~<\mathbf{\alpha >}\underline{\mathbf{\delta }}}^{%
\underline{\mathbf{\epsilon }}}+\theta _{\underline{\mathbf{\delta }}}^{~<%
\mathbf{\tau >}}~\Gamma _{\quad <\mathbf{\alpha ><\tau >}}^{~<\mathbf{\gamma
>}}  \notag
\end{eqnarray}
(we can prove this by a straightforward calculation).

Now we shall consider some possible relations between components of
d-connec\-ti\-ons $\gamma _{~<\mathbf{\alpha >}\underline{\mathbf{\delta }}%
}^{\underline{\mathbf{\epsilon }}}$ and $\Gamma _{\quad <\mathbf{\alpha
><\tau >}}^{~<\mathbf{\gamma >}}$ and derivations of $(\sigma _{<\mathbf{%
\alpha >}})^{\underline{\mathbf{\alpha }}\underline{\mathbf{\beta }}}$ . We
can write
\begin{eqnarray*}
\Gamma _{~<\mathbf{\beta ><\gamma >}}^{<\mathbf{\alpha >}} &=&l_{<\alpha
>}^{<\mathbf{\alpha >}}\bigtriangledown _{<\mathbf{\gamma >}}l_{<\mathbf{%
\beta >}}^{<\alpha >}=l_{<\alpha >}^{<\mathbf{\alpha >}}\bigtriangledown _{<%
\mathbf{\gamma >}}(\sigma _{<\mathbf{\beta >}})^{\underline{\epsilon }%
\underline{\tau }}l_{<\alpha >}^{<\mathbf{\alpha >}}\bigtriangledown _{<%
\mathbf{\gamma >}}((\sigma _{<\mathbf{\beta >}})^{\underline{\mathbf{%
\epsilon }}\underline{\mathbf{\tau }}}\delta _{\underline{\mathbf{\epsilon }}%
}^{~\underline{\epsilon }}\delta _{\underline{\mathbf{\tau }}}^{~\underline{%
\tau }}) \\
&=&l_{<\alpha >}^{<\mathbf{\alpha >}}\delta _{\underline{\mathbf{\alpha }}%
}^{~\underline{\alpha }}\delta _{\underline{\mathbf{\epsilon }}}^{~%
\underline{\epsilon }}\bigtriangledown _{<\mathbf{\gamma >}}(\sigma _{<%
\mathbf{\beta >}})^{\underline{\mathbf{\alpha }}\underline{\mathbf{\epsilon }%
}}+l_{<\alpha >}^{<\mathbf{\alpha >}}(\sigma _{<\mathbf{\beta >}})^{%
\underline{\mathbf{\epsilon }}\underline{\mathbf{\tau }}}(\delta _{%
\underline{\mathbf{\tau }}}^{~\underline{\tau }}\bigtriangledown _{<\mathbf{%
\gamma >}}\delta _{\underline{\mathbf{\epsilon }}}^{~\underline{\epsilon }%
}+\delta _{\underline{\mathbf{\epsilon }}}^{~\underline{\epsilon }%
}\bigtriangledown _{<\mathbf{\gamma >}}\delta _{\underline{\mathbf{\tau }}%
}^{~\underline{\tau }}) \\
&=&l_{\underline{\mathbf{\epsilon }}\underline{\mathbf{\tau }}}^{<\mathbf{%
\alpha >}}~\bigtriangledown _{<\mathbf{\gamma >}}(\sigma _{<\mathbf{\beta >}%
})^{\underline{\mathbf{\epsilon }}\underline{\mathbf{\tau }}}+l_{\underline{%
\mathbf{\mu }}\underline{\mathbf{\nu }}}^{<\mathbf{\alpha >}}\delta _{%
\underline{\epsilon }}^{~\underline{\mathbf{\mu }}}\delta _{\underline{\tau }%
}^{~\underline{\mathbf{\nu }}}(\sigma _{<\mathbf{\beta >}})^{\underline{%
\epsilon }\underline{\tau }}(\delta _{\underline{\mathbf{\tau }}}^{~%
\underline{\tau }}\bigtriangledown _{<\mathbf{\gamma >}}\delta _{\underline{%
\mathbf{\epsilon }}}^{~\underline{\epsilon }} \notag \\
 & & +\delta _{\underline{\mathbf{%
\epsilon }}}^{~\underline{\epsilon }}\bigtriangledown _{<\mathbf{\gamma >}%
}\delta _{\underline{\mathbf{\tau }}}^{~\underline{\tau }}),
\end{eqnarray*}%
where $l_{<\alpha >}^{<\mathbf{\alpha >}}=(\sigma _{\underline{\mathbf{%
\epsilon }}\underline{\mathbf{\tau }}})^{<\mathbf{\alpha >}}$ , from which
one follows
\begin{equation}
(\sigma _{<\mathbf{\alpha >}})^{\underline{\mathbf{\mu }}\underline{\mathbf{%
\nu }}}(\sigma _{\underline{\mathbf{\alpha }}\underline{\mathbf{\beta }}})^{<%
\mathbf{\beta >}}\Gamma _{~<\mathbf{\gamma ><\beta >}}^{<\mathbf{\alpha >}%
}=(\sigma _{\underline{\mathbf{\alpha }}\underline{\mathbf{\beta }}})^{<%
\mathbf{\beta >}}\bigtriangledown _{<\mathbf{\gamma >}}(\sigma _{<\mathbf{%
\alpha >}})^{\underline{\mathbf{\mu }}\underline{\mathbf{\nu }}}+\delta _{%
\underline{\mathbf{\beta }}}^{~\underline{\mathbf{\nu }}}\gamma _{~<\mathbf{%
\gamma >\underline{\alpha }}}^{\underline{\mathbf{\mu }}}+\delta _{%
\underline{\mathbf{\alpha }}}^{~\underline{\mathbf{\mu }}}\gamma _{~<\mathbf{%
\gamma >\underline{\beta }}}^{\underline{\mathbf{\nu }}}.  \notag
\end{equation}%
Connecting the last expression on \underline{$\mathbf{\beta }$} and
\underline{$\mathbf{\nu }$} and using an orthonormalized d-spinor basis when
$\gamma _{~<\mathbf{\gamma >\underline{\beta }}}^{\underline{\mathbf{\beta }}%
}=0$ (a consequence from (\ref{2.75a})) we have
\begin{equation}
\gamma _{~<\mathbf{\gamma >\underline{\alpha }}}^{\underline{\mathbf{\mu }}%
}=\frac{1}{N(n)+N(m_{1})+...+N(m_{z})}(\Gamma _{\quad <\mathbf{\gamma >~%
\underline{\alpha }\underline{\beta }}}^{\underline{\mathbf{\mu }}\underline{%
\mathbf{\beta }}}-(\sigma _{\underline{\mathbf{\alpha }}\underline{\mathbf{%
\beta }}})^{<\mathbf{\beta >}}\bigtriangledown _{<\mathbf{\gamma >}}(\sigma
_{<\mathbf{\beta >}})^{\underline{\mathbf{\mu }}\underline{\mathbf{\beta }}%
}),  \label{2.78a}
\end{equation}%
where
\begin{equation}
\Gamma _{\quad <\mathbf{\gamma >~\underline{\alpha }\underline{\beta }}}^{%
\underline{\mathbf{\mu }}\underline{\mathbf{\beta }}}=(\sigma _{<\mathbf{%
\alpha >}})^{\underline{\mathbf{\mu }}\underline{\mathbf{\beta }}}(\sigma _{%
\underline{\mathbf{\alpha }}\underline{\mathbf{\beta }}})^{\mathbf{\beta }%
}\Gamma _{~<\mathbf{\gamma ><\beta >}}^{<\mathbf{\alpha >}}.  \label{2.79a}
\end{equation}%
We also note here that, for instance, for the canonical and Berwald
connections and Christoffel d-symbols we can express d-spinor connection (%
\ref{2.79a}) through corresponding locally adapted derivations of components
of metric and N-connection by introducing corresponding coefficients instead
of $\Gamma _{~<\mathbf{\gamma ><\beta >}}^{<\mathbf{\alpha >}}$ in (\ref%
{2.79a}) and than in (\ref{2.78a}).

\subsection{ D--spinors of ha--space curvature and torsion}

The d-tensor indices of the commutator $\Delta _{<\alpha ><\beta >}$ can be
transformed into d-spinor ones:
\begin{eqnarray}
\Box _{\underline{\alpha }\underline{\beta }} &=&(\sigma ^{<\alpha ><\beta
>})_{\underline{\alpha }\underline{\beta }}\Delta _{\alpha \beta }=(\Box _{%
\underline{i}\underline{j}},\Box _{\underline{a}\underline{b}})=(\Box _{%
\underline{i}\underline{j}},\Box _{\underline{a}_{1}\underline{b}%
_{1}},...,\Box _{\underline{a}_{p}\underline{b}_{p}},...,\Box _{\underline{a}%
_{z}\underline{b}_{z}}),  \label{2.80a} \\
&&  \notag
\end{eqnarray}%
with h- and v$_{p}$-components,
\begin{equation*}
\Box _{\underline{i}\underline{j}}=(\sigma ^{<\alpha ><\beta >})_{\underline{%
i}\underline{j}}\Delta _{<\alpha ><\beta >}\mbox{ and }\Box _{\underline{a}%
\underline{b}}=(\sigma ^{<\alpha ><\beta >})_{\underline{a}\underline{b}%
}\Delta _{<\alpha ><\beta >},
\end{equation*}%
being symmetric or antisymmetric in dependence of corresponding
values of dimensions $n\,$ and $m_{p}$ (see eight-fold
parametizations. Considering the actions of operator
(\ref{2.80a}) on d-spinors $\pi ^{\underline{\gamma }}$ and $\mu
_{\underline{\gamma }}$ we
introduce the d-spinor curvature $X_{\underline{\delta }\quad \underline{%
\alpha }\underline{\beta }}^{\quad \underline{\gamma }}\,$ as to satisfy
equations
\begin{equation}
\Box _{\underline{\alpha }\underline{\beta }}\ \pi ^{\underline{\gamma }}=X_{%
\underline{\delta }\quad \underline{\alpha }\underline{\beta }}^{\quad
\underline{\gamma }}\pi ^{\underline{\delta }}\mbox{ and }\Box _{\underline{%
\alpha }\underline{\beta }}\ \mu _{\underline{\gamma }}=X_{\underline{\gamma
}\quad \underline{\alpha }\underline{\beta }}^{\quad \underline{\delta }}\mu
_{\underline{\delta }}.  \label{2.81a}
\end{equation}%
The gravitational d-spinor $\Psi _{\underline{\alpha }\underline{\beta }%
\underline{\gamma }\underline{\delta }}$ is defined by a corresponding
symmetrization of d-spinor indices:
\begin{equation}
\Psi _{\underline{\alpha }\underline{\beta }\underline{\gamma }\underline{%
\delta }}=X_{(\underline{\alpha }|\underline{\beta }|\underline{\gamma }%
\underline{\delta })}.  \label{2.82a}
\end{equation}%
We note that d-spinor tensors $X_{\underline{\delta }\quad \underline{\alpha
}\underline{\beta }}^{\quad \underline{\gamma }}$ and $\Psi _{\underline{%
\alpha }\underline{\beta }\underline{\gamma }\underline{\delta }}\,$ are
transformed into similar 2-spinor objects on locally isotropic spaces \cite%
{penr1,penr2} if we consider vanishing of the N-connection structure and a
limit to a locally isotropic space.

Putting $\delta _{\underline{\gamma }}^{\quad \mathbf{\underline{\gamma }}}$
instead of $\mu _{\underline{\gamma }}$ in (\ref{2.81a}) and using (\ref%
{2.82a}) we can express respectively the curvature and gravitational
d-spinors as
\begin{equation*}
X_{\underline{\gamma }\underline{\delta }\underline{\alpha }\underline{\beta
}}=\delta _{\underline{\delta }\underline{\mathbf{\tau }}}\Box _{\underline{%
\alpha }\underline{\beta }}\delta _{\underline{\gamma }}^{\quad \mathbf{%
\underline{\tau }}}\mbox{ and }\Psi _{\underline{\gamma }\underline{\delta }%
\underline{\alpha }\underline{\beta }}=\delta _{\underline{\delta }%
\underline{\mathbf{\tau }}}\Box _{(\underline{\alpha }\underline{\beta }%
}\delta _{\underline{\gamma })}^{\quad \mathbf{\underline{\tau }}}.
\end{equation*}

The d-spinor torsion $T_{\qquad \underline{\alpha }\underline{\beta }}^{%
\underline{\gamma }_1\underline{\gamma }_2}$ is defined similarly as for
d-tensors by using the d-spinor commutator (\ref{2.80a}) and equations
\begin{equation*}
\Box _{\underline{\alpha }\underline{\beta }}f=T_{\qquad \underline{\alpha }%
\underline{\beta }}^{\underline{\gamma }_1\underline{\gamma }%
_2}\bigtriangledown _{\underline{\gamma }_1\underline{\gamma }_2}f.
\end{equation*}

The d-spinor components $R_{\underline{\gamma }_1\underline{\gamma }_2\qquad
\underline{\alpha }\underline{\beta }}^{\qquad \underline{\delta }_1%
\underline{\delta }_2}$ of the curvature d-tensor $R_{\gamma \quad \alpha
\beta }^{\quad \delta }$ can be computed by using relations (\ref{2.79a}),
and (\ref{2.80a}) and (\ref{2.82a}) as to satisfy the equations
\begin{equation*}
(\Box _{\underline{\alpha }\underline{\beta }}-T_{\qquad \underline{\alpha }%
\underline{\beta }}^{\underline{\gamma }_1\underline{\gamma }%
_2}\bigtriangledown _{\underline{\gamma }_1\underline{\gamma }_2})V^{%
\underline{\delta }_1\underline{\delta }_2}=R_{\underline{\gamma }_1%
\underline{\gamma }_2\qquad \underline{\alpha }\underline{\beta }}^{\qquad
\underline{\delta }_1\underline{\delta }_2}V^{\underline{\gamma }_1%
\underline{\gamma }_2},
\end{equation*}
here d-vector $V^{\underline{\gamma }_1\underline{\gamma }_2}$ is considered
as a product of d-spinors, i.e. $V^{\underline{\gamma }_1\underline{\gamma }%
_2}=\nu ^{\underline{\gamma }_1}\mu ^{\underline{\gamma }_2}$. We find
\begin{eqnarray}
R_{\underline{\gamma }_{1}\underline{\gamma }_{2}\qquad \underline{\alpha }%
\underline{\beta }}^{\qquad \underline{\delta }_{1}\underline{\delta }_{2}}
&=&\left( X_{\underline{\gamma }_{1}~\underline{\alpha }\underline{\beta }%
}^{\quad \underline{\delta }_{1}}+T_{\qquad \underline{\alpha }\underline{%
\beta }}^{\underline{\tau }_{1}\underline{\tau }_{2}}\quad \gamma _{\quad
\underline{\tau }_{1}\underline{\tau }_{2}\underline{\gamma }_{1}}^{%
\underline{\delta }_{1}}\right) \delta _{\underline{\gamma
}_{2}}^{\quad \underline{\delta }_{2}} \notag \\ +\left(
X_{\underline{\gamma }_{2}~\underline{\alpha
}\underline{\beta }}^{\quad \underline{\delta }_{2}}+T_{\qquad \underline{%
\alpha }\underline{\beta }}^{\underline{\tau }_{1}\underline{\tau }%
_{2}}\quad \gamma _{\quad \underline{\tau }_{1}\underline{\tau }_{2}%
\underline{\gamma }_{2}}^{\underline{\delta }_{2}}\right) \delta _{%
\underline{\gamma }_{1}}^{\quad \underline{\delta }_{1}}.
\label{2.83a}
\end{eqnarray}

It is convenient to use this d-spinor expression for the curvature d-tensor
\begin{eqnarray*}
R_{\underline{\gamma }_1\underline{\gamma }_2\qquad \underline{\alpha }_1%
\underline{\alpha }_2\underline{\beta }_1\underline{\beta }_2}^{\qquad
\underline{\delta }_1\underline{\delta }_2} &=&\left( X_{\underline{\gamma }%
_1~\underline{\alpha }_1\underline{\alpha }_2\underline{\beta }_1\underline{%
\beta }_2}^{\quad \underline{\delta }_1}+T_{\qquad \underline{\alpha }_1%
\underline{\alpha }_2\underline{\beta }_1\underline{\beta }_2}^{\underline{%
\tau }_1\underline{\tau }_2}~\gamma _{\quad \underline{\tau }_1\underline{%
\tau }_2\underline{\gamma }_1}^{\underline{\delta }_1}\right) \delta _{%
\underline{\gamma }_2}^{\quad \underline{\delta }_2} \\
&&+\left( X_{\underline{\gamma }_2~\underline{\alpha }_1\underline{\alpha }_2%
\underline{\beta }_1\underline{\beta }_2}^{\quad \underline{\delta }%
_2}+T_{\qquad \underline{\alpha }_1\underline{\alpha }_2\underline{\beta }_1%
\underline{\beta }_2~}^{\underline{\tau }_1\underline{\tau }_2}\gamma
_{\quad \underline{\tau }_1\underline{\tau }_2\underline{\gamma }_2}^{%
\underline{\delta }_2}\right) \delta _{\underline{\gamma }_1}^{\quad
\underline{\delta }_1}
\end{eqnarray*}
in order to get the d--spinor components of the Ricci d-tensor
\begin{eqnarray}
& &R_{\underline{\gamma }_1\underline{\gamma }_2\underline{\alpha }_1%
\underline{\alpha }_2} = R_{\underline{\gamma }_1\underline{\gamma }_2\qquad
\underline{\alpha }_1\underline{\alpha }_2\underline{\delta }_1\underline{%
\delta }_2}^{\qquad \underline{\delta }_1\underline{\delta }_2}=X_{%
\underline{\gamma }_1~\underline{\alpha }_1\underline{\alpha }_2\underline{%
\delta }_1\underline{\gamma }_2}^{\quad \underline{\delta }_1}+
\label{2.84a} \\
&&T_{\qquad \underline{\alpha }_1\underline{\alpha }_2\underline{\delta }_1%
\underline{\gamma }_2}^{\underline{\tau }_1\underline{\tau }_2}~\gamma
_{\quad \underline{\tau }_1\underline{\tau }_2\underline{\gamma }_1}^{%
\underline{\delta }_1}+X_{\underline{\gamma }_2~\underline{\alpha }_1%
\underline{\alpha }_2\underline{\delta }_1\underline{\gamma }_2}^{\quad
\underline{\delta }_2}+T_{\qquad \underline{\alpha }_1\underline{\alpha }_2%
\underline{\gamma }_1\underline{\delta }_2~}^{\underline{\tau }_1\underline{%
\tau }_2}\gamma _{\quad \underline{\tau }_1\underline{\tau }_2\underline{%
\gamma }_2}^{\underline{\delta }_2}  \notag
\end{eqnarray}
and this d-spinor decomposition of the scalar curvature:
\begin{eqnarray}
q\overleftarrow{R} &=&R_{\qquad \underline{\alpha }_1\underline{\alpha }_2}^{%
\underline{\alpha }_1\underline{\alpha }_2}=X_{\quad ~\underline{~\alpha }%
_1\quad \underline{\delta }_1\underline{\alpha }_2}^{\underline{\alpha }_1%
\underline{\delta }_1~~\underline{\alpha }_2}+T_{\qquad ~~\underline{\alpha }%
_2\underline{\delta }_1}^{\underline{\tau }_1\underline{\tau }_2\underline{%
\alpha }_1\quad ~\underline{\alpha }_2}~\gamma _{\quad \underline{\tau }_1%
\underline{\tau }_2\underline{\alpha }_1}^{\underline{\delta }_1}
\label{2.85a} \\
&&+X_{\qquad \quad \underline{\alpha }_2\underline{\delta }_2\underline{%
\alpha }_1}^{\underline{\alpha }_2\underline{\delta }_2\underline{\alpha }%
_1}+T_{\qquad \underline{\alpha }_1\quad ~\underline{\delta }_2~}^{%
\underline{\tau }_1\underline{\tau }_2~~\underline{\alpha }_2\underline{%
\alpha }_1}\gamma _{\quad \underline{\tau }_1\underline{\tau }_2\underline{%
\alpha }_2}^{\underline{\delta }_2}.  \notag
\end{eqnarray}

Putting (\ref{2.84a}) and (\ref{2.85a}) into (\ref{2.34}) and,
correspondingly, (\ref{2.35a}) we find the d--spinor components of the
Einstein and $\Phi _{<\alpha ><\beta >}$ d--tensors:
\begin{eqnarray}
\overleftarrow{G}_{<\gamma ><\alpha >} &=&\overleftarrow{G}_{\underline{%
\gamma }_1\underline{\gamma }_2\underline{\alpha }_1\underline{\alpha }%
_2}=X_{\underline{\gamma }_1~\underline{\alpha }_1\underline{\alpha }_2%
\underline{\delta }_1\underline{\gamma }_2}^{\quad \underline{\delta }%
_1}+T_{\qquad \underline{\alpha }_1\underline{\alpha }_2\underline{\delta }_1%
\underline{\gamma }_2}^{\underline{\tau }_1\underline{\tau }_2}~\gamma
_{\quad \underline{\tau }_1\underline{\tau }_2\underline{\gamma }_1}^{%
\underline{\delta }_1}  \notag \\
&&+X_{\underline{\gamma }_2~\underline{\alpha }_1\underline{\alpha }_2%
\underline{\delta }_1\underline{\gamma }_2}^{\quad \underline{\delta }%
_2}+T_{\qquad \underline{\alpha }_1\underline{\alpha }_2\underline{\gamma }_1%
\underline{\delta }_2~}^{\underline{\tau }_1\underline{\tau }_2}\gamma
_{\quad \underline{\tau }_1\underline{\tau }_2\underline{\gamma }_2}^{%
\underline{\delta }_2}-  \notag \\
&&\frac 12\varepsilon _{\underline{\gamma }_1\underline{\alpha }%
_1}\varepsilon _{\underline{\gamma }_2\underline{\alpha }_2}[X_{\quad ~%
\underline{~\beta }_1\quad \underline{\mu }_1\underline{\beta }_2}^{%
\underline{\beta }_1\underline{\mu }_1~~\underline{\beta }_2}+T_{\qquad ~~%
\underline{\beta }_2\underline{\mu }_1}^{\underline{\tau }_1\underline{\tau }%
_2\underline{\beta }_1\quad ~\underline{\beta }_2}~\gamma _{\quad \underline{%
\tau }_1\underline{\tau }_2\underline{\beta }_1}^{\underline{\mu }_1}+
\notag \\
&&X_{\qquad \quad \underline{\beta }_2\underline{\mu }_2\underline{\nu }_1}^{%
\underline{\beta }_2\underline{\mu }_2\underline{\nu }_1}+T_{\qquad
\underline{\beta }_1\quad ~\underline{\delta }_2~}^{\underline{\tau }_1%
\underline{\tau }_2~~\underline{\beta }_2\underline{\beta }_1}\gamma _{\quad
\underline{\tau }_1\underline{\tau }_2\underline{\beta }_2}^{\underline{%
\delta }_2}]  \label{2.86a}
\end{eqnarray}
and
\begin{eqnarray}
& &\Phi _{<\gamma ><\alpha >} =\Phi _{\underline{\gamma }_1\underline{\gamma
}_2\underline{\alpha }_1\underline{\alpha }_2}=\frac
1{2(n+m_1+...+m_z)}\varepsilon _{\underline{\gamma }_1\underline{\alpha }%
_1}\varepsilon _{\underline{\gamma }_2\underline{\alpha }_2}[X_{\quad ~%
\underline{~\beta }_1\quad \underline{\mu }_1\underline{\beta }_2}^{%
\underline{\beta }_1\underline{\mu }_1~~\underline{\beta }_2}+  \notag \\
&&T_{\qquad ~~\underline{\beta }_2\underline{\mu }_1}^{\underline{\tau }_1%
\underline{\tau }_2\underline{\beta }_1\quad ~\underline{\beta }_2}~\gamma
_{\quad \underline{\tau }_1\underline{\tau }_2\underline{\beta }_1}^{%
\underline{\mu }_1}+X_{\qquad \quad \underline{\beta }_2\underline{\mu }_2%
\underline{\nu }_1}^{\underline{\beta }_2\underline{\mu }_2\underline{\nu }%
_1}+T_{\qquad \underline{\beta }_1\quad ~\underline{\delta }_2~}^{\underline{%
\tau }_1\underline{\tau }_2~~\underline{\beta }_2\underline{\beta }_1}\gamma
_{\quad \underline{\tau }_1\underline{\tau }_2\underline{\beta }_2}^{%
\underline{\delta }_2}]-  \notag \\
&&\frac 12[X_{\underline{\gamma }_1~\underline{\alpha }_1\underline{\alpha }%
_2\underline{\delta }_1\underline{\gamma }_2}^{\quad \underline{\delta }%
_1}+T_{\qquad \underline{\alpha }_1\underline{\alpha }_2\underline{\delta }_1%
\underline{\gamma }_2}^{\underline{\tau }_1\underline{\tau }_2}~\gamma
_{\quad \underline{\tau }_1\underline{\tau }_2\underline{\gamma }_1}^{%
\underline{\delta }_1}+  \notag \\
&&X_{\underline{\gamma }_2~\underline{\alpha }_1\underline{\alpha }_2%
\underline{\delta }_1\underline{\gamma }_2}^{\quad \underline{\delta }%
_2}+T_{\qquad \underline{\alpha }_1\underline{\alpha }_2\underline{\gamma }_1%
\underline{\delta }_2~}^{\underline{\tau }_1\underline{\tau }_2}\gamma
_{\quad \underline{\tau }_1\underline{\tau }_2\underline{\gamma }_2}^{%
\underline{\delta }_2}].  \label{2.87a}
\end{eqnarray}

The components of the conformal Weyl d-spinor can be computed by putting
d-spinor values of the curvature (\ref{2.83a}) and Ricci (\ref{2.84a})
d-tensors into corresponding expression for the d-tensor (\ref{2.33}). We
omit this calculus in this work.

\chapter{ Ha-Spinors and Field Interactions}

The problem of formulation gravitational and gauge field equations on
different types of locally anisotropic spaces is considered, for instance,
in \cite{ma94,bej,asa88} and \cite{vg}. In this Chapter we shall introduce
the basic field equations for gravitational and matter field la-interactions
in a generalized form for generic higher order anisotropic spaces.

\section{Scalar field ha--interactions}

Let $\varphi \left( u\right) =(\varphi _1\left( u\right) ,\varphi _2\left(
u\right) \dot{,}...,\varphi _k\left( u\right) )$ be a complex k-component
scalar field of mass $\mu $ on ha-space $\mathcal{E}^{<z>}.$ The d-covariant
generalization of the conformally invariant (in the massless case) scalar
field equation \cite{penr1,penr2} can be defined by using the d'Alambert
locally anisotropic operator \cite{ana94,vst96} $\Box =D^{<\alpha
>}D_{<\alpha >}$, where $D_{<\alpha >}$ is a d-covariant derivation on $%
\mathcal{E}^{<z>}$ and constructed, for simplicity, by using Christoffel
d--symbols (all formulas for field equations and conservation values can be
deformed by using corresponding deformations d--tensors $P_{<\beta ><\gamma
>}^{<\alpha >}$ from the Cristoffel d--symbols, or the canonical
d--connection to a general d-connection into consideration):

\begin{equation}
(\Box +\frac{n_{E}-2}{4(n_{E}-1)}\overleftarrow{R}+\mu ^{2})\varphi \left(
u\right) =0,  \label{2.88a}
\end{equation}%
where $n_{E}=n+m_{1}+...+m_{z}.$We must change d-covariant derivation $%
D_{<\alpha >}$ into $^{\diamond }D_{<\alpha >}=D_{<\alpha >}+ieA_{<\alpha >}$
and take into account the d-vector current
\begin{equation*}
J_{<\alpha >}^{(0)}\left( u\right) =i(\left( \overline{\varphi }\left(
u\right) D_{<\alpha >}\varphi \left( u\right) -D_{<\alpha >}\overline{%
\varphi }\left( u\right) )\varphi \left( u\right) \right)
\end{equation*}%
if interactions between locally anisotropic electromagnetic field ( d-vector
potential $A_{<\alpha >}$ ), where $e$ is the electromagnetic constant, and
charged scalar field $\varphi $ are considered. The equations (\ref{2.88a})
are (locally adapted to the N-connection structure) Euler equations for the
Lagrangian
\begin{eqnarray}
&&\mathcal{L}^{(0)}\left( u\right) =\sqrt{|g|}\left[ g^{<\alpha ><\beta
>}\delta _{<\alpha >}\overline{\varphi }\left( u\right) \delta _{<\beta
>}\varphi \left( u\right) -\left( \mu ^{2}+\frac{n_{E}-2}{4(n_{E}-1)}\right)
\overline{\varphi }\left( u\right) \varphi \left( u\right) \right] ,
\label{2.89a} \\
&&  \notag
\end{eqnarray}%
where $|g|=detg_{<\alpha ><\beta >}.$

The locally adapted variations of the action with Lagrangian (\ref{2.89a})
on variables $\varphi \left( u\right) $ and $\overline{\varphi }\left(
u\right) $ leads to the locally anisotropic generalization of the
energy-momentum tensor,
\begin{eqnarray}
E_{<\alpha ><\beta >}^{(0,can)}\left( u\right) &=&\delta _{<\alpha >}%
\overline{\varphi }\left( u\right) \delta _{<\beta >}\varphi \left( u\right)
+\delta _{<\beta >}\overline{\varphi }\left( u\right) \delta _{<\alpha
>}\varphi \left( u\right)  \label{2.90a} \\
&&-\frac{1}{\sqrt{|g|}}g_{<\alpha ><\beta >}\mathcal{L}^{(0)}\left( u\right)
,  \notag
\end{eqnarray}%
and a similar variation on the components of a d-metric (\ref{dmetrichcv})
leads to a symmetric metric energy-momentum d-tensor,
\begin{eqnarray}
&&E_{<\alpha ><\beta >}^{(0)}\left( u\right) =E_{(<\alpha ><\beta
>)}^{(0,can)}\left( u\right)  \label{2.91a} \\
&&-\frac{n_{E}-2}{2(n_{E}-1)}\left[ R_{(<\alpha ><\beta >)}+D_{(<\alpha
>}D_{<\beta >)}-g_{<\alpha ><\beta >}\Box \right] \overline{\varphi }\left(
u\right) \varphi \left( u\right) .  \notag
\end{eqnarray}%
Here we note that we can obtain a nonsymmetric energy-momentum d-tensor if
we firstly vary on $G_{<\alpha ><\beta >}$ and than impose the constraint of
compatibility with the N-connection structure. We also conclude that the
existence of a N-connection in dv-bundle $\mathcal{E}^{<z>}$ results in a
nonequivalence of energy-momentum d-tensors (\ref{2.90a}) and (\ref{2.91a}),
nonsymmetry of the Ricci tensor, nonvanishing of the d-covariant derivation
of the Einstein d-tensor, $D_{<\alpha >}\overleftarrow{G}^{<\alpha ><\beta
>}\neq 0$ and, in consequence, a corresponding breaking of conservation laws
on higher order anisotropic spaces when $D_{<\alpha >}E^{<\alpha ><\beta
>}\neq 0\,$. The problem of formulation of conservation laws on locally
anisotropic spaces is discussed in details and two variants of its solution
(by using nearly autoparallel maps and tensor integral formalism on locally
anisotropic and higher order multispaces) are proposed in \cite{vst96}.

In this Chapter we present only straightforward generalizations of field
equations and necessary formulas for energy-momentum d-tensors of matter
fields on $\mathcal{E}^{<z>}$ considering that it is naturally that the
conservation laws (usually being consequences of global, local and/or
intrinsic symmetries of the fundamental space-time and of the type of field
interactions) have to be broken on locally anisotropic spaces.

\section{ Proca equations on ha--spaces}

\index{Proca}

Let consider a d-vector field $\varphi _{<\alpha >}\left( u\right) $ with
mass $\mu ^2$ (locally anisotropic Proca field ) interacting with exterior
la-gravitational field. From the Lagrangian
\begin{equation}
\mathcal{L}^{(1)}\left( u\right) =%
\sqrt{\left| g\right| }\left[ -\frac 12{\overline{f}}_{<\alpha ><\beta
>}\left( u\right) f^{<\alpha ><\beta >}\left( u\right) +\mu ^2{\overline{%
\varphi }}_{<\alpha >}\left( u\right) \varphi ^{<\alpha >}\left( u\right) %
\right] ,  \label{2.92a}
\end{equation}
where $f_{<\alpha ><\beta >}=D_{<\alpha >}\varphi _{<\beta >}-D_{<\beta
>}\varphi _{<\alpha >},$ one follows the Proca equations on higher order
anisotropic spaces
\begin{equation}
D_{<\alpha >}f^{<\alpha ><\beta >}\left( u\right) +\mu ^2\varphi ^{<\beta
>}\left( u\right) =0.  \label{2.93a}
\end{equation}
Equations (\ref{2.93a}) are a first type constraints for $\beta =0.$ Acting
with $D_{<\alpha >}$ on (\ref{2.93a}), for $\mu \neq 0$ we obtain second
type constraints
\begin{equation}
D_{<\alpha >}\varphi ^{<\alpha >}\left( u\right) =0.  \label{2.94a}
\end{equation}

Putting (\ref{2.94a}) into (\ref{2.93a}) we obtain second order field
equations with respect to $\varphi _{<\alpha >}$ :
\begin{equation}
\Box \varphi _{<\alpha >}\left( u\right) +R_{<\alpha ><\beta >}\varphi
^{<\beta >}\left( u\right) +\mu ^2\varphi _{<\alpha >}\left( u\right) =0.
\label{2.95a}
\end{equation}
The energy-momentum d-tensor and d-vector current following from the (\ref%
{2.95a}) can be written as
\begin{eqnarray*}
E_{<\alpha ><\beta >}^{(1)}\left( u\right) &=&-g^{<\varepsilon ><\tau
>}\left( {\overline{f}}_{<\beta ><\tau >}f_{<\alpha ><\varepsilon >}+{%
\overline{f}}_{<\alpha ><\varepsilon >}f_{<\beta ><\tau >}\right) \\
&&+\mu ^2\left( {\overline{\varphi }}_{<\alpha >}\varphi _{<\beta >}+{%
\overline{\varphi }}_{<\beta >}\varphi _{<\alpha >}\right) -\frac{g_{<\alpha
><\beta >}}{\sqrt{\left| g\right| }}\mathcal{L}^{(1)}\left( u\right)
\end{eqnarray*}
and
\begin{equation*}
J_{<\alpha >}^{\left( 1\right) }\left( u\right) =i\left( {\overline{f}}%
_{<\alpha ><\beta >}\left( u\right) \varphi ^{<\beta >}\left( u\right) -{%
\overline{\varphi }}^{<\beta >}\left( u\right) f_{<\alpha ><\beta >}\left(
u\right) \right) .
\end{equation*}

For $\mu =0$ the d-tensor $f_{<\alpha ><\beta >}$ and the Lagrangian (\ref%
{2.92a}) are invariant with respect to locally anisotropic gauge transforms
of type
\begin{equation*}
\varphi _{<\alpha >}\left( u\right) \rightarrow \varphi _{<\alpha >}\left(
u\right) +\delta _{<\alpha >}\Lambda \left( u\right) ,
\end{equation*}
where $\Lambda \left( u\right) $ is a d-differentiable scalar function, and
we obtain a locally anisot\-rop\-ic variant of Maxwell theory.

\section{ Higher order anisotropic Dirac equations}

Let denote the Dirac d--spinor field on $\mathcal{E}^{<z>}$ as $\psi \left(
u\right) =\left( \psi ^{\underline{\alpha }}\left( u\right) \right) $ and
consider as the generalized Lorentz transforms the group of automorphysm of
the metric $G_{<\widehat{\alpha }><\widehat{\beta }>}$ (see (\ref{dmetrichcv}%
)).The d--covariant derivation of field $\psi \left( u\right) $ is written
as
\begin{equation}
\overrightarrow{\bigtriangledown _{<\alpha >}}\psi =\left[ \delta _{<\alpha
>}+\frac{1}{4}C_{\widehat{\alpha }\widehat{\beta }\widehat{\gamma }}\left(
u\right) ~l_{<\alpha >}^{\widehat{\alpha }}\left( u\right) \sigma ^{\widehat{%
\beta }}\sigma ^{\widehat{\gamma }}\right] \psi ,  \label{2.96a}
\end{equation}%
where coefficients $C_{\widehat{\alpha }\widehat{\beta }\widehat{\gamma }%
}=\left( D_{<\gamma >}l_{\widehat{\alpha }}^{<\alpha >}\right) l_{\widehat{%
\beta }<\alpha >}l_{\widehat{\gamma }}^{<\gamma >}$ generalize for ha-spaces
the corresponding Ricci coefficients on Riemannian spaces \cite{foc}. Using $%
\sigma $-objects $\sigma ^{<\alpha >}\left( u\right) $ (see (\ref{2.44a})
and (\ref{2.60a})--(\ref{2.62a})) we define the Dirac equations on
ha--spaces:
\begin{equation*}
(i\sigma ^{<\alpha >}\left( u\right) \overrightarrow{\bigtriangledown
_{<\alpha >}}-\mu )\psi =0,
\end{equation*}%
which are the Euler equations for the Lagrangian
\begin{eqnarray}
\mathcal{L}^{(1/2)}\left( u\right) &=&\sqrt{\left| g\right| }\{[\psi
^{+}\left( u\right) \sigma ^{<\alpha >}\left( u\right) \overrightarrow{%
\bigtriangledown _{<\alpha >}}\psi \left( u\right) \notag \\ & & -(\overrightarrow{%
\bigtriangledown _{<\alpha >}}\psi ^{+}\left( u\right) )\sigma ^{<\alpha
>}\left( u\right) \psi \left( u\right) ]
-\mu \psi ^{+}\left( u\right) \psi \left( u\right) \},
\label{2.97a}
\end{eqnarray}%
where $\psi ^{+}\left( u\right) $ is the complex conjugation and
transposition of the column $\psi \left( u\right) .$

From (\ref{2.97a}) we obtain the d--metric energy-momentum d-tensor
\begin{eqnarray*}
& & E_{<\alpha ><\beta >}^{(1/2)}= \frac i4[\psi ^{+}\left( u\right) \sigma
_{<\alpha >}\left( u\right) \overrightarrow{\bigtriangledown _{<\beta >}}%
\psi \left( u\right) +\psi ^{+}\left( u\right) \sigma _{<\beta >}\left(
u\right) \overrightarrow{\bigtriangledown _{<\alpha >}}\psi \left( u\right)
\\
& & -(\overrightarrow{\bigtriangledown _{<\alpha >}}\psi ^{+}\left( u\right)
)\sigma _{<\beta >}\left( u\right) \psi \left( u\right) -(\overrightarrow{%
\bigtriangledown _{<\beta >}}\psi ^{+}\left( u\right) )\sigma _{<\alpha
>}\left( u\right) \psi \left( u\right) ]
\end{eqnarray*}
and the d-vector source
\begin{equation*}
J_{<\alpha >}^{(1/2)}\left( u\right) =\psi ^{+}\left( u\right) \sigma
_{<\alpha >}\left( u\right) \psi \left( u\right) .
\end{equation*}
We emphasize that locally anisotropic interactions with exterior gauge
fields can be introduced by changing the higher order anisotropic partial
derivation from (\ref{2.96a}) in this manner:
\begin{equation*}
\delta _\alpha \rightarrow \delta _\alpha +ie^{\star }B_\alpha ,
\end{equation*}
where $e^{\star }$ and $B_\alpha $ are respectively the constant d-vector
potential of locally anisotropic gauge interactions on higher order
anisotropic spaces (see \cite{vg} and the next section).

\section{ D--spinor Yang--Mills fields}

\index{Yang--Mills}

We consider a dv--bundle $\mathcal{B}_E,~\pi _B:\mathcal{B\rightarrow E}%
^{<z>}$ on ha--space $\mathcal{E}^{<z>}.\mathcal{\,}$Additionally to
d-tensor and d-spinor indices we shall use capital Greek letters, $\Phi
,\Upsilon ,$ $\Xi ,\Psi ,...$ for fibre (of this bundle) indices (see
details in \cite{penr1,penr2} for the case when the base space of the
v-bundle $\pi _B$ is a locally isotropic space-time). Let $%
\underline{\bigtriangledown }_{<\alpha >}$ be, for simplicity, a
torsionless, linear connection in $\mathcal{B}_E$ satisfying conditions:
\begin{eqnarray*}
& &\underline{\bigtriangledown }_{<\alpha >} :{\ \Upsilon }^\Theta
\rightarrow {\ \Upsilon }_{<\alpha >}^\Theta \quad \left[ \mbox{or }{\ \Xi }%
^\Theta \rightarrow {\ \Xi }_{<\alpha >}^\Theta \right] , \\
& & \underline{\bigtriangledown }_{<\alpha >}\left( \lambda ^\Theta +\nu
^\Theta \right) =\underline{\bigtriangledown }_{<\alpha >}\lambda ^\Theta +%
\underline{\bigtriangledown }_{<\alpha >}\nu ^\Theta , \\
& & \underline{\bigtriangledown }_{<\alpha >}~(f\lambda ^\Theta )=\lambda
^\Theta \underline{\bigtriangledown }_{<\alpha >}f+f\underline{%
\bigtriangledown }_{<\alpha >}\lambda ^\Theta ,\quad f\in {\ \Upsilon }%
^\Theta ~[\mbox{or }{\ \Xi }^\Theta ],
\end{eqnarray*}
where by ${\ \Upsilon }^\Theta ~\left( {\ \Xi }^\Theta \right) $ we denote
the module of sections of the real (complex) v--bundle $\mathcal{B}_E$
provided with the abstract index $\Theta .$ The curvature of connection $%
\underline{\bigtriangledown }_{<\alpha >}$ is defined as
\begin{equation*}
K_{<\alpha ><\beta >\Omega }^{\qquad \Theta }\lambda ^\Omega =\left(
\underline{\bigtriangledown }_{<\alpha >}\underline{\bigtriangledown }%
_{<\beta >}-\underline{\bigtriangledown }_{<\beta >}\underline{%
\bigtriangledown }_{<\alpha >}\right) \lambda ^\Theta .
\end{equation*}

For Yang-Mills fields as a rule one considers that $\mathcal{B}_E$ is
enabled with a unitary (complex) structure (complex conjugation changes
mutually the upper and lower Greek indices). It is useful to introduce
instead of $K_{<\alpha ><\beta >\Omega }^{\qquad \Theta }$ a Hermitian
matrix $F_{<\alpha ><\beta >\Omega }^{\qquad \Theta }=i$ $K_{<\alpha ><\beta
>\Omega }^{\qquad \Theta }$ connected with components of the Yang-Mills
d-vector potential $B_{<\alpha >\Xi }^{\quad \Phi }$ according the formula:

\begin{equation}
\frac 12F_{<\alpha ><\beta >\Xi }^{\qquad \Phi }=\underline{\bigtriangledown
}_{[<\alpha >}B_{<\beta >]\Xi }^{\quad \Phi }-iB_{[<\alpha >|\Lambda
|}^{\quad \Phi }B_{<\beta >]\Xi }^{\quad \Lambda },  \label{2.98a}
\end{equation}
where the locally anisotropic space indices commute with capital Greek
indices. The gauge transforms are written in the form:

\begin{eqnarray*}
B_{<\alpha >\Theta }^{\quad \Phi } &\mapsto &B_{<\alpha >\widehat{\Theta }%
}^{\quad \widehat{\Phi }}=B_{<\alpha >\Theta }^{\quad \Phi }~s_\Phi ^{\quad
\widehat{\Phi }}~q_{\widehat{\Theta }}^{\quad \Theta }+is_\Theta ^{\quad
\widehat{\Phi }}\underline{\bigtriangledown }_{<\alpha >}~q_{\widehat{\Theta
}}^{\quad \Theta }, \\
F_{<\alpha ><\beta >\Xi }^{\qquad \Phi } &\mapsto &F_{<\alpha ><\beta >%
\widehat{\Xi }}^{\qquad \widehat{\Phi }}=F_{<\alpha ><\beta >\Xi }^{\qquad
\Phi }s_\Phi ^{\quad \widehat{\Phi }}q_{\widehat{\Xi }}^{\quad \Xi },
\end{eqnarray*}
where matrices $s_\Phi ^{\quad \widehat{\Phi }}$ and $q_{\widehat{\Xi }%
}^{\quad \Xi }$ are mutually inverse (Hermitian conjugated in the unitary
case). The Yang-Mills equations on torsionless locally anisotropic spaces %
\cite{vg} (see details in the next Section) are written in this form:
\begin{eqnarray}
\underline{\bigtriangledown }^{<\alpha >}F_{<\alpha ><\beta >\Theta
}^{\qquad \Psi } &=&J_{<\beta >\ \Theta }^{\qquad \Psi },  \label{2.99a} \\
\underline{\bigtriangledown }_{[<\alpha >}F_{<\beta ><\gamma >]\Theta
}^{\qquad \Xi } &=&0.  \notag
\end{eqnarray}
We must introduce deformations of connection of type $\underline{%
\bigtriangledown }_\alpha ^{\star }~\longrightarrow \underline{%
\bigtriangledown }_\alpha +P_\alpha ,$ (the deformation d-tensor $P_\alpha $
is induced by the torsion in dv-bundle $\mathcal{B}_E)$ into the definition
of the curvature of gauge ha--fields (\ref{2.98a}) and motion equations (\ref%
{2.99a}) if interactions are modeled on a generic higher order anisotropic
space.

\section{D--spinor Einstein--Cartan Theory}

\index{Einstein--Cartan}

The Einstein equations in some models of higher order anisotropic
supergravity have been considered in \cite{vlasg,vbook}. Here we note that
the Einstein equations and conservation laws on v--bundles provided with
N-connection structures were studied in detail in \cite%
{ma87,ma94,ana86,ana87,vodg,voa,vcl96}. In Ref. \cite{vg} we proved that the
locally anisotropic gravity can be formulated in a gauge like manner and
analyzed the conditions when the Einstein gravitational locally anisotropic
field equations are equivalent to a corresponding form of Yang-Mills
equations. Our aim here is to write the higher order anisotropic
gravitational field equations in a form more convenient for theirs
equivalent reformulation in higher order anisotropic spinor variables.

\subsection{Einstein ha--equations}

We define d-tensor $\Phi _{<\alpha ><\beta >}$ as to satisfy conditions
\begin{equation*}
-2\Phi _{<\alpha ><\beta >}\doteq R_{<\alpha ><\beta >}-\frac
1{n+m_1+...+m_z}%
\overleftarrow{R}g_{<\alpha ><\beta >}
\end{equation*}
which is the torsionless part of the Ricci tensor for locally isotropic
spaces \cite{penr1,penr2}, i.e. $\Phi _{<\alpha >}^{~~<\alpha >}\doteq 0$.\
The Einstein equations on higher order anisotrop\-ic spaces
\begin{equation}
\overleftarrow{G}_{<\alpha ><\beta >}+\lambda g_{<\alpha ><\beta >}=\kappa
E_{<\alpha ><\beta >},  \label{2.34a}
\end{equation}
where
\begin{equation*}
\overleftarrow{G}_{<\alpha ><\beta >}=R_{<\alpha ><\beta >}-\frac 12%
\overleftarrow{R}g_{<\alpha ><\beta >}
\end{equation*}
is the Einstein d--tensor, $\lambda $ and $\kappa $ are correspondingly the
cosmological and gravitational constants and by $E_{<\alpha ><\beta >}$ is
denoted the locally anisotropic energy--momentum d--tensor, can be rewritten
in equivalent form:
\begin{equation}
\Phi _{<\alpha ><\beta >}=-\frac \kappa 2(E_{<\alpha ><\beta >}-\frac
1{n+m_1+...+m_z}E_{<\tau >}^{~<\tau >}~g_{<\alpha ><\beta >}).  \label{2.35a}
\end{equation}

Because ha--spaces generally have nonzero torsions we shall add to (\ref%
{2.35a}) (equivalently to (\ref{2.34a})) a system of algebraic d--field
equations with the source $S_{~<\beta ><\gamma >}^{<\alpha >}$ being the
locally anisotropic spin density of matter (if we consider a variant of
higher order anisotropic Einstein--Cartan theory ):
\begin{equation}
T_{~<\alpha ><\beta >}^{<\gamma >}+2\delta _{~[<\alpha >}^{<\gamma
>}T_{~<\beta >]<\delta >}^{<\delta >}=\kappa S_{~<\alpha ><\beta
>.}^{<\gamma >}  \label{2.36a}
\end{equation}
From (\ref{2.36a}) one follows the conservation law of higher order
anisotropic spin matter:
\begin{equation*}
\bigtriangledown _{<\gamma >}S_{~<\alpha ><\beta >}^{<\gamma >}-T_{~<\delta
><\gamma >}^{<\delta >}S_{~<\alpha ><\beta >}^{<\gamma >}=E_{<\beta ><\alpha
>}-E_{<\alpha ><\beta >}.
\end{equation*}

\subsection{Einstein--Cartan d--equations}

Now we can write out the field equations of the Einstein--Cartan theory in
the d-spinor form. So, for the Einstein equations (\ref{2.34}) we have

\begin{equation*}
\overleftarrow{G}_{\underline{\gamma }_1\underline{\gamma }_2\underline{%
\alpha }_1\underline{\alpha }_2}+\lambda \varepsilon _{\underline{\gamma }_1%
\underline{\alpha }_1}\varepsilon _{\underline{\gamma }_2\underline{\alpha }%
_2}=\kappa E_{\underline{\gamma }_1\underline{\gamma }_2\underline{\alpha }_1%
\underline{\alpha }_2},
\end{equation*}
with $\overleftarrow{G}_{\underline{\gamma }_1\underline{\gamma }_2%
\underline{\alpha }_1\underline{\alpha }_2}$ from (\ref{2.86a}), or, by
using the d-tensor (\ref{2.87a}),

\begin{equation*}
\Phi _{\underline{\gamma }_1\underline{\gamma }_2\underline{\alpha }_1%
\underline{\alpha }_2}+(\frac{\overleftarrow{R}}4-\frac \lambda
2)\varepsilon _{\underline{\gamma }_1\underline{\alpha }_1}\varepsilon _{%
\underline{\gamma }_2\underline{\alpha }_2}=-\frac \kappa 2E_{\underline{%
\gamma }_1\underline{\gamma }_2\underline{\alpha }_1\underline{\alpha }_2},
\end{equation*}
which are the d-spinor equivalent of the equations (\ref{2.35a}). These
equations must be completed by the algebraic equations (\ref{2.36a}) for the
d-torsion and d-spin density with d-tensor indices changed into
corresponding d--spinor ones.

\subsection{Higher order an\-i\-sot\-rop\-ic gravitons}

Let a massless d-tensor field $h_{<\alpha ><\beta >}\left( u\right) $ is
interpreted as a small perturbation of the locally anisotropic background
metric d-field $g_{<\alpha ><\beta >}\left( u\right) .$ Considering, for
simplicity, a torsionless background we have locally anisotropic
Fierz--Pauli equations
\begin{equation*}
\Box h_{<\alpha ><\beta >}\left( u\right) +2R_{<\tau ><\alpha ><\beta ><\nu
>}\left( u\right) ~h^{<\tau ><\nu >}\left( u\right) =0
\end{equation*}
and d--gauge conditions
\begin{equation*}
D_{<\alpha >}h_{<\beta >}^{<\alpha >}\left( u\right) =0,\quad h\left(
u\right) \equiv h_{<\beta >}^{<\alpha >}(u)=0,
\end{equation*}
where $R_{<\tau ><\alpha ><\beta ><\nu >}\left( u\right) $ is curvature
d-tensor of the locally anisotropic background space (these formulae can be
obtained by using a perturbation formalism with respect to $h_{<\alpha
><\beta >}\left( u\right) $ developed in \cite{gri}; in our case we must
take into account the distinguishing of geometrical objects and operators on
higher order anisotropic spaces).

Finally, we remark that all presented geometric constructions contain those
elaborated for generalized Lagrange spaces \cite{ma87,ma94} (for which a
tangent bundle $TM$ is considered instead of a v-bundle $\mathcal{E}^{<z>}$
) and for constructions on the so called osculator bundles with different
prolongations and extensions of Finsler and Lagrange metrics \cite{mirata}.
We also note that the higher order Lagrange (Finsler) geometry is
characterized by a metric of type (dmetrichcv) with components parametized
as $g_{ij}=\frac 12\frac{\partial ^2\mathcal{L}}{\partial y^i\partial y^j}$ $%
\left( g_{ij}=\frac 12\frac{\partial ^2\Lambda ^2}{\partial y^i\partial y^j}%
\right) $ and $h_{a_pb_p}=g_{ij},$ where $\mathcal{L=L}$ $%
(x,y_{(1)},y_{(2)},....,y_{(z)})$ $\left( \Lambda =\Lambda \left(
x,y_{(1)},y_{(2)},....,y_{(z)}\right) \right) $ is a Lagrangian $\left( %
\mbox{Finsler metric}\right) $ on $TM^{(z)}$ (see details in \cite%
{ma87,ma94,mat,bej}). 

\end{document}